\documentclass[11pt]{book}

\usepackage[a4paper,twoside,hmarginratio=1:1]{geometry} 

\usepackage[english]{babel}
\usepackage[utf8]{inputenc}
\usepackage{floatrow} 

\usepackage{graphicx}
\usepackage{xcolor,colortbl}
\usepackage[font=small,labelfont=bf,hypcap=true]{caption} 
\usepackage[hypcap=true]{subcaption}
\usepackage[chapter]{algorithm}
\usepackage{algpseudocode}

\definecolor{mygreen}{cmyk}{0.82,0.11,1,0.25}
\definecolor{mydarkgreen}{cmyk}{0.82,0.11,1,0.45}
\definecolor{nominalgreen}{rgb}{0,0.5,0}
\definecolor{ncred}{rgb}{0.88,0,0}

\makeatletter
\newcounter{HALG@line}
\renewcommand{\theHALG@line}{\thealgorithm.\arabic{ALG@line}}
\makeatother

\newcommand{\plusline}{%
  \renewcommand{\ALG@lineno} {+%
    }%
}

\usepackage{setspace}

\usepackage{booktabs}
\usepackage{multirow}
\usepackage{arydshln}
\usepackage{amsfonts}
\usepackage{amsmath}
\usepackage{array}
\usepackage{bm}
\usepackage{mathtools}

\usepackage{enumerate}

\usepackage[square,numbers,sort&compress]{natbib}
\usepackage{doi}

\usepackage{pdfpages}
\usepackage[intoc]{nomencl}
\usepackage{etoolbox}

\makenomenclature{}

\usepackage[bf,topmarks,calcwidth,pagestyles]{titlesec}
\definecolor{gray75}{gray}{0.45}
\newcommand{\hsp}{\hspace{20pt}}
\titleformat{\chapter}
    [hang]
    {\Huge\bfseries}
    {\thechapter\hsp\textcolor{gray75}{$|$}\hsp}
    {0pt}
    {\Huge\bfseries\raggedright}
    [\color{gray75}\rule{\textwidth}{0.1cm}]

\usepackage{amsthm}
\usepackage{thmtools}
\declaretheoremstyle[numberwithin=chapter,thmbox=L]{Lthm}
\declaretheoremstyle[numberwithin=chapter,
                    shaded={rulecolor=mygreen!30, textwidth=0.95\linewidth, margin=0.8em, rulewidth=2pt, bgcolor=white!97!mygreen},
                    headformat={{\color{mydarkgreen} \NAME~\NUMBER~:\NOTE}\hfill\smallskip\linebreak},
                    headpunct={},
                    notebraces={}{},
                    notefont=\bfseries,
                ]{sthm}

\declaretheorem[style=sthm]{remark}
\declaretheorem[style=sthm]{hypothesis}
\declaretheorem[style=sthm]{notation}

\usepackage{fancyhdr}
\fancypagestyle{plain}{%
    \fancyhead{} 
}

\makeatletter
\def\cleardoublepage{\clearpage\if@twoside \ifodd\c@page\else
\begingroup \mbox{}
\thispagestyle{empty}
\newpage
\if@twocolumn\mbox{}\newpage\fi
\endgroup\fi\fi}
\makeatother

\usepackage[capitalize]{cleveref}

\crefformat{algorithm}{Alg.~#2#1#3}

\crefformat{line}{Step~#2#1#3}
\crefrangeformat{line}{Steps~#3#1#4 to~#5#2#6}
\crefmultiformat{line}{Steps~#2#1#3}{ and~#2#1#3}{, #2#1#3}{ and~#2#1#3}

\crefname{remark}{Remark}{Remarks}
\crefname{hypothesis}{Hypothesis}{Hypotheses}
\crefname{notation}{Notation}{Notations}

\usepackage{tikz}
\usepackage[export]{adjustbox}
\tikzstyle{every picture}+=[remember picture]
\usetikzlibrary{arrows.meta,positioning,fit,shapes.geometric,calc}
\pgfdeclarelayer{bg}    
\pgfsetlayers{bg,main}  

\tikzstyle{rect} = [rectangle, draw, text centered, minimum height=4em]
\tikzstyle{oval} = [ellipse, draw, text centered, minimum height=4em]
\tikzstyle{io} = [trapezium, trapezium left angle=95, trapezium right angle=85, draw, text width=4em, text centered, minimum height=3em, inner sep=0pt, minimum width=8em]
\tikzstyle{process} = [rectangle, draw, text width=8em, text centered, rounded corners, minimum height=3em]
\tikzstyle{arrow} = [draw, -{Latex[length=4mm,width=3mm]}]
\tikzstyle{smallarrow} = [draw, -{Latex[length=2mm,width=1mm]}]
\tikzstyle{dasharrow} = [draw,dashed, -{Latex[length=4mm,width=3mm]}]
\tikzstyle{decision} = [diamond, draw, text badly centered, inner sep=0pt]
\tikzset{-|-/.style={to path={(\tikztostart) -| ($(\tikztostart)!#1!(\tikztotarget)$) |- (\tikztotarget) \tikztonodes }  },-|-/.default=0.5,
	|-|/.style={to path={(\tikztostart) |- ($(\tikztostart)!#1!(\tikztotarget)$) -| (\tikztotarget)\tikztonodes} }, |-|/.default=0.5, }

\DeclareMathOperator*{\argmin}{arg\,min}
\DeclareMathOperator*{\argmax}{arg\,max}
\DeclareMathOperator*{\argminmax}{arg\,minmax}

\newcommand{\norm}[2]{{\left| \left| #1 \right| \right|}_{#2}}
\newcommand{\innerproduct}[2]{{\left\langle#1,#2\right\rangle}}

\newcommand{\PrimalSpace}[0]{\bm{\mathcal{V}}}
\newcommand{\PrimalRB}[0]{\mathbf{\Phi}}
\newcommand{\PrimalRBF}[0]{\bm{\phi}}

\newcommand{\DualSpace}[0]{\mathcal{W}_{+}}
\newcommand{\DualRB}[0]{\mathbf{\Theta}}
\newcommand{\DualRBF}[0]{\theta}

\newcommand{\DualDict}[0]{\mathbf{D}_d}
\newcommand{\DualDictActive}[0]{\mathbf{D}_{d,\mathcal{I}}}
\newcommand{\PrimalDict}[0]{\mathbf{D}_p}

\newcommand{\ActiveSet}[0]{\mathcal{A}}

\newcommand{\ndofprimal}[0]{N^h_u}
\newcommand{\ndofdual}[0]{N^h_\lambda}
\newcommand{\nrdofprimal}[0]{N^r_u}
\newcommand{\nrdofdual}[0]{N^r_\lambda}

\newcommand{\WarpedCoord}[0]{\zeta}
\newcommand{\WarpedSpace}[0]{\widetilde{\Gamma}}

\usepackage{hyperref}
\hypersetup{
    pdftitle={Low-rank and sparse approximations for contact mechanics},
    pdfauthor={Kiran Sagar Kollepara},
    pdfsubject={Doctoral thesis},
    pdfkeywords={Model Order Reduction, Low-rank methods, Contact Mechanics, Dictionary methods, Non-linear dimensionality reduction},
    bookmarks=true,
    bookmarksopen=true,
    bookmarksopenlevel=1,
    bookmarksnumbered=true,
    pdfpagemode=UseOutlines,
    colorlinks=true,
    linkcolor=blue,
    citecolor=ncred,
}

\begin{document}

\includepdf[pages={1}]{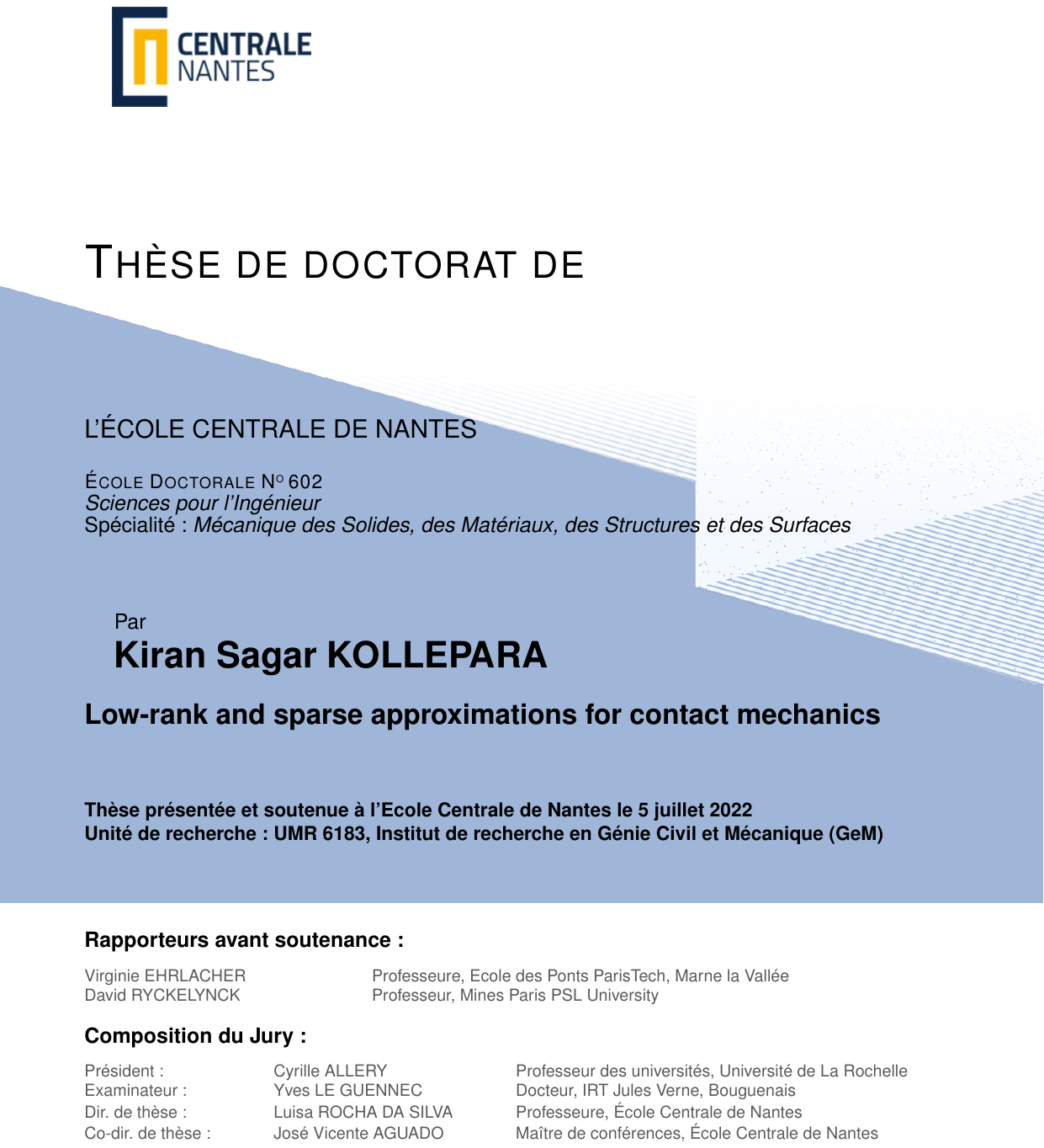}

\cleardoublepage{}
\includepdf[pages={1}]{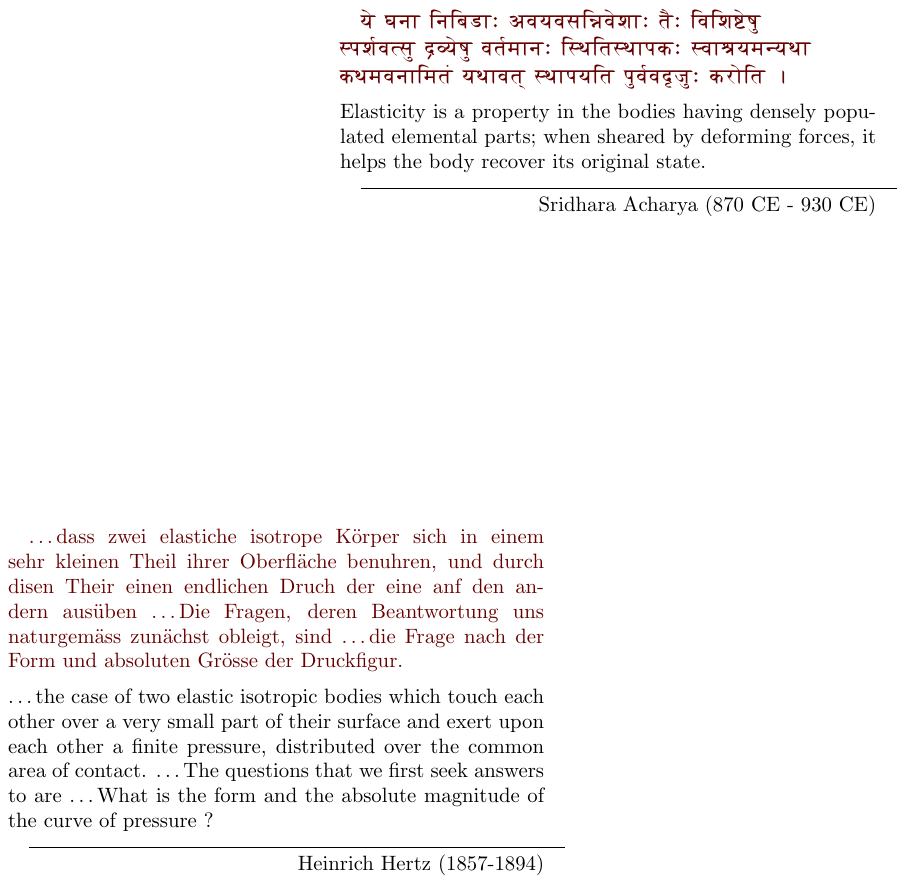}

\cleardoublepage{}
\includepdf[pages={1}]{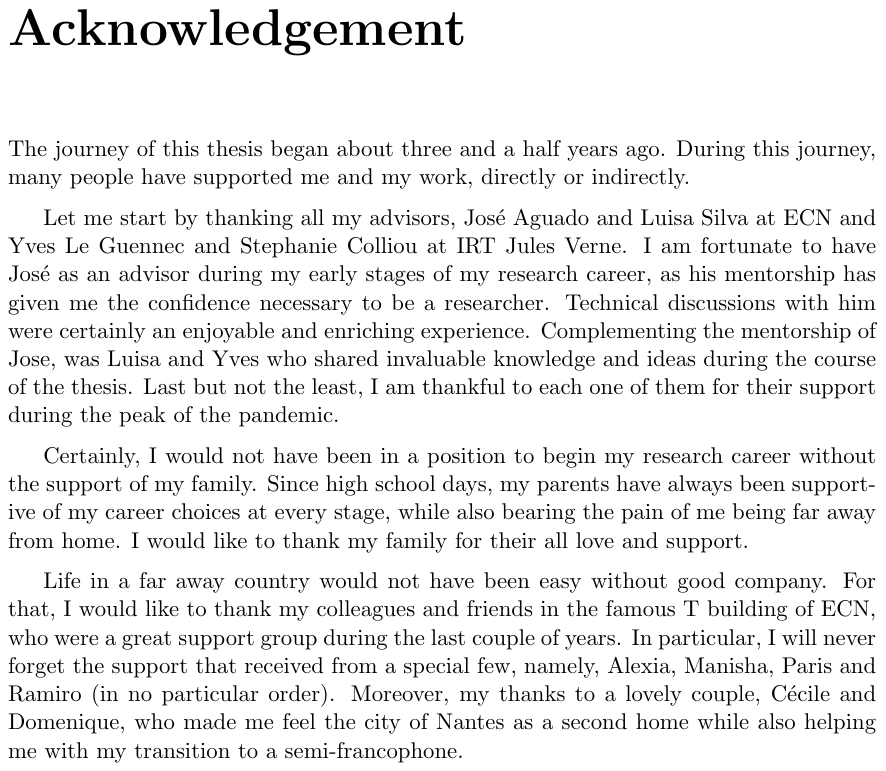}

\cleardoublepage{}
\frontmatter

\pdfbookmark{\contentsname}{toc}
\tableofcontents

\cleardoublepage{}
\phantomsection{}
\addcontentsline{toc}{chapter}{\listfigurename}
\listoffigures

\cleardoublepage{}
\phantomsection{}
\addcontentsline{toc}{chapter}{\listalgorithmname}
\listofalgorithms

\cleardoublepage{}
\phantomsection{}
\addcontentsline{toc}{chapter}{\listtablename}
\listoftables

\cleardoublepage{}
\phantomsection{}
\printnomenclature[6em]


\nomgroup{AOperators}
\nomenclature[A]{$k,\ g$}{Distance functions}
\nomenclature[A]{$b,\ d$}{Bilinear and Linear forms of distance function in weak form}
\nomenclature[A]{$a,\ f$}{Bilinear and Linear forms of energy functional}
\nomenclature[A]{$\bm{u},\ \bm{v}$}{Continuous/Discretized Displacement}
\nomenclature[A]{$\lambda, \eta, \bm{\lambda}, \bm{\eta}$}{Continuous/Discretized Lagrange Multiplier (Contact pressure)}
\nomenclature[A]{$\PrimalSpace,\ \DualSpace$}{Displacement solution space and Contact pressure solution cone}
\nomenclature[A]{$\mathcal{L}$}{Lagrangian functional}
\nomenclature[A]{$\mathcal{K}$}{Feasible region of an inequality constrained problem}

\nomgroup{CParametric}
\nomenclature[C]{$\mu$}{Parameter}
\nomenclature[C]{$\mathcal{P}$}{Parametric space}

\nomgroup{FEDiscretization}
\nomenclature[F]{$\langle\cdot\rangle^h$}{Discrete version of the entity $\langle\cdot\rangle$}
\nomenclature[F]{$\bm{N},\ M$}{Displacement and contact pressure shape functions}
\nomenclature[F]{$\mathbf{K},\ \bm{f}$}{Discretized operators of bilinear and linear terms $a$ and $f$}
\nomenclature[F]{$\mathbf{C},\ \bm{g}$}{Discretized operators of bilinear and linear terms $k$ and $g$}
\nomenclature[F]{$\bm{r}(\cdot)$}{Residual of a discretized problem}
\nomenclature[F]{$\mathbf{K}_{\texttt{mono}}$}{Operator to generate the residual with monolithic dictionary}
\nomenclature[F]{$\ndofprimal,\ \ndofdual$}{Number of FE dofs for $\bm{u}$ and $\lambda$}

\nomgroup{JContactRelated}
\nomenclature[J]{$(\bm{x}, \bar{\bm{x}})$}{Contact pairs}
\nomenclature[J]{$\Gamma$}{Potential contact surface}
\nomenclature[J]{$\Gamma_c, \ \Gamma_{/c}$}{Surface in contact, and out of contact (gap)}
\nomenclature[J]{$\ActiveSet,\ \mathcal{I}$}{Active sets}
\nomenclature[J]{$\ActiveSet^{\mathsf{c}},\ \mathcal{I}^{\mathsf{c}}$}{Complements of active sets}

\nomgroup{Reduced}
\nomenclature[R]{$\langle\cdot\rangle^r$}{Reduced space of $\langle\cdot\rangle$}
\nomenclature[R]{$\PrimalRB,\ \DualRB$}{Primal and Dual reduced bases}
\nomenclature[R]{$\PrimalRBF,\ \DualRBF$}{Primal and Dual reduced bases vectors/functions}
\nomenclature[R]{$\widehat{\langle\cdot\rangle}$}{Reduced operators and dofs}
\nomenclature[R]{$\Pi_\circ$}{Projection operator on the subspace/subcone of basis $\circ$}
\nomenclature[R]{$\mathbf{S}_\mathtt{tr}$}{Training set snapshots}
\nomenclature[R]{$\mathbf{D}$}{Dictionary, Monolithic Dictionary}
\nomenclature[R]{$\PrimalDict,\ \DualDict$}{Primal and Dual Dictionaries}
\nomenclature[R]{$\mathbf{B}$}{An orthogonal basis of a given snapshot matrix or its randomized combination}
\nomenclature[R]{$\mathbf{R}$}{Matrix containing random entries between $[0,1]$}
\nomenclature[R]{$\nrdofprimal,\ \nrdofdual$}{Number of reduced dofs for $\bm{u}$ and $\lambda$}
\nomenclature[R]{$\bm{d}, \widetilde{\bm{d}}$}{Element in a dictionary, Leave-one-out candidate in a dictionary}

\nomgroup{TErrors}
\nomenclature[T]{$H(n)$}{Nested error between $n$-th and $n+1$-th level nested training set}
\nomenclature[T]{$\varepsilon_{\texttt{CrPen}}$}{Cross penetration error in convexity test of the feasible region}
\nomenclature[T]{$\varepsilon_{\texttt{CHLS}}$}{Convex hull least square error for the test of convex subset hypothesis}
\nomenclature[T]{$C(m)$}{Compactness of rank-$m$}
\nomenclature[T]{$G(m)$}{Generalization ability of rank-$m$ basis}
\nomenclature[T]{$S(m)$}{Specificity of rank-$m$ basis}

\nomgroup{VRegression}
\nomenclature[V]{$\bm{\alpha},\ \bm{\gamma}$}{Coefficient vector multiplying with columns of a matrix/basis/dictionary}
\nomenclature[V]{$\bm{x}$}{Vector/Signal to be approximated}

\nomgroup{XNonlinearInterpolation}
\nomenclature[X]{$\varphi$}{Non-linear transformation to a higher-dimensional warped space/domain}
\nomenclature[X]{$I_k$}{DTW alignments computed for a vector $\bm{\lambda}_k$}
\nomenclature[X]{$\WarpedCoord$}{Co-ordinate in the warped space}
\nomenclature[X]{$\widetilde{\langle\cdot\rangle}$}{Quantity $\langle\cdot\rangle$ in warped space}
\nomenclature[X]{$\circ_{\texttt{int}}$}{Quantity $\circ$ interpolated in warped space}
\nomenclature[X]{$\mathbf{D}_p^{\texttt{adapt}}$}{Dictionary adapted/enriched using non-linear interpolations}
\nomenclature[X]{$\alpha$}{Parameter for interpolating in high-dimensional space}

\nomgroup{YOtherSymbols}
\nomenclature[Y]{$\circ[:,S],\ \circ[S]$}{Slicing of a vector/matrix $\circ$ using indices in set $S$ (Python style)}
\nomenclature[Y]{$\delta$}{Truncation tolerance for an orthogonal basis}
\nomenclature[Y]{$\varepsilon$}{Machine precision, Tolerance on residual}
\nomenclature[Y]{${(z)}^{+,\tau}$}{$z$ if $z>\tau$, $0$ otherwise}
\nomenclature[Y]{${(z)}^-$}{$z$ if $z<0$, $0$ otherwise}
\nomenclature[Y]{${\langle\cdot\rangle}^\dagger$}{Moore-Penrose pseudo-inverse}
\nomenclature[Y]{$\bm{0},\ \bm{1}$}{Matrix or Vector consisting of zeros or ones of appropriate size}


\mainmatter{}
\chapter*{Abstract}
\markboth{Abstract}{Abstract}
\addcontentsline{toc}{chapter}{Abstract}
Non-conformance decision-making processes in high-precision manufacturing of engineering structures are often delayed due to numerical simulations that are needed for analyzing the defective  parts and assemblies. Various engineering assemblies often involve interfaces between parts that can only be simulated using the modeling of contact. Therefore, efficient parametric reduced order models (ROMs) are necessary for performing contact mechanics simulations in near real-time scenarios. 

Typical strategies for reducing the computational cost of contact mechanics models use low-rank approximations. The underlying hypothesis is the existence of a low-dimensional subspace for the displacement field and a low-dimensional subcone for the contact pressure, as a result of non-negativity constraints. However, the contact pressure exhibits a local nature, as the position of contact can vary with parameters like loading or geometry. In this thesis, the adequacy of low-rank approximations for contact mechanics is investigated and alternative routes based on sparse regression techniques are explored.

It is shown that the local nature leads to loss of linear separability of contact pressure, thereby limiting the reconstruction accuracy of low-rank methods. The applicability of the low-rank assumption to contact pressure is analyzed using three different criteria, namely compactness, generalization and specificity.

Subsequently, the use of over-complete dictionaries containing a large number of snapshots to mitigate the inseparability issues is investigated. Two strategies to solve the dictionary-based approximation problem are devised: one based on a greedy active-set method where the elements from the contact pressure dictionary are selected greedily and another approach based on convex hull approximations eliminating the need to explicitly enforce non-penetration constraints in convex problems.

Lastly, Dynamic Time Warping (DTW) is studied as a possible non-linear interpolation method that permits the exploration of the non-linear manifold. This allows the synthesis of snapshots not computed in the training set with low complexity, thereby reducing the burden of creating over-complete dictionaries in the offline phase.

\chapter*{Introduction}
\markboth{Introduction}{Introduction}
\addcontentsline{toc}{chapter}{Introduction}

\let\oldfig\thefigure
\renewcommand{\thefigure}{\MakeUppercase{\roman{figure}}}

In the manufacturing of engineering structures, process variabilities often lead to deviations from the tolerance specifications on the features of a part or an assembly. These deviations can influence the stresses at the contact interfaces in the assembly, potentially leading to reduced performance or even failure of the assembled structure. Economic constraints often provide the motivation to salvage the ``non-conforming'' hardware. Depending on the nature of the deviations, a non-conforming part or assembly can be either: \textit{accepted}, \textit{reworked} to improve compliance or \textit{rejected}. The decision is usually based on further analyses or tests on the non-conforming part. Naturally, our understanding of the mechanics of the contact interfaces in assemblies is crucial in this decision-making process.

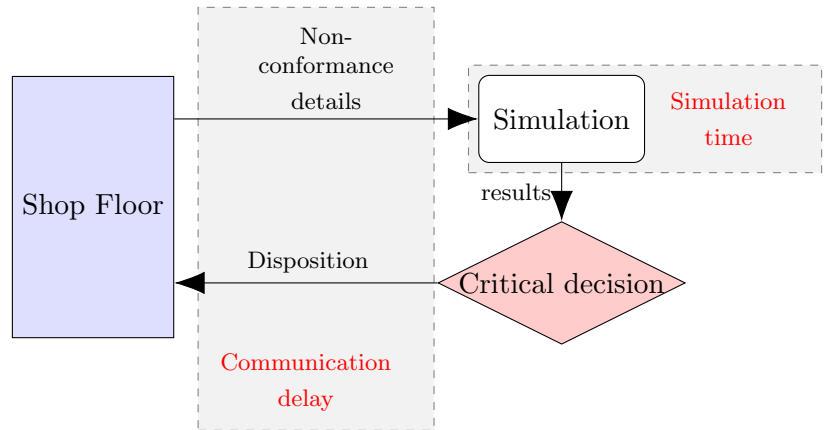
\begin{figure}[htpb!]
	\centering
	\begin{tikzpicture}[node distance=1cm, auto]
		\node [process, text width=5em, fill=white] (simulation) {Simulation};
		\node [rect, minimum height=9em, fill=blue!13,below left=-3em and 4cm of simulation] (shop) {Shop Floor};

		\node[decision, aspect=2, below= 2em of simulation, fill=red!20] (decide){Critical decision} ;
		\path [arrow] (simulation.south) to node [midway,left]() {\footnotesize results} (decide.north);
		\path [arrow] (shop.east |- simulation.west) to node [midway, above, align=center, text width=6em](ncdetails) {\footnotesize Non-conformance \\ details} (simulation.west);
		\path [arrow] (decide.west) to node [midway, above, align=center, text width=4em](disposition) {\footnotesize Disposition } (decide.west -| shop.east);
		\node[below=0.8cm of disposition, text width=6em, align=center] (CT) {\footnotesize \color{red} Communication \\ delay};
		\node[right=0cm of simulation, text width=5em, align=center] (ST) {\footnotesize \color{red} Simulation \\ time};
		\begin{pgfonlayer}{bg}
			\node[fit=(ncdetails)(CT),draw=gray,dashed,fill=gray!10] {};
			\node[fit=(simulation)(ST),draw=gray,dashed,fill=gray!10] {};
		\end{pgfonlayer}
	\end{tikzpicture}
    \caption[Decision-making process using numerical simulations.]{Decision-making process using numerical simulations. Additional time delays can be caused by communication and simulation. In an ideal situation: Communication delay + Simulation time $\ll$ Manufacturing time }
	\label{fig:decision-making}
\end{figure}

Analytical and empirical solutions are usually not available for complex mechanical/structural configurations, while experimental testing is cost-prohibitive. Therefore, numerical modelling is often the most practical route to aid the decision-making process. However, the numerical simulations consume a significant amount of time and computational resources in building the model and subsequent resolution, especially for assemblies that include dynamic contacts generating additional complexities. Also, the communication delay between the shop floor and the simulation team adds to the delays. These bottlenecks have been schematically described in \cref{fig:decision-making}. Therefore, it is necessary to accelerate assembly simulations in order to comply with desired decision times in the manufacturing cycle.

Real-time simulation tools have been of interest in various decision-making processes across fields like design and health-monitoring, among others. In such applications, fast simulations need to be performed for parametric values that are fed to the simulator in real-time. However, the complexity of a model usually increases exponentially with the number of parameters in the model. Therefore, high-fidelity models are often replaced with surrogate models that offer fast simulations, sometimes referred to as \emph{digital twins}. Reduced Order Modelling (ROM) a.k.a. Model Order Reduction, is a well-known framework to create such models. ROMs are generally built to efficiently solve a parametric problem in a \emph{multi-query} context. Therefore, ROM is a promising framework for moving towards efficient simulations of non-conforming engineering hardware.

For a real-time simulation of non-conforming assemblies, the task is two-fold. First is the modelling of geometric variabilities in individual parts. Shape parametrized ROMs have been explored~\cite{Lauzeral2019paper,Marconi2020,Navarro2022} which can be used for parametrizing geometric variabilities. The second, and more challenging, is the application of ROM approach to contact mechanics. As contact mechanics deals with multi-body problems with interaction between the boundaries of each body, it introduces geometric non-linearities. These non-linearities are different in nature compared to, say, non-linear material behaviour or loading conditions which have been widely studied by the ROM community. Therefore, contact ROMs still need extensive research before they can be used in industrial real-time applications, such as the acceleration of the decision-making process discussed previously.

Therefore, in this thesis, the current limitations of contact ROMs have been investigated and corresponding mitigation strategies have been proposed. The thesis begins with an introduction to contact mechanics and the state of the art in ROMs applied to contact mechanics (Chapter 1). Then, the limitations of the so-called \emph{low-rank} approach are  investigated (Chapter 2). Strategies to counter these limitations are developed using over-complete dictionaries (Chapter 3) and non-linear interpolations (Chapter 4).


\renewcommand{\thefigure}{\oldfig}

\chapter{State of the Art}\label{ch:literature}
In this chapter, an overview of concepts and state of the art related to contact mechanics and reduced order modelling are presented. Ideas in contact mechanics are first introduced using a two-body problem where inequality constraints apply to the energy minimization statement of the standard mechanics problem. Some of the classical methods to solve inequality constrained problems are discussed. As one of the principal sources of complexity in such problems is contact detection, some methods applicable to contacts between one-dimensional boundaries are introduced. In the second part of the chapter, the concept of Reduced Order Modelling and its application on unconstrained problems is discussed. In the final part of the chapter, the two elements are combined, where the developments in reduced models of contact mechanics are presented. 

\section{Overview of Contact Mechanics}
Many computational mechanics problems involve interaction between multiple bodies with various types of physics, say thermal or mechanical, where each kind of physical interaction adds a challenge to the modelling process. Contact mechanics, the field that considers mechanical interaction in multi-body problems, finds applications in a wide variety of fields as almost all mechanical system, natural or man-made, displays some sort of contact phenomena. A few examples that fall in this category include: foundation studies in civil engineering, human joints in biomechanics, roller joints in mechanisms, adhesive and frictional contacts in tribology and crash simulations in automotive industry. 

Historical works in contact mechanics include studies by that of~\citet{Hertz1882} on analytical solutions of contact between spherical frictionless surfaces of elastic bodies, that of~\citet{Coulomb1973} on frictional contact, and many more. However, the work by~\citet{Signorini1959} on the contact problems with elastic bodies, which consisted of the Boundary Value problem with non-penetration constraints, is central to approach of contact mechanics in the computational field. The problem was titled ``problems with ambiguous boundary conditions.'' by Signorini, aptly describing the nature of the problem. Analytical solutions to certain contact problems have been computed in~\cite{Johnson1985,Johnson1992,Goryacheva1998}.

\subsection{Parametrized contact mechanics problems}
\label{sec:model_prob}
Usually, an engineering problem is described as \emph{parametrized} if a parametric solution is sought i.e.\, the solution of the problem is sought for a range of parameter values. For a mechanical problem, the parameters can be any attribute of the mechanical system such as material properties, geometry, loading or boundary conditions, etc. Parametric solutions usually find applications in design optimization, real-time monitoring and inverse problems (see \cref{sec:rom_literature}). In this section, parametric mechanical problems with non-penetration constraints on displacement field are discussed. Linear elastic and small deformation problems are considered, keeping the evaluation of internal energy straightforward. The inequality constraints may be a result of the presence of either an obstacle or a second body in the domain, with possibility of contact between different surfaces in the deformed configuration. 

\begin{figure}
	\begin{subfigure}[t]{0.48\linewidth}
		\centering
		\includegraphics[scale=0.7,bb=5 48 148 248]{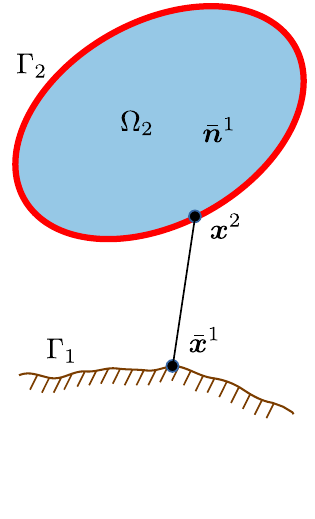}
		\caption{Type 1: Body-obstacle problem}
		\label{fig:body-obstacle}
	\end{subfigure}
	\begin{subfigure}[t]{0.48\linewidth}
		\centering
		\includegraphics[scale=0.7,bb=1 1 148 248]{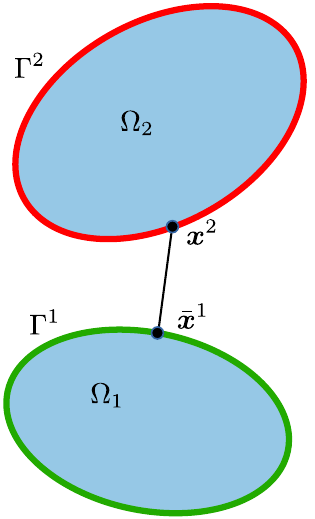}
		\caption{Type 2: Two-body problem}
		\label{fig:two-body}
	\end{subfigure}
    \caption[Kinematic description of mechanical problems with possibility of contact]{Kinematic description of mechanical problems with possibility of contact. Points $\bm{x}^2$ will possibly come into contact with point $\bar{\bm{x}}^1$, and hence considered as a contact pair. The distance between the two points in the initial and deformed configuration is given by $g(\mu,\bm{v})(\bm{x}^2)$ and $\left (g(\mu,\bm{v})(\bm{x}^2) - k(\mu,\bm{v};\bm{v})(\bm{x}^2) \right )$ respectively}
	\label{fig:contact-schematic}
\end{figure}

For convenience, the contact problems will be referred to as Type-1, referring to a problem involving a deformable body and an immovable obstacle, and Type-2 referring to problems involving two deformable bodies. A schematic of the two types of contact problems is shown in \cref{fig:contact-schematic}. Contact phenomena like friction and adhesiveness are neglected. The resolution of contact mechanics problems depends on computing the distance between surface of bodies involved. As evaluation of distance is not a trivial task because of dynamic contact pairing, different strategies exist based on underlying simplifications. The details to these approaches can be found in~\cite{Wriggers2006,Fischer2006,Yastrebov2013}.

The generic weak form of a mechanical problem with a parameter $\mu$ involving contacts (for both Type-1 and Type-2) is described by the following inequality constrained minimization problem:
\begin{align}
    \begin{aligned}
        & \bm{u} = \argmin_{\bm{v}\in \PrimalSpace} \frac{1}{2} a(\mu;\bm{v},\bm{v}) - f(\mu;\bm{v}) &  \\
        \text{s.t. \ \ \ }& k(\mu,\bm{v};\bm{v})(\bm{x}) \leq g(\mu,\bm{v})(\bm{x})  & \text{on } \Gamma^2
    \end{aligned}
    \label{eq:constr_minimization}
\end{align}
where
\begin{itemize}
    \item $a(\mu;\bm{v},\bm{v})$ is strain energy function and $f(\mu;\bm{v})$ is the work contribution by the external forces. The operators are defined $\forall \bm{v} \in \PrimalSpace$ and $\forall \mu \in \mathcal{P}$, where $\PrimalSpace$ be an appropriate function space for displacement $\bm{v}$ and $\mathcal{P}$ is the parametric space for the parameter $\mu$. Such problems with parametric dependence are usually referred to as parametrized problems. Note that $\mu$ may not be necessarily a scalar and the parametric space $\mathcal{P}$ can be possibly multidimensional.
    \item $k(\mu,\bm{v};\bm{v})(\bm{x})$ and $g(\mu,\bm{v})(\bm{x})$ are the distance functions that indicate separation between bodies, defined conveniently on surface $\Gamma^2$. The first term indicates the contribution of displacement field $\bm{v}$ to the distance and the second term indicates the distance in the undeformed configuration. For a Type-2 problem, the distance functions can be expressed as:
        \begin{align}
            \begin{aligned}
                &k(\mu,\bm{v};\bm{v})(\bm{x}) &&= \left[ \bm{v}^2(\bm{x}^2) - \bm{v}^1(\bar{\bm{x}}^1) \right]  \cdot \bm{n}\\
                &g(\mu,\bm{v})(\bm{x}) &&=  \left [ \bm{x}^2 - \bar{\bm{x}}^1 \right ] \cdot \bm{n}
            \end{aligned}
            \label{eq:distance_functions}
        \end{align}
        The pair $(\bm{x}^2, \bar{\bm{x}}^1)$  and the normal $\bm{n}$ are determined on the discrete domain using various methods, some of which are given in upcoming \cref{sec:contact_pairing}. For a Type-1 problem, the displacement field of the obstacle $\bm{v}^1$ vanishes.
\end{itemize}

\begin{notation}
    \label{not:strong_form}
    Notations to distinguish strong and weak forms are adopted from~\cite{Benaceur2020}. Linear and bilinear operators that appear in weak formulations are expressed in the format $o (\cdot;\cdot,\cdot)$, where the non-linear dependencies and linear dependencies are separated by the semi-colon ';' . On the other hand, strong form functions  are expressed in the format $o (\cdot;\cdot,\cdot)(\cdot)$, distinguished from weak forms using an additional argument which is the geometric position. The same applies to linear forms expressed as $o (\cdot;\cdot)$ and $o (\cdot;\cdot)(\cdot)$, where the linear dependency is on the second argument.
\end{notation}

\subsection{Contact Formulations}\label{sec:contact_formulations}
To treat the contact problem~\eqref{eq:constr_minimization} numerically, the contribution of contact to the total energy of the system is computed. This additional energy term is defined differently in each contact formulation, however, the basic idea is to introduce contact forces, directly or indirectly, that will prevent penetrations of contact surfaces. Among numerous formulations available, three of them are introduced in this section. Two of them are based on introducing the Lagrange Multiplier, which computes the contact forces directly, while the other method penalizes the penetrations.
\subsubsection{Lagrange Multiplier Method}\label{sec:lmm}
One of classical methods is an attempt to accurately satisfy the inequality constraints of the contact problem (or any minimization problem) is the Lagrange Multiplier Method (LMM). The resulting optimization problem can be expressed as a saddle point problem of the Lagrangian functional:
\begin{align}
    (\bm{u}, \lambda ) = \argminmax_{\bm{v}\in \PrimalSpace, \eta\in \DualSpace } \mathcal{L}(\bm{v},\eta)
    \label{eq:saddle_point}
\end{align}
where  
\begin{align}
    \mathcal{L}(\bm{v},\eta) = \frac{1}{2}  a(\mu;\bm{v},\bm{v}) - f(\mu;\bm{v}) + b(\mu,\bm{u};\eta,\bm{u}) - d(\mu,\bm{u};\eta)
    \label{eq:lagrangian}
\end{align}
where the distance function is expressed by following weak forms:
\begin{gather}
    \begin{aligned}
        b(\mu,\bm{u};\eta,\bm{u}) & = \int_{\Gamma^2} \eta \ k(\mu,\bm{u};\bm{u})(\bm{x})  \ \partial \Gamma  \\
        d(\mu,\bm{u};\eta) & = \int_{\Gamma^2} \eta \ g(\mu,\bm{u})(\bm{x})  \ \partial \Gamma
    \end{aligned}
    \label{eq:dist_func_weak}
\end{gather}
where $\DualSpace$ is an appropriate non-negative function cone for contact pressure defined on the surface $\Gamma^2$. The inequality $\lambda \geq 0 $ admits a physical meaning because $\lambda$ is equivalent to contact pressure. A negative contact pressure implies traction, in other words, adhesion between contact surfaces which contradicts the simplifying assumptions.

\noindent The KKT conditions associated to the optimizer $(\bm{u}, \lambda) \in \PrimalSpace \times \DualSpace$ of~\eqref{eq:saddle_point} can be expressed as follows:
\begin{subequations}
	\begin{align}
        a(\mu;\bm{v},\bm{u}) - f(\mu;\bm{v}) + b(\mu,\bm{u};\lambda,\bm{u}) &= 0 , && \bm{v} \in \PrimalSpace \\
        b(\mu,\bm{u};\eta,\bm{u}) - d(\mu,\bm{u};\eta) 	   &  \leq 0 ,  && \eta \in \DualSpace \\
        b(\mu,\bm{u};\lambda,\bm{u}) - d(\mu,\bm{u};\lambda)   &   = 0 && \label{eq:contact_energy_kkt_comp_slack}
	\end{align}  \label{eq:kkt_weak_form}
\end{subequations}

\noindent Notice that the last KKT conditions~\eqref{eq:contact_energy_kkt_comp_slack} implies that contribution of contact energy term will be zero at the solution. This makes sense because the non-penetration constraint can prevent the system from reaching the true minima of the energy functional, but it cannot add energy to the system.

\subsubsection{Penalty Method}
To constrain the solution in the feasible region defined by the inequality constraints, violations of inequality constraints can be penalized by adding a penalty term to the energy of the system. The resulting system is an unconstrained problem, unlike the LMM problem. 
\begin{align}
	\begin{aligned}
        \bm{u}  &=& \arg \min_{\bm{v}\in \PrimalSpace} \frac{1}{2}  a(\mu;\bm{v},\bm{v}) - f(\mu;\bm{v})  & + \int_{\Gamma^2} \tau \left[\left(k(\mu,\bm{v};\bm{v})(\bm{x}) - g(\mu,\bm{v})(\bm{x})\right)^{+}\right]^2  \ \partial \Gamma \\
	\end{aligned} \label{eq:contact_penalty}
\end{align}
where $
	(z)^+ :=\left \lbrace\begin{matrix}
		z & \text{if } z>0 \\
		0 & \text{otherwise}
	\end{matrix}  \right . $

The solution of~\eqref{eq:contact_penalty} may not necessarily satisfy the inequality constraints exactly, but may be nearly compliant for sufficiently large values of the penalty parameter $\tau$. For very large values of $\tau$, the penalty solution should approach the LMM solution, however, the numerical system are ill-conditioned in this regime of $\tau$.

The penalty parameter $\tau$ can be interpreted as a linear spring that resists penetration. Therefore, a natural variation of the penalty method can be built by using a non-linear penalty function where the penetration and the resultant contact pressure have a non-linear relation. This approach is useful to reduce the penetration effects, but also in cases where a constitutive law relating contact pressure and penetration is known. However, the problem of ill-conditioning still persists when the constitutive law enters a highly stiff regime~\cite[Chapter~5,~6]{Wriggers2006}.

\subsubsection{Augmented Lagrange Multiplier Method}
To circumvent the ill-conditioning problems of penalty method, Augmented Lagrange Multiplier Method~\cite{Hestenes1969} problems have been developed. This class of method use a combination of penalty and LMM\@. The basic idea is to establish contact forces using the penalty method and then transfer these forces, in multiple steps, into the Lagrange Multiplier, thereby reducing illegal penetrations. The advantage of using this method is that the solution will have high compliance to the constraints irrespective of the magnitude of the penalty parameter. At a given step $k$, the following problem is minimized:
\begin{gather}
	\begin{align*}
        \bm{u}^k  &=& \arg \min_{\bm{v}\in \PrimalSpace} \frac{1}{2}  a(\mu;\bm{v},\bm{v}) - f(\mu;\bm{v})  & + \int_{\Gamma^2} \tau \left[\left(k(\mu,\bm{v};\bm{v})(\bm{x}) - g(\mu,\bm{v})(\bm{x})\right)^{+}\right]^2 \ \partial \Gamma\\[1.0em]
        &  & &+b(\mu,\bm{v};\lambda^k,\bm{v}) - d(\mu,\bm{v};\lambda^k)
	\end{align*} \\
    \begin{align*}
        \intertext{with the update step:}
        \lambda^{k+1}(\bm{x}) &= \lambda^k(\bm{x}) + \tau \left[\left(k(\mu,\bm{v};\bm{v})(\bm{x}) - g(\mu,\bm{v})(\bm{x})\right)^{+}\right]  \\
        \intertext{and initial conditions:}
        \lambda^{0}(\bm{x}) &= 0
    \end{align*}
\end{gather}
As $\lambda^k$ does not contain the true Lagrange Multiplier during the intermediate steps, this method still allows some penetration. However, as the Lagrange Multipliers are updated with new contact pressure forces, the penetration reduces in the subsequent steps. With sufficient number of steps, penetration diminishes to very low-levels and the solution approaches that of the LMM solution.
\subsection{Discrete Formulations}\label{sec:contact_discretization}
In this section, the usual discretization strategies based on Finite Element (FE) are briefly introduced. One of the foremost problems in contact mechanics is the detection of contact pairs. For large-scale problems, this step can prove to be expensive. Then comes the part of building the finite element (FE) operators that can be used to enforce or verify the constraints. As the contact problems in upcoming chapters will be solved by using the LMM approach, the discussion on discretization schemes will be limited to this method. Discretization of the penalty and Augmented LMM can be carried out in a similar manner. Details of all approaches can be found in~\cite{Wriggers2006,Yastrebov2013}.
\subsubsection{Contact Detection and Pairing}\label{sec:contact_pairing}
\begin{figure}[htpb!]
    \centering
    \includegraphics[width=0.5\linewidth,bb=16 6 415 254]{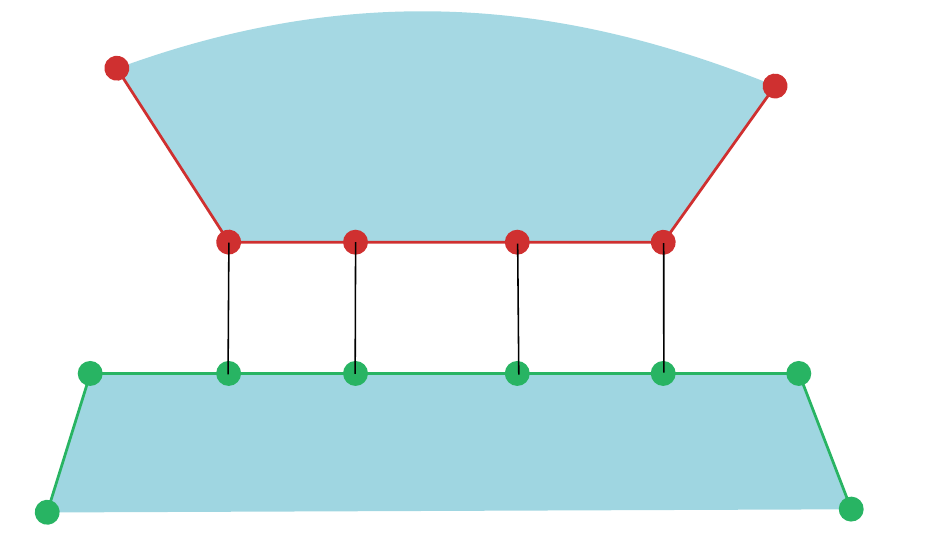}
    \caption{Schematic of node-to-node contact pairing for conforming surface meshes}
	\label{fig:n2n}
\end{figure}

Establishing the region of contact between two bodies is one of the primary challenges in contact problem. In simplified cases of small displacement, conforming meshes can be used, thereby allowing the use of \emph{node-to-node} contacts. In this case, both points in the contact pair $\bm{x}^2$ and $\bar{\bm{x}}^1$ are nodes over respective surfaces, as shown in \cref{fig:n2n}. Naturally, the contact pair can be predetermined based on the discretization. Obviously, the use of conforming mesh limits the applicability of this method. For example, it cannot be applied to sliding surfaces.

\begin{figure}[htpb!]
	\begin{subfigure}[t]{0.48\linewidth}
		\centering
		\includegraphics[width=1\linewidth,bb=-1 0 448 259]{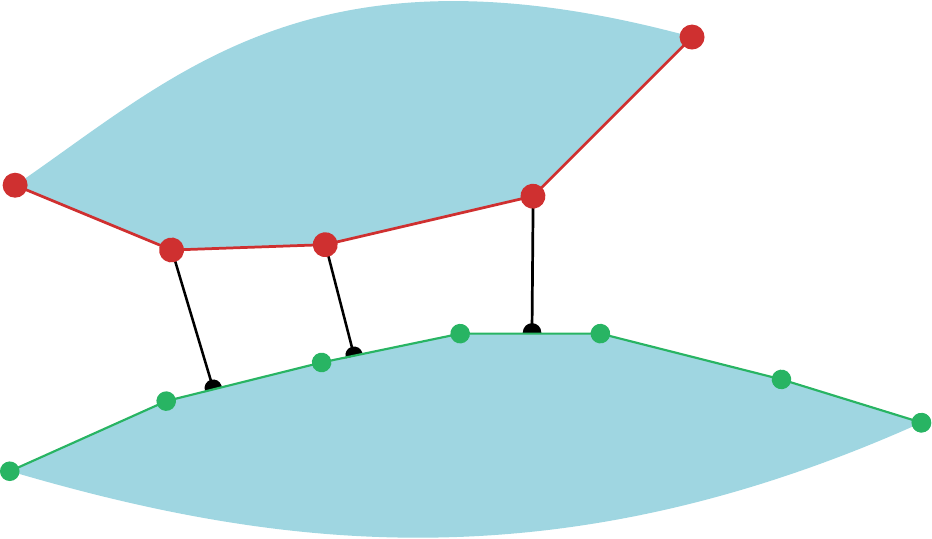}
		\caption{Closest point projection}
		\label{fig:cpp}
	\end{subfigure}
	\begin{subfigure}[t]{0.48\linewidth}
		\centering
		\includegraphics[width=1\linewidth,bb=-1 0 448 259]{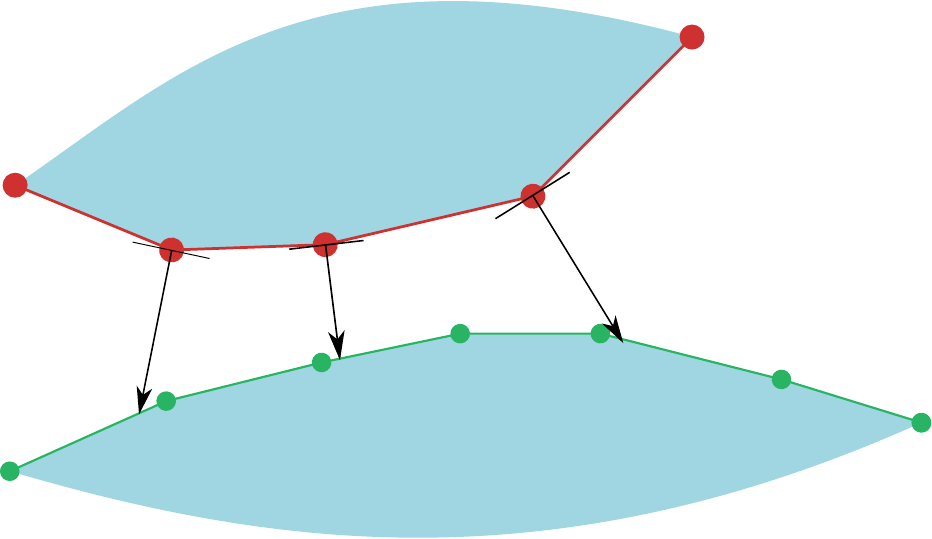}
		\caption{Ray tracing}
		\label{fig:ray_tr}
	\end{subfigure}
    \caption{Schematics of node-to-segment methods for finding contact pairs}
	\label{fig:contact-detection}
\end{figure}

\emph{Node-to-segment} method is used in cases of non-conforming meshes with 1D boundaries, for e.g., 2D problems with possibility of contact on the edges of bodies. The contact pairings $(\bm{x}^2, \bar{\bm{x}}^1)$ are not straightforward to compute as it was in node-to-node method. The contact pair evaluation is one of the sources of complexities in contact problems, especially for large displacement problems as contact pairs can be a strong function of displacement field. Established methods of computing contact pairs include:
\begin{itemize}
    \item \textbf{Closest Point Projection}: For each node $\bm{x}^2$ on $\Gamma^2$, the closest point $\bar{\bm{x}}^1$ on $\Gamma^1$ is chosen as the pair~\cite{Hallquist1985}. Naturally, the distance functions are evaluated using~\eqref{eq:distance_functions} with $\bm{n} = \bar{\bm{n}}^1$, the normal at the point $\bar{\bm{x}}^1$. A limitation of this method is the possibility of non-unique closest point in certain geometries~\cite[Chapter~2]{Yastrebov2013}.
    \item \textbf{Ray Tracing}: At each node $\bm{x}^2$ of $\Gamma^2$, the local normal $\bm{n}^2$ is extended to intersect $\Gamma^1$ at $\bar{\bm{x}}^1$~\cite{Poulios2015}. Hence, distance functions are evaluated using $\bm{n} = -\bm{n}^2$ for this method. This method may not find the nearest point, but avoids the problems of non-uniqueness and projecting on irregular points. 
\end{itemize}
The two methods are shown in \cref{fig:contact-detection}. The normals $\bm{n}^1$ and $\bm{n}^2$ are evaluated in deformed configuration. The pairing ($\bar{\bm{x}}^1$, $\bm{x}^2$) is hence dependent on the displacement field in both methods. This forces the distance function $k$ to have an implicit non-linear dependence on displacement, apart from the explicit linear dependence that is evident in the expression. 

\subsubsection{Discretization using Finite Elements}\label{sec:mixed_form_FE}
\noindent The discrete equivalent of the Lagrangian functional in~\eqref{eq:lagrangian} can be derived by introducing finite element spaces $\PrimalSpace^h\subset \PrimalSpace$ and $\DualSpace^h \subset \DualSpace$:
\begin{align}
    \mathcal{L}^h(\bm{v}, \bm{\eta}) = \frac{1}{2} \bm{v}^T \mathbf{K}(\mu) \bm{v} - \bm{v}^T \bm{f}(\mu) &+ \bm{\eta}^T \left[ \mathbf{C}(\mu,\bm{v}) \bm{v} - \bm{g}(\mu,\bm{v})  \right] \label{eq:discrete_lagrangian}
\end{align} 

\noindent where $\bm{v} \in \mathbb{R}^{\ndofprimal}$ and $\bm{\eta} \in \mathbb{R}^{\ndofdual}$. Here, $\ndofprimal$ and $\ndofdual$ indicate the number of respective finite element degrees of freedoms (dofs). Also, note that the notation $\bm{u}$ and $\bm{v}$ are used for both continuous and discrete versions of displacement fields. 

\noindent Let $\lbrace \bm{N}_1, \bm{N}_2 \ldots \bm{N}_{\ndofprimal} \rbrace$  and   $\lbrace M_1, M_2 \ldots M_{\ndofdual} \rbrace$ denote the basis for the FE space $\PrimalSpace^h$ and cone $\DualSpace^h$ respectively. Here, for a system with $n$ displacement dofs per node, $\bm{N}_i : \Omega \to \mathbb{R}^n$ and $M_j : \Gamma \to \mathbb{R}$. The discrete operators $\mathbf{K}(\mu)$ and $\mathbf{C}(\mu,\bm{v})$ can be defined as:

\begin{alignat*}{3}
	&\mathbf{K}(\mu)_{ij} &&= a(\mu;\bm{N}_i,\bm{N}_j) \\
	&\mathbf{C}(\mu,\bm{v})_{ij} &&=  b(\mu,\bm{v};M_i, \bm{N}_j)
	\intertext{and the vectors $\bm{f}$ and $\bm{g}$ are defined as:}
	&\bm{f}(\mu)_i &&= f(\mu;\bm{N}_i) \\
	&\bm{g}(\mu,\bm{v})_i &&= d(\mu,\bm{v};M_i )
\end{alignat*}
The bilinear and linear terms $b$ and $d$ are non-linear w.r.t.\ the displacement field $\bm{v}$ and the parameter $\mu$, which is reflected in discrete form as $\mathbf{C}(\mu,\bm{v})$ and $\bm{g}(\mu,\bm{v})$ respectively. Using the above operators, the discrete form of the KKT conditions in~\eqref{eq:kkt_weak_form} for the optimal solution $(\bm{u}, \bm{\lambda}) \in \mathbb{R}^{\ndofprimal} \times \mathbb{R}^{\ndofdual}$ can be written as:

\begin{subequations}
	\begin{align}
        \mathbf{K}(\mu) \bm{u} - \bm{f}(\mu) + \mathbf{C}^T(\mu,\bm{u}) \bm{\lambda} &= \bm{0} \label{eq:kkt_discrete_balance} \\
        \bm{\lambda} &\geq \bm{0} \label{eq:kkt_discrete_nonnegativity}\\
        \mathbf{C}(\mu,\bm{u}) \bm{u} - \bm{g}(\mu,\bm{u}) &\leq \bm{0} \label{eq:kkt_discrete_nonpen}\\
        \bm{\lambda}^T [\mathbf{C}(\mu,\bm{u}) \bm{u} - \bm{g}(\mu,\bm{u})] &= 0 \label{eq:kkt_discrete_comp_slack}
	\end{align} \label{eq:kkt_discrete}
\end{subequations}

Apart from the non-linearity of the constraint operators $\mathbf{C}$ and $\bm{g}$, additional non-linearities appear because the regions where the inequality constraints are active and inactive are unknown. Active regions correspond to the regions where contact is established, or in other words, constraints satisfy the equality. These non-linearities are potential contributors to the overall computational complexity. \\

\noindent \textbf{Stability of the saddle point problem}: The stability conditions of mixed formulations usually appear in the form of an inf-sup condition. The inf-sup constant of the bilinear operator $b$ in~\eqref{eq:lagrangian} must be positive for every pair of test functions in the two solution spaces. Not every combination of finite-element spaces $\PrimalSpace^h$ and $\DualSpace^h$ satisfies the stability condition. In this thesis, linear elements for displacement and node-centered zero order (piecewise-constant) shape functions for  contact pressure are used (based on Reference~\cite{Benaceur2020}), which generates smooth solutions for two-dimensional problems. More details on inf-sup stability of contact problems can be found in~\cite{ElAbbasi2001,Wriggers2006}.

\subsection{Resolution schemes}
Even for contact problems involving small deformations with linear elastic constitutive law, geometric non-linearities can still occur if large displacements are involved. This is because the contact pairs change with the displacement field. Even for small displacement problems where contact pairs can be approximately precomputed in the reference configuration, the region of contact cannot be predetermined. Linearization schemes are needed to solve the discrete problem, that resolve both the contact region and the contact pairs in an iterative manner. Two such schemes: Fixed-point and Newton-Raphson iterations are described here. Python notation is used to indicate slicing operations. 

\subsubsection{Fixed point method}
Fixed point schemes are used to find solutions to non-linear equations by posing them in the fixed-point form ${x}=f({x})$, leading to a fairly simple iterative form ${x}^{p+1} = f({x}^{p})$. For a discrete system with multiple dofs, the linearized form usually looks like $\mathbf{A}(\bm{x}^{p}) \bm{x}^{p+1} = \bm{b}(\bm{x}^{p})$. These linearized forms can be iterated until convergence, without computing derivatives or gradients of the non-linear equation.

Fixed point methods for an inequality constrained problem can be implemented using active set approach~\cite{Nocedal2006,Luenberger2016}, in which the constraints are divided into active and inactive sets. For the discretized contact problem, the balance form~\eqref{eq:kkt_discrete_balance} can be used in the fixed-point method:
\begin{align}
    \mathbf{K}(\mu) \bm{u} +  \mathbf{C}^T(\mu, {\bm{u}}) \bm{\lambda} = \bm{f}(\mu) \label{eq:contact_discrete_equilibrium}
\end{align}

\noindent Given the current estimate of active set $\ActiveSet$, the activation state of the inequality constraints can be written as:
\begin{subequations}
    \begin{align}
        \bm{\lambda}[\ActiveSet^\mathsf{c}] &\gets \bm{0}  \\
        \mathbf{C}_\ActiveSet(\mu,  {\bm{u}}) \bm{u} - \bm{g}_\ActiveSet(\mu, \bm{u}) &= \bm{0} \label{eq:contact_active_constraint}
    \end{align}
\end{subequations}
\noindent where  
\begin{align*}
    \mathbf{C}_\ActiveSet(\mu,  {\bm{u}}^p) &= \mathbf{C}(\mu,  {\bm{u}}^p)[\ActiveSet,:]  \\
    \bm{g}_\ActiveSet(\mu,  {\bm{u}}^p) &= \bm{g}(\mu,  {\bm{u}}^p)[\ActiveSet]  \\
\end{align*}

\noindent and $\ActiveSet^{\mathsf{c}}$ is the complement of set $\ActiveSet$. 

These conditions apply because each entry of the vectors $\mathbf{C}(\mu, \bm{u}) \bm{u}-\bm{g}(\mu,\bm{v})$ and $\bm{\lambda}$ correspond to a nodal constraint and a nodal Lagrange Multiplier respectively. As the active set contains the list of active nodes, it can be used for filtering the entries of constraint operators and enforcing the equality conditions i.e.,\ where contact between bodies is established. The displacement field and the active Lagrange multipliers at iterative level $p+1$ can be obtained using the linear system obtained from the equilibrium condition~\eqref{eq:contact_discrete_equilibrium} and the active constraint~\eqref{eq:contact_active_constraint}


\begin{align}
	\label{eq:contact_fixedpoint}
	&\begin{bmatrix}
		\mathbf{K}(\mu) &  \mathbf{C}_\ActiveSet^T(\mu, {\bm{u}}^p) \\
		\mathbf{C}_\ActiveSet(\mu,  \bm{u}^p)  & \mathbf{0}
	\end{bmatrix}
	\begin{bmatrix}
		\bm{u}^{p+1} \\
		\bm{\lambda}^{p+1}_\ActiveSet
	\end{bmatrix} =  \begin{bmatrix}
		\bm{f}(\mu) \\
		\bm{g}_\ActiveSet(\mu, \bm{u}^p)
	\end{bmatrix}
\end{align}
\noindent and the discrete contact pressure can be updated as:
\begin{align*}
    \bm{\lambda}^{p+1}[\ActiveSet] &\gets \bm{\lambda}^{p+1}_\ActiveSet  \\
    \bm{\lambda}^{p+1}[\ActiveSet^\mathsf{c}] &\gets \bm{0} 
\end{align*}

The sign of entries in the expression $\mathbf{C}(\mu, \bm{u}) \bm{u} - \bm{g}(\mu,\bm{v})$ and in the vector $\bm{\lambda}$ correspond to satisfaction or violation of inequality constraints at the respective node's support. The active set for the next iteration is obtained by verifying the constraints:
\begin{align}
    \ActiveSet_i \gets \left\{ i \ \left| \left(\mathbf{C}(\mu, \bm{u}^{p})\bm{u}^{p+1} - \bm{g}(\mu, \bm{u}^{p})\right)_i \geq 0  \text{ and }  \bm{\lambda}_i > 0 \right . \right\}
\label{eq:active_set_update}
\end{align}

\noindent Moreover, the contact pairs are also be updated using the latest displacement field and the new operators $\mathbf{C}(\mu, \bm{u}^{p+1})$ and $\bm{g}(\mu, \bm{u}^{p+1})$ are computed. This process is continued until convergence.

\subsubsection{Newton-Raphson}
The Newton-Raphson iterative method (a.k.a. Newton's method) is a derivative-based iterative method for root finding and optimization problems. For finding the roots of an equation $f(x) = 0$, the Newton-Raphson update at an iterative level $p$ can be written as $\Delta x^{p+1} = -  \frac{f(x^p)}{f'(x^p)}$, where $\Delta x^{p+1} = x^{p+1} - x^p$. For a discrete system with residual $\bm{b}(\bm{x}) - \mathbf{A}(\bm{x}) \bm{x}$, the Newton-Raphson update step can be written as: 
\begin{align*}
    \mathbf{A}_T(\bm{x}) \Delta \bm{x}^{p+1} = \bm{b}(\bm{x}^p) - \mathbf{A}(\bm{x}^p) \bm{x}^p
\end{align*}
\noindent where $\mathbf{A}_T(\bm{x})$ is the tangent matrix of the residual.

The Newton-Raphson method for contact problems is based on computing the tangent matrix $\mathbf{K}_T^c$ of the residual of equilibrium equation~\eqref{eq:contact_discrete_equilibrium} and then computing incremental updates. The contribution to the tangent matrix $\mathbf{K}_T$ comes from of the contact part of the energy functional or the equilibrium equation. It consists of the second order terms that appear due to variations in contact pairing. The computation of tangent matrix for the non-linear system is quite detailed and can be found in Reference~\cite[Chapter~9,~10]{Wriggers2006}. The iterative update at a level $p+1$ can be written as:

\begin{align*}
\begin{aligned}
	\begin{bmatrix}
		\mathbf{K}(\mu) + \mathbf{K}_T^c(\mu,  {\bm{u}}^p)  &  \mathbf{C}_\ActiveSet^T(\mu, {\bm{u}}^p) \\
		\mathbf{C}_\ActiveSet(\mu,  {\bm{u}}^p)  & \mathbf{0}
	\end{bmatrix}
	& \begin{bmatrix}
		\Delta \bm{u}^{p+1} \\
		\Delta \bm{\lambda}^{p+1}_\ActiveSet
	\end{bmatrix} \\
    =  & \begin{bmatrix}
        \bm{f}(\mu) \\
        \bm{g}_\ActiveSet(\mu, \bm{u}^p)
	\end{bmatrix} -  \begin{bmatrix}
		\mathbf{K} &  \mathbf{C}_\ActiveSet^T(\mu, {\bm{u}}^p) \\
		\mathbf{C}_\ActiveSet(\mu,  {\bm{u}}^p)  & \mathbf{0}
	\end{bmatrix}
	\begin{bmatrix}
		\bm{u}^p \\
		\bm{\lambda}^p_\ActiveSet
	\end{bmatrix}
\end{aligned}
\end{align*}
The stiffness matrix $\mathbf{K}(\mu)$ appears in addition to the contact tangent matrix because the stiffness matrix is the tangent of the elastic energy functional. After computing the increment, the state vector can be updated as:
\begin{align*}
	\begin{bmatrix}
		\bm{u}^{p+1} \\
		\bm{\lambda}^{p+1}_\ActiveSet
	\end{bmatrix} = 
	\begin{bmatrix}
		\bm{u}^{p} \\
		\bm{\lambda}^{p}_\ActiveSet
	\end{bmatrix} + 
	\begin{bmatrix}
		\Delta \bm{u}^{p+1} \\
		\Delta \bm{\lambda}^{p+1}_\ActiveSet
	\end{bmatrix}
\end{align*}

\noindent This is followed by the active set update in~\eqref{eq:active_set_update}.

\section{Reduced Order Modelling for parametric problems}\label{sec:rom_literature}
In recent applications of numerical methods to engineering design, decision-making processes are increasingly dependent on simulation results~\cite{Guratzsch2010,Modesto2015,Mosavi2015,Bigoni2020}. Various configurations of the product can be explored and compared by optimization algorithms to develop an optimal configuration using modern simulation tools~\cite{Chinesta2013,Nguyen2014,Woldemariam2019}. For the integration of simulation tools into the decision-making process, simulation time for each configuration becomes a primary factor in determining the swiftness of the process. In this context of \emph{multi-query} simulations, a high computational complexity of traditional simulation methods might discourage comprehensive exploration of the design space for economic reasons.

Traditional methods of simulations were restricted to the design phase of engineering developments due to the heavy computational costs. On the other hand, modern simulation tools capable of real-time simulations have promising applications not only in the design phase but also during the life-cycle of engineering products, for example, for operations and maintenance~\cite{Hartmann2018}. The idea behind such applications is to create real-time simulations whose predictions can be used to initiate corrective or preventive actions.  At the heart of the development of real-time simulations lies the approach of \emph{Reduced Order Modelling} (ROM).

Real-time simulators need to mimic the relevant physics of the system being studied taking into account the parameters that are fed to the simulator in real-time. High-fidelity numerical methods usually have very high computational complexity for such applications. The real-time query can be addressed by solving a high-dimensional problem, in which the real-time parameters are treated as coordinates of the problem (like space and time). Such a numerical model is usually referred to as a parametrized model and is usually very expensive to resolve using standard numerical methods. ROM is a framework for reducing the complexity of high-dimensional problems. ROMs have found successful applications in a wide variety of problems such as design optimization~\cite{McBane2021}, uncertainty quantification~\cite{Freitas2016}, inverse problems~\cite{Gonzalez2012}, optimal control~\cite{Karcher2017}, real-time monitoring~\cite{Modesto2015}, and others.

Various methodologies have been adopted to create ROMs to address different bottlenecks. Some are based on reducing the number of dofs of the system, whereas others focus on the efficient evaluation of integrals over a small part of the domain. However, a common feature among most ROM methods has been the computation of a training set of \emph{snapshots} i.e.\ a set of high-fidelity solutions based on which a reduced model is built. A brief description of a few classes of reduction methods is given here:
\begin{itemize}
    \item \emph{Low-rank methods}: These methods primarily focus on finding a low-rank solution subspace. Usually, this is done by computing the basis of the subspace using Singular Value Decomposition (SVD) or Proper Orthogonal Decomposition (POD). Low-rank methods are further detailed in \cref{sec:low_rank}.
    \item \emph{Hyper-reduction methods}: Primary aim of hyper-reduction methods is to reduce the number of unknowns as well as quadrature points to compute weak form terms by generating a reduced integration domain (RID), which is generated by selecting appropriate dofs~\cite{Ryckelynck2005}. This feature of hyper-reduction is used in~\cite{Ryckelynck2009} to reduce the effort in the computation of internal variables of non-linear constitutive models, where the RID is made up of regions containing maximum energy of POD basis vectors. On the other hand, the quadrature rules for RID in Reference~\cite{Hernandez2017}, are computed using the least-square procedure so that all training samples are exactly integrated. Cheaper evaluation of weak forms over the RID leads to an efficient model for non-linear problems.
    \item \emph{Non-intrusive methods}: Some methods have been developed that do not need access and modification of the high-fidelity source codes, and are termed as \emph{non-intrusive} methods. These are highly useful if proprietary source codes or complex constitutive models are used. Also, it can be applied to different kind of physics without changing the non-intrusive framework. As the numerical model and its operators are not accessible, many non-intrusive methods rely on fitting an interpolation in the parametric domain. Radial basis functions (RBF) have been popularly used in this context. In these works, a POD-based solution subspace is computed and an RBF-based interpolation for the reduced coefficients is developed for the POD subspace~\cite{Audouze2009,Audouze2013,Xiao2015}. Another approach is that of sparse subspace learning where a predefined set of hierarchical parametric basis functions with special quasi-interpolative properties, allowing enrichment of the reduced model at minimal cost~\cite{Borzacchiello2017}. The CUR decomposition method takes a different approach for generating the parametric interpolation. A set of selected spatial and temporal snapshots are used to build a surrogate model with coefficients that fit the snapshots in the training set~\cite{LeGuennec2018}.
    \item \emph{Projection-based a priori method}: The Proper Generalized Decomposition (PGD) treats the parameters of the problem as coordinates and yet circumvents the \emph{curse of dimensionality} by using separated representations for solving the PDE~\cite{Ammar2007}. This is accomplished by means of a greedy approach, where the current solution is updated a with rank-$1$ tensor product of separated functions in each step. The rank-$1$ updates are computing using the separated form of PDE, one direction at a time, where each direction corresponds to one or more coordinates of the problem, for instance, position, time, or a parameter, or a combination of these. This is described as an \emph{a priori} method because a parametric solution is built and stored, resulting in quick accessibility of the solution. Efficient computation of separated form of non-linear terms has been implemented in the PGD framework using cross-approximations~\cite{Aguado2019}.
\end{itemize}

Even though these approach have different names and methodologies, all of them, directly or indirectly, assume the low-rank hypothesis. In the next section, low-rank methods are described for generic unconstrained problems.

\subsection{Low-rank methods}
\label{sec:low_rank}
In many numerical problems, the underlying structure lies in a subspace whose dimensionality is quite small compared to the full dimensionality of the solution space used in the numerical model. Examples include thermal problems, transport problems, elasticity and fluid flows among many others~\cite{Quarteroni2015,Hesthaven2016}. Low-rank methods can be used to reduce the complexity of such problems. The idea behind low-rank ROMs is to split the cost of computation into two stages. 

\begin{itemize}
    \item The \emph{offline} stage, where most of the computational complexity is resolved, consists of computing a set of solutions in the parametric space i.e.\ the training set, and then computing the reduced subspace that defines the span of training set. The idea is that space spanned by the training set snapshots is a good approximation to the true solution subspace. This is usually done by computing a reduced basis (RB) corresponding to the span of the training set. This step is typically performed only once to create the reduced model. 
    \item The \emph{online} step is performed every time a new query is submitted, in which the on-demand solution is computed using the reduced model which is cheaper to evaluate than the high-fidelity model.
\end{itemize}

Using notation from the previous section, let the high-fidelity (FE) solution space be denoted by $\PrimalSpace^h$ and the reduced subspace be denoted by $\PrimalSpace^r = \text{span}(\PrimalRB)$, where $\PrimalRB$ denotes the corresponding RB\@. In online stage, solutions in the subspace $\PrimalSpace^r$ are sought to approximate the high-fidelity solution. In discrete sense, $\PrimalRB$ is a matrix with each column corresponding to a basis vector. Then, the approximation in the subspace can be expressed as:
\begin{align*}
	& & \bm{u} \approx \PrimalRB \widehat{\bm{u}} & & \widehat{\bm{u}}  \in \mathbb{R}^{\nrdofprimal} & &
\end{align*}

\noindent All reduced degrees of freedom and reduced operators will be indicated by a hat \ $\widehat{\cdot}$    

For constructing the RB $\PrimalRB$ in offline stage, a finite subset of the parametric space $\mathcal{P}_\texttt{tr} \in \mathcal{P}$ is explored and an approximate subspace $\PrimalSpace^r$ is extracted from the set of high-fidelity snapshots from the training set. The number of degrees of freedom associated to the reconstruction problem in online stage is same as the cardinality of basis $\PrimalRB$, and it influences the efficiency of online phase. This is the key to low-rank methods, as it is based on the assumption that the intrinsic structure of the system behaviour lies in a low-rank subspace. If this assumption does not hold, it may be difficult to generate a RB that approximates well the system behaviour in the entire parametric space.

Next, two well known methods of finding the low-rank subspace: the Proper Orthogonal Decomposition and the Reduced Basis Methods are described briefly, as these methods form the backbone of many reduced models~\cite{Quarteroni2015,Hesthaven2016}. This is followed by the necessity and tools for efficient computation of terms that are non-linear w.r.t.\ the parameter and/or the solution. 

\subsubsection{Proper Orthogonal Decomposition}
In Proper Orthogonal Decomposition (POD) framework, high-fidelity snapshots are generated over the training set $\mathcal{P}_{\texttt{tr}}$ are collected and arranged as columns of a matrix. Let $\mathbf{S}_{\texttt{tr}} = \{\bm{u}_s \}_{s=1}^{N_S} \in \mathbb{R}^{\ndofprimal \times N_S}$ be a matrix containing the training set snapshots. Here, $N_S$ is the number of snapshots and $\ndofprimal$ is the number of finite element dofs. Typically, $N_S \ll \ndofprimal$ for two reasons: first is that the snapshot computation can be expensive and second is the fact that discovery of a low rank subspace only needs a relatively small number of snapshots  well distributed in the parametric domain. The RB can generated using the Singular Value Decomposition (SVD) or the Principal Component Analysis (PCA) methods. The SVD can be expressed as:
\begin{align*}
    \mathbf{S}_{\texttt{tr}} = \PrimalRB \mathbf{\Sigma} \mathbf{\Psi}^T
\end{align*}

\noindent where $\PrimalRB \in \mathbb{R}^{\ndofprimal \times N_S} $ and $\mathbf{\Psi} \in \mathbb{R}^{N_S \times N_S}$ are unitary matrices containing orthonormal vectors. $\mathbf{\Sigma} \in \mathbb{R}^{N_S \times N_S}$ is a diagonal matrix containing singular values. 

The decomposition given by the SVD has many interesting properties. The columns of $\PrimalRB$ form a orthonormal basis for column-space of $\mathbf{S}_{\texttt{tr}}$ and therefore can be used to approximate the solution subspace. Moreover, the columns of $\PrimalRB$, called as singular vectors, are arranged in the decreasing order of contribution to the total energy in $\mathbf{S}_{\texttt{tr}}$. Their contribution is indicated by the singular values in $\mathbf{\Sigma}$, which can be used to truncate the basis $\PrimalRB$. In fact, the singular values $ \lbrace \sigma_i \rbrace_{i=1}^{N_S}$ are arranged in decreasing order, where $\sigma_i$ is the $i$-th diagonal value of $\mathbf{\Sigma}$, exhibit many interesting properties in relation to the matrix $\mathbf{S}_{\texttt{tr}}$, some of which are listed below:

\begin{itemize}
    \item The square of singular values sum up to the total energy in the snapshot matrix, i.e. $\norm{\mathbf{S}_{\texttt{tr}}}{F} = \sum_{i=1}^{N_s} \sigma_i^2$
    \item The energy contribution by first $n$ singular vectors to the snapshot matrix is given by $\sum_{i=1}^{n} \sigma_i^2$.
    \item If $\PrimalRB^{'}$ and $\mathbf{\Psi}^{'}$ contain the first $n$ columns of $\PrimalRB$ and $\mathbf{\Psi}$, and also $\mathbf{\Sigma}^{'} \in \mathbb{R}^{n \times n}$ contains the first $n$ singular vectors $ \lbrace \sigma_i \rbrace_{i=1}^{n}$ in its diagonal; then the approximation
        \begin{align*}
            \mathbf{S}_{\texttt{tr}} \approx \PrimalRB^{'} \mathbf{\Sigma}^{'} \mathbf{\Psi}^{'T}
        \end{align*}
        is the best $n$-rank approximation in the sense of Frobenius-norm, with an approximation error $\sum_{i=n+1}^{N_S} \sigma_i^2$. A related property is that  $\PrimalRB^{'}$ and $\mathbf{\Psi}^{'}$ form the best $n$-rank subspaces of the column and row space, respectively, of the matrix $\mathbf{S}_{\texttt{tr}}$ in sense of Frobenius-norm. 
\end{itemize}

\noindent The above properties can be used for truncation of the basis $\PrimalRB$ by applying a tolerance $\delta$ to the fraction of energy contribution of first $n$ singular vectors. The first $n$ singular vectors of $\PrimalRB$ that contribute more $1-\delta$ fraction of energy are kept. This number $n$ can be calculated using:
\begin{align*}
    \min_n \ 1- \frac{\sum_{i=1}^{n} \sigma_i^2}{\sum_{i=1}^{N_s} \sigma_i^2}  < \delta
\end{align*}

The truncation leads to a subspace of the full solution space, whose rank is a function of tolerance $\delta$. This allows us to regulate the dimensionality but at the cost of losing diminished accuracy. Generation of the reduced model using the reduced basis $\PrimalRB$ will be shown in \cref{sec:rom_contact_lowrank}. 



\subsubsection{Reduced Basis Methods}
Even though POD methods offer a robust technique to compute the RB, they do not offer a strategy to sample the parametric space. On the other hand, the Reduced Basis Methods (RBM) are usually equipped with error estimators that allow greedy sampling of the parametric space~\cite{Prudhomme2001,Prudhomme2002}. The idea of greedy methods is to arrive at the objective by using the locally best way forward at each step, which is typically done by solving a minimization problem. An error estimator $\Delta(\mu)$ is an analytical error bound for the true reconstruction error as follows:
\begin{align*}
    \norm{\bm{u}(\mu) - \PrimalRB \widehat{\bm{u}}(\mu)}{2} \leq \Delta(\mu) \quad \forall \mu \in \mathcal{P}
\end{align*}
The error estimator is computed using the discrete residual $\bm{r}_h)$ of the high-fidelity system. For an unconstrained discrete problem with the residual:
\begin{align*}
    \bm{r}_h(\bm{u}, \mu) = \bm{f}(\mu) - \mathbf{K}(\mu) \bm{u}
\end{align*}
the error bound $\Delta(\mu)$ can be expressed as:
\begin{align*}
    \Delta(\mu) &= \norm{\mathbf{K}^{-1}(\mu) \bm{r}_h (\bm{u}, \PrimalRB \widehat{\bm{u}}(\mu))}{2} \\
                &\leq \norm{\mathbf{K}^{-1}(\mu)}{2} \norm{\bm{r}_h (\bm{u}, \PrimalRB \widehat{\bm{u}}(\mu))}{2} \\
                &= \frac{1}{\sigma_{min}(\mathbf{K}(\mu))} \norm{\bm{r}_h (\bm{u}, \PrimalRB \widehat{\bm{u}}(\mu))}{2}
\end{align*}
Further details pertaining to efficient computation of the residual norm and the least singular value can be found in~\cite[Chapter~3]{Quarteroni2015}. The error estimator has negligible cost of evaluation compared to the  true error, as the true numerical solution $\bm{u}$ is unknown. Therefore, at given greedy level $k$, the estimated best value of the parameter to be added is given by $\mu^{k+1} = \argmax \Delta(\mu)$ and is cheaper to compute. Thus, the training step is enriched with a new snapshot $\bm{u}(\mu^{k+1})$ and the RB $\PrimalRB$ is updated using methods like Gram-Schmidt process. This avoids the computation of the basis from scratch at each step $k$. The process can be stopped when $\max \Delta(\mu)$ is under a certain tolerance. The RBMs are advantageous in terms of reliability and efficiency of the reduction process, though the performance of the process, of course, depends closeness of the error estimator with the true error (see~\cite[Chapter~7]{Quarteroni2015}).

\subsubsection{Computation of non-linear terms}
\label{sec:affine_form}
Operators such as $\mathbf{K}(\mu)$ and $\bm{f}(\mu)$ can be a non-linear function of $\mu$ and therefore, their computation in each iterative step of a linearized reduced problem can prove to be expensive. Despite the introduction of low dimensional spaces $\PrimalSpace^r$, the complexity of computing these operators depends on the number of dofs in high-fidelity model. Commonly used methods to efficiently compute non-linear terms include the computation of affine decomposition using the Empirical Interpolation Method (EIM) and the computation of an reduced integration domain (RID) using hyper reduction methods.

The affine decompositions of the non-linear terms can be written as:
\begin{align*}
    \mathbf{K}(\mu) &\approx \sum_{i=1}^{N_K} \alpha^K_i(\mu) \mathbf{K}_i \\
    \bm{f}(\mu) &\approx \sum_{i=1}^{N_f} \alpha^f_i(\mu) \bm{f}_i
\end{align*}
\noindent where the scalar values $\big\{ \alpha^K_i \big\}$ and $\big\{ \alpha^f_i \big\}$ are the only quantities dependent on $\mu$. Therefore, the cost of online computation for each query of $\mu$ is reduced to the cost of computing these scalar values. 

It is possible that simple dependencies on $\mu$ admit affine decompositions that are easier to compute (see elastic energy computation in~\cite{Benaceur2020}). However, this is not the case in general and also there are many cases where operators show dependence on the solution $\bm{u}$ itself, for e.g.\ large deformation problems. In such cases, an affine decomposition is even more valuable as the resolution of the reduced problem will involve an iterative procedure demanding evaluation of non-linear terms in each iterative step. Usually, the EIM framework is used for computing the affine decomposition in such cases~\cite{Barrault2004}. EIM based affine decomposition is also equipped with so-called magic points in the non-parametric domain. Evaluation of the non-linear term only at the magic points is sufficient for computing the coefficients $\big\{ \alpha^K_i \big\}$ and $\big\{ \alpha^f_i \big\}$. The cost of computation is reduced drastically as the parameter independent modes $\mathbf{K}_i$ and $\bm{f}_i$ can be computed offline. For further details on computation of EIM based affine forms, see~\cite[Chapter~10]{Quarteroni2015}. 

Similar reduction in computation cost can also be achieved by hyper-reduction methods, where the non-linear terms are evaluated only a selected subset of the quadrature points in the entire integration domain. Recalling the property of standard Galerkin approach in finite elements, each row in matrix $\mathbf{K}(\mu)$ (and each entry in $\bm{f}(\mu)$) corresponds to integration over the local support of a particular node. The hyper-reduced evaluation of non-linear operators can be written as $\mathbf{P} \mathbf{K}(\mu)$ (and $\mathbf{P} \bm{f}(\mu)$). The rectangular matrix $\mathbf{P}$ each row containing only one non-zero entry equal to $1$, at a position corresponding to a node in the RID~\cite{Ryckelynck2009}. In other words, weak forms are computed only on supports surrounding the nodes in the RID. The RID has be constructed by selecting specific dofs or quadrature points in the entire domain. Various methods have been applied to this end, including free energy indicators~\cite{Ryckelynck2009}, optimized cubature methods~\cite{Hernandez2017,Lauzeral2019}, and Discrete Empirical Interpolation Method~\cite{Fauque2018}.

\section{ROMs for parametrized contact mechanics problems}\label{sec:rom_contact_literature}
ROMs have been applied extensively to unconstrained problems such as thermal and mechanics problems~\cite{Quarteroni2015,Hesthaven2016} and also to mixed problems with equality constraints such as the incompressibility condition in  Stokes and Navier-Stokes flow problems~\cite{Quarteroni2007}. Application of ROMs to  variational problems with inequality constraints (also referred to as variational inequalities) has been more recent~\cite{Haasdonk2012,Benaceur2020,Balajewicz2016,Bader2016,Burkovska2015,Fauque2018}. Inequality constraints appear in mechanical problems with obstacles or multi-body mechanical problems where there is a possibility of contact between bodies and obstacles, or with other bodies. Moreover, the region of contact is unknown a priori. Although the scope of this thesis is limited to contact mechanics, several other applications of variational inequalities are found in porous media flow problems~\cite{Oden1980}, cavitation problems in lubrication systems~\cite{Kinderlehrer2000}, anti-plane frictional problems~\cite{Sofonea2009} and even in financial trading problems~\cite{Burkovska2015}.

As seen in the \cref{sec:contact_formulations}, inequality constraint problems are often posed in mixed form with the KKT conditions that force the Lagrange Multipliers to be non-negative. Enforcing this condition in ROMs is non-trivial as traditional methods of computing RBs do not preserve the non-negative nature of the input information.

One of the first works on reducing the contact mechanics problem~\cite{Haasdonk2012}, proposed the idea of using the contact pressure snapshots directly to define a non-negative subcone. Both displacement and contact pressure snapshots were generated in a greedy fashion using error estimators. Fundamental aspects of the reduced problem, such as the existence and uniqueness of the solution and inf-sup stability were also explored.

Compression of snapshots to create a reduced basis was studied in~\cite{Balajewicz2016}, where Non-Negative Matrix Factorization (NNMF) was used to compute a basis with user-specified cardinality but does not provide any means to specify truncation tolerance for the reduced basis. Also, error estimators for greedy sampling of parametric space are developed.

Algorithms to sort a precomputed set of snapshots in order of importance to create a compact basis are explored using projection methods in~\cite{Burkovska2015} and~\cite{Benaceur2020}, using an Angle-Greedy and Cone-Projected Greedy (CPG) procedures respectively. The former does not take into account non-negative restrictions, whereas the latter uses a \emph{cone} projection involving non-negative coefficients. This makes CPG more efficient in capturing the contact pressure subcone.

The hyper-reduction approach in~\cite{Fauque2018} defines a subdomain of the contact problem containing the most important points. The reduction is achieved with the usage of a POD basis for displacement and resolution of the weak form on the reduced integration domain. Reconstruction of displacement solution on the full domain is relatively straightforward using the POD basis. On the other hand, the reconstruction of contact pressure involves solving a non-negative least square problem using the snapshots.

A Craig-Bampton based resolution of the contact problem was discussed in~\cite{Manvelyan2021}  where the reduction of the displacement field was achieved using the Krylov subspace method. The Lagrange Multiplier method was not reduced under the assumption that the number of contact dofs remains small.

All of these works are based on Lagrange Multiplier approach to solve the inequality constrained reduced problem. An exception to this trend is~\cite{Bader2016} where penalty approach is used.  Instead, EIM is used for efficient computation of barrier functions that are used to restrict the solution in the feasible domain. 

Most contributions build reduced models on problems with static contact pairs, meaning contact pairs do not change with the state of the system. However,~\cite{Benaceur2020} considers dynamic contact pairs, with node-to-segment formulation. To efficiently evaluate non-linearities due to dynamic pairing, they use Empirical Interpolation Method to define an affine decomposition of distance functions.

\subsection{Low-rank approach for contact problems}
\label{sec:rom_contact_lowrank}
As seen above, many of the ROM approach to contact mechanics are based on the low-rank approach. As the validity of low-rank hypothesis will be discussed in detail in \cref{ch:lowrank}, a generic low-rank approach is detailed here, based on ideas of~\cite{Balajewicz2016,Benaceur2020}. The reduction of the displacement and contact pressure fields is discussed, followed by defining the reduced problem  using the reduced spaces for the two fields. As the contact problem is a mixed problem, the displacement and contact pressure field will also be referred to as \emph{primal} and \emph{dual} fields \\[-0.5em]

\noindent \textbf{Reduction of displacement}: The displacement field $\bm{u}$ for a contact mechanics does not need any special treatment. Any methods that can be applied for reducing a mechanics problem without contacts can also be used for the displacement. Left singular vectors of displacement snapshots were used in many references in ROMs of contact mechanics. The same will be the case in this thesis unless specified, and the primal RB will be derived using POD\@.\\[-0.5em]

\noindent  \textbf{Reduction of contact pressure}: Computing the dual basis is more complicated than the primal basis. The dual field i.e.\ the Lagrange multiplier must satisfy the non-negativity constraints. In case of the contact problems without adhesive and cohesive surfaces, this constraint admits a physical meaning as negative contact pressure cannot be admitted, as discussed in \cref{sec:mixed_form_FE}. An orthogonal basis generated using projection based methods cannot satisfy such constraints, as it is surely bound to contain negative entries. To ensure such constraints, one way is to define a subcone for the dual field instead of a subspace. A function cone, unlike a function space, is spanned by a set of non-negative basis functions and non-negative coefficients. In fact, $\DualSpace$ and $\DualSpace^h$ in \cref{sec:contact_discretization} are also cones. By extension, a subcone, unlike a subspace, must be equipped by a non-negative RB and must be spanned by non-negative coefficients. 
\begin{align*}
	\DualSpace^r = \text{span}^{+}(\DualRB) = \sum_{i=1}^{N} \widehat{\lambda}_i \DualRBF_i  \ \ , \ \ \ \widehat{\lambda}_i \geq 0
\end{align*}
where $\DualRB$ and $\{\DualRBF_i\}_{i=1}^{N}$ are the dual RB and dual RB functions. $\widehat{\lambda}_i$s are the dual reduced dofs.

\noindent Once the reduced subspace $\PrimalSpace^r$ and subcone $\DualSpace^r$ are available, the reduced problem can be generated from~\eqref{eq:kkt_weak_form}, replacing the continuous spaces $\PrimalSpace$ and $\DualSpace$ with the reduced spaces $\PrimalSpace^r$ and $\DualSpace^r$. The reduced KKT conditions on the solution $(\widehat{\bm{u}}, \widehat{\bm{\lambda}}) \in \mathbb{R}^{\nrdofprimal} \times \mathbb{R}^{\nrdofprimal}$ can be expressed as:

\begin{subequations}
	\label{eq:kkt_reduced}
	\begin{align}
        \widehat{\mathbf{K}}(\mu) \widehat{\bm{v}} - \widehat{\bm{v}}^T \widehat{\bm{f}}(\mu) + \widehat{\mathbf{C}}^T(\mu,\widehat{\bm{v}}) \widehat{\bm{\lambda}} &= \bm{0}  \\
        \widehat{\bm{\lambda}} & \geq \bm{0} \label{eq:kkt_reduced_nonnegativity}  \\
        \widehat{\mathbf{C}}(\mu,\widehat{\bm{u}}) \widehat{\bm{u}} - \widehat{\bm{g}}(\mu,\widehat{\bm{u}}) & \leq \bm{0} \label{eq:kkt_reduced_constraint} \\
        \widehat{\bm{\lambda}}^T (\widehat{\mathbf{C}}(\mu,\widehat{\bm{u}}) \widehat{\bm{u}} - \widehat{\bm{g}}(\mu,\widehat{\bm{u}})) & = 0 \label{eq:kkt_reduced_comp_slack}
	\end{align}
\end{subequations}
where the discrete reduced operators are built by introducing the reduced basis functions $\PrimalRBF$ and $\DualRBF$
\begin{alignat*}{3}
	&\widehat{\mathbf{K}}(\mu)_{ij} &&= \ a(\mu;\PrimalRBF_i,\PrimalRBF_j) \\
	&\widehat{\mathbf{C}}(\mu,\widehat{\bm{v}})_{ij} &&= \ b(\mu,\bm{v}^r;\DualRBF_i, \PrimalRBF_j) \\
	&\widehat{\bm{f}}(\mu)_i &&= \ f(\mu;\PrimalRBF_i) \\
    &\widehat{\bm{g}}(\mu,\widehat{\bm{v}})_i &&= \ d(\mu,\bm{v}^r;\DualRBF_i) 
    \intertext{with}
    &\bm{v}^r = \PrimalRB \widehat{\bm{v}} &&
\end{alignat*}
\noindent It is evident that the reduced form~\eqref{eq:kkt_reduced} has the same structure as the discrete form generated using finite elements in~\eqref{eq:kkt_discrete}, as both equations are built using a Galerkin formulation and by introducing discrete solution spaces, albeit of different dimensionality.

Non-negative nature of the dual RB $\DualRB$ forbids the use of orthogonal decompositions of Lagrange Multiplier snapshot matrix $\mathbf{\Lambda}$. To this end, non-negativity preserving decompositions have been explored by various authors. Non-negative matrix factorization (NNMF) method~\cite{Lee1999} can be used to decompose a non-negative matrix, such as $\mathbf{\Lambda}$ into two low-rank non-negative matrices $\mathbf{W},\ \mathbf{H}$, such that $\mathbf{\Lambda} \approx \mathbf{W} \mathbf{H}$. The left-hand matrix $\mathbf{W}$ is used as the dual RB in~\cite{Balajewicz2016}.

Another possibility is to directly use the snapshot vectors as the basis vectors, instead of computing the NNMF. This was the approach of~\cite{Haasdonk2012}, where the possibility of non-unique dual solution is also discussed, as snapshot vectors are not guaranteed to be linearly independent. Also, the number of dual dofs increases with the number of snapshots, preventing optimal reduction of the system. A greedy snapshot selection method based on the criteria of maximizing the volume of the reduced cone, namely the Cone-Projected Greedy (CPG) algorithm, was proposed in~\cite{Benaceur2020}. This algorithm creates a dual RB by greedily selecting snapshots from the snapshot matrix, attempting to create a more compact basis in comparison to using the full snapshot matrix. The algorithm is based on projection $\Pi_\DualRB$ of a vector $\bm{\lambda}$ on a vector cone $\DualRB$,
\begin{align}
	\Pi_\DualRB(\bm{\lambda}) :\approx \DualRB \bm{\alpha} \ , \ \ \text{where } \bm{\alpha} \text{ is } \argmin_\gamma \norm{\bm{\lambda}- \DualRB \bm{\gamma}}{} \ \forall \bm{\gamma} \geq  \bm{0}
	\label{eq:cp}
\end{align}
The greedy algorithm evaluates the error between each snapshot and its projection on the cone spanned by previously selected snapshots, and then adds that snapshot with maximum cone projection error. The process continues until the cone projection error is within a set tolerance.~\cite{Burkovska2015} had also proposed a similar algorithm called Angle-Greedy algorithm, but the projection error is calculated based on computation of the angle between the candidate snapshot vector and the space, and not the cone, spanned by the previously selected snapshots i.e.\ it does not place non-negativity constrain on the coefficients $\gamma$ in~\eqref{eq:cp}. This is not the best way of selecting snapshots for variational inequality problems since solutions must be sought in the reduced cone, and not the entire reduced space.

\begin{remark}[Constraint satisfying primal subspace in flow problems]{}{}\label{rem:div_free}
    If we shift our gaze to other types of constrained ROMs, an interesting treatment of the mixed problems appears in incompressible flow problems. In such problems, the velocity snapshots are divergence-free and the divergence operator in the mass conservation equation is linear; offering the possibility to compute a divergence-free subspace for the velocity field. Consequently, the reduced incompressible flow problem is as an unconstrained problem, as all candidates in velocity subspace naturally satisfy the constraints~\cite{Veroy2005,Liberge2010}. Apart from naturally satisfying the constraints, some elegant properties of vector algebra\footnote{The product of a divergence-free velocity field and the gradient of pressure vanishes due to the orthogonality of the solenoidal and irrotational vector fields} result in disappearance of bilinear term that includes pressure.

    Unfortunately, it is not straightforward to extend the same idea to contact mechanics problems (or variational inequalities), i.e.\ to compute a primal subspace that satisfies the non-penetration (or the inequality) constraints. This is because linear combinations of non-penetrating snapshots do not satisfy the non-penetration condition in general. As a consequence, the inequality constraints need be enforced explicitly while solving the reduced problem, in general.
\end{remark}

\noindent \textbf{Construction of non-linear operators:} For an efficient ROM, the construction of nonlinear operators must also be cheap. The construction of operators $\mathbf{C}(\mu,\widehat{\bm{u}})$ and $\bm{g}(\mu,\widehat{\bm{u}})$ in each iterative step can prove expensive. Affine decompositions of the nonlinear constraint operators using the EIM are discussed in~\cite{Benaceur2020}. The affine form was achieved by computing the EIM based interpolation of the distance functions $k(\mu, \bm{v}, \cdot)(\bm{x})$ and $g(\mu, \bm{v})(\bm{x})$ (recall \cref{not:strong_form}):
\begin{align*}
    k(\mu,\bm{v}; \cdot)(\bm{x}) &\approx \sum_{p=1}^{N_k} \alpha^k_p(\mu,\bm{v}) \kappa_p(\cdot)(\bm{x})\\
    g(\mu, \bm{v})(\bm{x}) &\approx \sum_{q=1}^{N_g} \alpha^g_q(\mu,\bm{v}) \gamma_q(\bm{x})
\end{align*}
where the interpolation is exact at the set of magic points $\{\bm{x}^k_i\}$ and $\{\bm{x}^g_i\}$. Note that $\kappa_p$ is linear w.r.t.\ its argument which can be any member in $\PrimalSpace$. Therefore, the coefficients can be computed so that they satisfy the following conditions:
\begin{gather}
    \begin{aligned}
        k(\mu,\bm{v}; \cdot)(\bm{x}^k_m) &= \sum_{i=1}^{N_k} \alpha^k_p(\mu,\bm{v}) \kappa_p(\cdot)(\bm{x}^k_m) \ &\forall m=\{1,2,\dots N_k\}\\
        g(\mu, \bm{v})(\bm{x}^g_n) &= \sum_{q=1}^{N_g} \alpha^g_q(\mu,\bm{v}) \gamma_q(\bm{x}^k_n) \ &\forall n=\{1,2,\dots N_g\}
    \end{aligned}
    \label{eq:dist_funcs_eim_interp}
\end{gather}
Therefore, the distance function needs to be evaluated at most $N_k + N_g$ points, instead of the whole contact surface. Also, $\kappa_i(\cdot)$ and $\gamma_i$ are computed just once using the training set and the corresponding discrete operators can also be built during the offline stage using the definition of weak forms in~\eqref{eq:dist_func_weak}.
\begin{gather}
    \begin{aligned}
        \widehat{\mathbf{C}}^p_{ij} &= \int_{\Gamma^2} \DualRBF_i \kappa_p(\PrimalRBF_j)(\bm{x})  \ \partial \Gamma  \\
        \widehat{\bm{g}}^q_i &= \int_{\Gamma^2} \DualRBF_i \gamma_q(\bm{x})  \ \partial \Gamma
    \end{aligned}
    \label{eq:dist_funcs_eim_operators}
\end{gather}
where $b_p$ and $d_q$ correspond to the weak forms of individual EIM modes $\kappa_p$ and $\gamma_q$ that are independent of the state of the system. Thus, during the online stage, the non-linear operators can be approximated as:
\begin{align*}
    \widehat{\mathbf{C}}(\mu,\bm{v}) &= \sum_{p=1}^{N_k} \alpha^k_p(\mu,\bm{v}) \widehat{\mathbf{C}}^p \\
    \widehat{\bm{g}}(\mu,\bm{v}) &= \sum_{q=1}^{N_g} \alpha^g_q(\mu,\bm{v}) \widehat{\bm{g}}^q 
\end{align*}
where only $\alpha^k_p$ and $\alpha^g_q$ are evaluated efficiently in online stage using~\eqref{eq:dist_funcs_eim_interp}.

Another approach for efficient computation of non-linear operators is the hyper-reduction approach~\cite{Fauque2018}. In this approach, a small subset of the full domain $\Omega_A \subset \Omega_1 \cup \Omega_2$ is computed that serves as the reduced integration domain (RID). The RID $\Omega_A$ consists of the most ``important'' integration points computed by performing EIM on the displacement RB. Then, integration of the operators $\widehat{\mathbf{C}}(\mu, \bm{v})$ and $\widehat{\bm{g}}(\mu, \bm{v})$ is carried out on a subset of the potential contact surface, computed using the intersection of the RID and the contact surface $\Omega_A \cap \Gamma_2$. Thus, the distance functions are computed only on a small subset of the contact surface.

\section{Scope of the thesis}
The current literature consists of various methodologies to reduce the contact mechanics problem, especially to the Lagrange Multiplier approach. In this thesis, the following contributions are made:
\begin{itemize}
    \item As seen in this chapter, the reducibility of contact problems has not been sufficiently investigated in current literature. This is particularly important as contact mechanics problems exhibit a local nature in terms of the contact zone. As the region of contact is dynamic for many contact problems, significant changes in the contact zone can prove consequential to reducibility. Therefore, in \cref{ch:lowrank}, we explore the validity of so-called \emph{low-rank hypothesis} and robustness of the reduced models using various metrics. The chapter concludes that the contact pressure field lacks linear separability, an essential feature for reducibility.
    \item The lack of low-rank structure implies that the offline phase might be expensive as small training sets will be unable to explore the underlying structure sufficiently. Even if the user is willing to compute large training sets, it is still not possible to compute a robust low-rank subspace for the contact pressure. In \cref{ch:sparse}, we focus on algorithms to handle dictionaries containing a large number of snapshots using ideas from sparse regression methods.
    \item As the linear inseparability of contact pressure field is established, we explore the idea of non-linear interpolations between contact pressure snapshots in \cref{ch:nlInterp}. Interpolation is done in a warped space generated using Dynamic Time Warping algorithm. We show that reconstructions of contact pressure can be more reliably performed with smaller training sets, demonstrating the potential of non-linear dimensionality reduction methods for contact mechanics problems.
\end{itemize}

\chapter{Limitations of low-rank approach to contact problems}\label{ch:lowrank}

Strategies for reducing complex numerical models to facilitate building parametric models requires the acquisition of underlying physics using precomputed data. A brief overview of these methods was given in \cref{sec:rom_literature}. Commonly used ROM methods are based on the paradigm of the so-called \textit{linear dimensionality reduction}. The foundation of this framework rests on the low-rank hypothesis, which states that for a high-dimensional quantity of interest, often a representative low-dimensional space exists and can be extracted using using linear transformations. The successful applications of ROM methods like POD, based on linear dimensionality reduction methods like PCA, to many parametric models can be attributed to the underlying low-dimensional behaviour of these systems.

In \cref{sec:rom_literature,sec:rom_contact_literature}, the wide variety of low-rank applications were mentioned, including constrained problems such as incompressible flows and contact mechanics problems. Reduced models for contact mechanics are important because the inequality constraints are geometrically non-linear, thereby introducing additional computational challenges (see \crefrange{sec:contact_formulations}{sec:contact_discretization}). However, the assumption of low-rank behaviour has implications on reduced models of contact mechanics problems, which are a class of variational inequality problems.

This chapter focuses on applicability of the low-rank hypothesis to the variational inequality problems, specifically contact mechanics problems. The effect of this hypothesis on various contact mechanics problems are demonstrated. A generic methodology of low-rank approach in contact mechanics problems was discussed in \cref{sec:rom_contact_literature}. In this chapter, the limitation of low-rank approach will be demonstrated on reduced models of various contact problems. It will be shown that the contact pressure field is linearly inseparable and therefore does not satisfy low-rank hypothesis. To endorse this argument quantitatively, the so-called \textit{validation metrics} will be introduced to used to endorse the arguments about linear inseparability using specific numerical examples. Finally, concluding arguments about linear inseparability and a few perspectives on circumventing the limitations due to lack of low-rank structure are given.

\section{Low-rank contact problem using active set method}
The non-linear reduced problem in~\eqref{eq:kkt_reduced} is based on the low-rank hypothesis, as it introduces the low-rank primal and dual bases. This approach is similar to~\cite{Haasdonk2012,Balajewicz2016,Benaceur2020}. Two sources of non-linearities exist in this problem: the resolution of active/inactive constraints and the dynamic contact pairs. Irrespective of the linearization algorithm utilized, both non-linearities depend on the current estimate of the displacement field at a given linearization step. 

Here,~\eqref{eq:kkt_reduced} is solved using the fixed-point iterations and the active-set method applied to the reduced constraints. The linearized problem at a fixed-point iteration level $p$ can be expressed as follows:

\begin{align}
	\label{eq:contact_red_fixedpoint}
	\begin{bmatrix}
		\widehat{\mathbf{K}} &  \widehat{\mathbf{C}}_{\ActiveSet}^T(\mu, \widehat{\bm{u}}^p) \\
        \widehat{\mathbf{C}}_{\ActiveSet}(\mu,  \widehat{\bm{u}}^p)  & \mathbf{0}
	\end{bmatrix}  
	\begin{bmatrix}
		\widehat{\bm{u}}^{p+1} \\
		\widehat{\bm{\lambda}}^{p+1}_{\ActiveSet}
	\end{bmatrix} =  \begin{bmatrix}
		\widehat{\bm{f}} \\
		\widehat{\bm{g}}_{\ActiveSet}(\mu, \bm{u}^p)
	\end{bmatrix}
\end{align}

\noindent The subscript $A$ indicates the active set of reduced contact constraints, i.e.\ the constraints that satisfy the equality and force the solution to lie on the boundary of the feasible region. Note that~\eqref{eq:contact_red_fixedpoint} is similar to the fixed-point system in~\eqref{eq:contact_fixedpoint}. The fixed point algorithm stated in \cref{alg:online_phase}, was implemented in Python using NumPy package~\cite{numpy} for vectorial operations. \\

\begin{algorithm}[!htb]
	\caption{Online phase}\label{alg:online_phase}
	\begin{algorithmic}[1]
		\State{Input: Queried value of parameter $\mu$}
		\State{Given: Primal basis $\PrimalRB$ and dual basis $\DualRB$}
		\Statex{\hspace{3em} Reduced operators $\widehat{\mathbf{K}}, \widehat{\bm{f}}$ (possible to build offline)}
		\State{Initiate boolean array $\ActiveSet$  with one random element set to True.}
        \While{$\widehat{\bm{u}}$ not converge}
			\State{Build constraint operators $\mathbf{C}(\mu,\widehat{\bm{u}}^p), \ \bm{g}(\mu,\widehat{\bm{u}}^p)$ using FEM} \label{step:operators}
			\State{Project constraint operators on RBs} 
			\Statex{\centerline{$ \begin{array}{ll}
						\widehat{\mathbf{C}}(\mu,  \widehat{\bm{u}}^p) &= \DualRB^T{\mathbf{C}}(\mu,  \widehat{\bm{u}}^p) \PrimalRB  \\
						\widehat{\bm{g}}(\mu,  \widehat{\bm{u}}^p) &= \DualRB^T{\bm{g}}(\mu,  \widehat{\bm{u}}^p) 
                    \end{array} $}}
			\State{Filter-out rows that are not in active set: } 
			\Statex{\centerline{$\begin{array}{ll}\widehat{\mathbf{C}}_{\ActiveSet}(\mu,  \widehat{\bm{u}}^p) &= \widehat{\mathbf{C}}(\mu,  \widehat{\bm{u}}^p)[\ActiveSet,:] \\
                    \widehat{\bm{g}}_{\ActiveSet}(\mu,  \widehat{\bm{u}}^p) &= \widehat{\bm{g}}(\mu,  \widehat{\bm{u}}^p)[\ActiveSet]\end{array}$} }
            \State{Solve system~\eqref{eq:contact_red_fixedpoint}}
            \State{Set $\widehat{\bm{\lambda}}^{p+1}[\ActiveSet] \gets \widehat{\bm{\lambda}}^{p+1}_{\ActiveSet}$ and $\widehat{\bm{\lambda}}^{p+1}[\ActiveSet^{\mathsf{c}}] \gets \bm{0}$}
            \State{Update active constraints set}
			\Statex{\hspace{6em} $\ActiveSet[i] = \left \lbrace \begin{array}{lll}
                                    \texttt{False} & \text{if} & \widehat{\bm{\lambda}}^{p+1}_i < 0 \\
                                    \texttt{True} & \text{if} & (\widehat{\mathbf{C}} \widehat{\bm{u}}^{p+1} - \widehat{\bm{g}} \geq 0)_i
                            \end{array} \right. $ } \Comment{violations of~\eqref{eq:kkt_reduced_nonnegativity} and~\eqref{eq:kkt_reduced_constraint}}
        \EndWhile{}
        \State{Reconstruct $\bm{u} = \PrimalRB \widehat{\bm{u}}$ and $\bm{\lambda} = \DualRB \widehat{\bm{\lambda}}$}
        \State{Output: $\bm{u}, \bm{\lambda}$}
	\end{algorithmic}
\end{algorithm}

\noindent \textbf{Construction of operators:} The implication of dynamic contact pairing between two bodies of a contact problem is that the non-linear operators that enforce the non-penetration constraint must be updated for each new estimate of displacement field. On the other hand, for generating an efficient ROM, the construction of nonlinear operators must be inexpensive. This is not the case in \cref{alg:online_phase}, where the full order operators $\mathbf{C}(\mu,\widehat{\bm{u}})$ and $\bm{g}(\mu,\widehat{\bm{u}})$ are constructed in each iteration (\cref{step:operators}). Affine decompositions of the nonlinear constraint operators using the EIM are discussed in~\cite{Benaceur2020}. EIM decomposition of the constraint operators further splits computational complexity  into offline and online stages, by selecting the so-called ``magic points'' that are a small subset of all potential contact nodes where the nonlinear terms are evaluated in the online phase. The hyper-reduction method also permits computation of operators using fewer points~\cite{Fauque2018}. This is crucial to solve problems where the contact pairs are strongly dependent on the displacement field. However, the efficient construction of constraint operators is not covered in this work, and the focus is maintained on the non-negative dual space and its low-rank approximation.

\section{Illustrative case}\label{sec:illustration}
In this section, the Type-1 problem of a 1D elastic rope-obstacle problem from~\cite{Haasdonk2012,Balajewicz2016,Bader2016} is considered. The geometry of the problem enjoys further simplification of static contact pairs. In other words, each point on the elastic rope can come in contact with a unique point on the obstacle, thereby stripping the distance functions $k$ and $g$ of their non-linear dependence on the displacement field $u$.

The parametrized model of the elastic rope with obstacle function can be expressed as:
\begin{equation}
	\begin{aligned}
		\nu \nabla^2 u(x) &= f && \text{ on }  x \in [0,1] \\
		u(0) = u(1) &= 0 & \\
		u(x) &\geq g(x,\gamma) && \text{ on } \gamma \in [-0.5,0.5]
	\end{aligned}
\end{equation}

\noindent where the parameterized obstacle function is defined as:
\begin{equation}
	\begin{aligned}
		g(x,\gamma) = -0.2 (\sin(\pi x) - \sin(3 \pi x)) -0.5 + \gamma x 
	\end{aligned}
	\label{eq:obstacle}
\end{equation}

\noindent and the quantities $\nu = 30$ and $f=250$ are independent of the parameter $\gamma$.

\begin{figure}[!htb]
	\centering
	\begin{subfigure}[t]{0.48\linewidth}
		\includegraphics[width=1\linewidth]{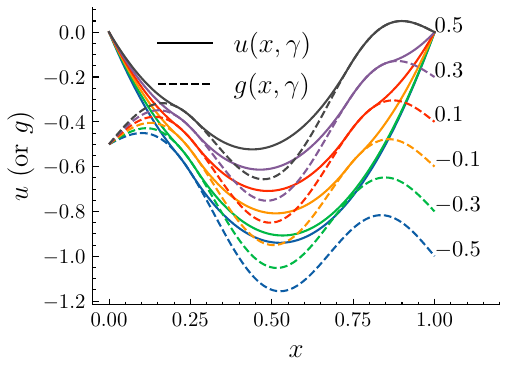}
		\caption{Deformation $u(x,\gamma)$ and obstacle $g(x,\gamma)$ snapshots}
		\label{fig:membrane_snaps_deformation}
	\end{subfigure}
	\begin{subfigure}[t]{0.48\linewidth}
		\includegraphics[width=1\linewidth]{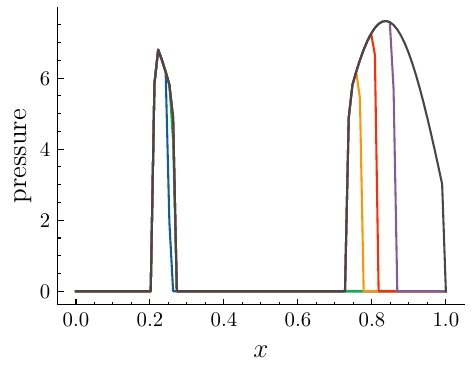}
		\caption{Contact pressure snapshots}
		\label{fig:membrane_snaps_cp} 
	\end{subfigure}
    \caption[Sample snapshots for elastic rope-obstacle problem]{Sample snapshots for elastic rope-obstacle problem with $\gamma \in \mathcal{P}_{\texttt{tr}} \subset [-0.5,0.5]$.}\label{fig:membrane_snaps}
\end{figure}

Snapshots are generated in the training set, $\mathcal{P}_\texttt{tr} \in \mathcal{P}$ consisting of 10 equidistant points in the parametric space $\mathcal{P}$ defined in~\eqref{eq:obstacle}. The deformation and contact pressure snapshots are shown in \cref{fig:membrane_snaps}. The primal RB is built using the POD approach and dual RB is created using the cone-projected greedy algorithm of~\cite{Benaceur2020}.

\begin{figure}[!htb]
	\centering
	\begin{subfigure}[t]{0.48\linewidth}
        \includegraphics[width=0.9\linewidth]{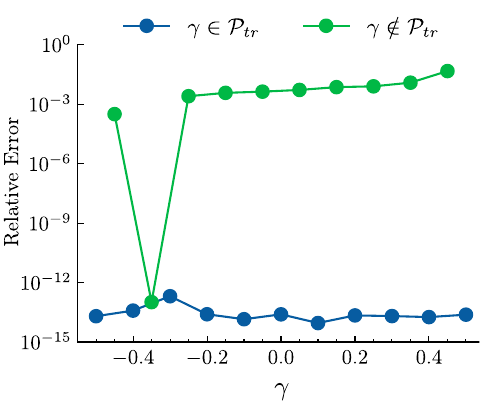}
		\caption{Full rank}
		\label{fig:membrane_recon_error_full_RB}
	\end{subfigure}
	\begin{subfigure}[t]{0.48\linewidth}
        \includegraphics[width=0.9\linewidth]{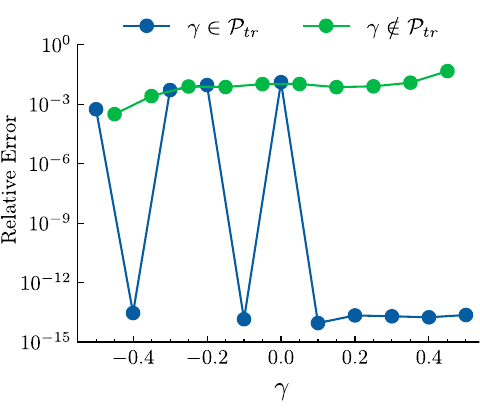}
		\caption{Full rank primal and truncated dual RB}
		\label{fig:membrane_recon_error_truncated_RB}
	\end{subfigure}
    \caption{Primal reconstruction ($\mathcal{H}^1$) errors for points in and outside the training set for the rope-obstacle problem.}\label{fig:membrane_recon_error}
\end{figure}

The reconstruction errors are shown in \cref{fig:membrane_recon_error} for the two cases of full and truncated dual RB\@. As expected, the points in the training set show very low reconstruction error, and the points outside the training set show a moderate error. The only exception at $\gamma=-0.35$ where reconstruction with full RB is highly accurate because the contact area does not change in the regime $\gamma \in [-0.4, -0.3]$ and the training set has two snapshots in this region. When dual RB is truncated, some points in the training set achieve the same error level as that of points outside the training set. These points correspond to the same dual snapshots that were eliminated by the truncation procedure of the CPG algorithm. Another relevant observation is that the number of active dofs of the reduced dual solution $\widehat{\bm{\lambda}}$ is small for all reconstruction cases. For points in the training set, $\widehat{\bm{\lambda}}$ has exactly 1 active dof, whereas for the points outside the training set, it has a maximum of 2 active dofs. This is due to the fact that the dual basis is composed primarily of snapshots and  does not undergo any compression like in the case of the primal basis generated by POD\@.

\begin{table}[!htb]
    \small
	\centering
	\begin{tabular}{llcccc}
		\toprule
		&& Primal basis & Dual basis     \\
		&& (POD) & (CP-greedy)   \\
		\midrule 
		\multirow{2}{*}{Full} & Tolerance & $\varepsilon$ & $\varepsilon $ \\
		\cdashline{2-4} 
		& Rank & 11 & 11  \\
		\midrule
		\multirow{2}{*}{Truncated} & Tolerance & $\varepsilon$ & $2\times10^{-1}$ \\
		\cdashline{2-4}
		& Rank & 11 & 7 \\
		\bottomrule
	\end{tabular}
    \caption[Truncation tolerances and ranks for both bases for the rope-obstacle problem.]{Truncation tolerances and ranks for both bases for the rope-obstacle problem. $\varepsilon$ indicates numerical precision}\label{tab:membrane_basis_ranks}
\end{table}

\begin{figure}[!htb]
	\centering
	\includegraphics[width=0.5\linewidth]{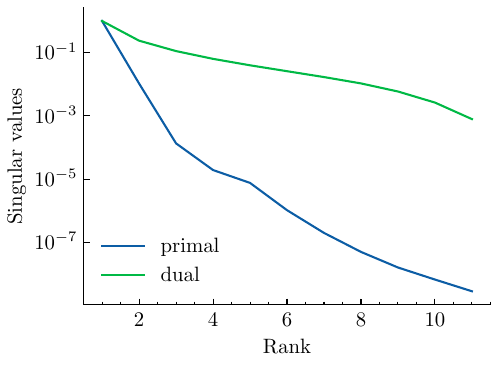}
    \caption[Decay of singular values for the training set of the elastic rope-obstacle problem.]{Decay of singular values (cumulative and normalized) for the training set of the elastic rope-obstacle problem. A similar figure is also given in Reference~\cite[Figure~4]{Fauque2018}}\label{fig:membrane_error_decay}
\end{figure}

In general, any efficiently generated reduced basis has its vectors arranged in order of importance. In case the basis consists of left singular vectors, generated by SVD, the decreasing importance of basis vectors is indicated by the singular values. This is a very useful indication of the rank of the subspace. Even though the dual basis is not computed using SVD, it can still be useful to compute singular values and assess their decay. This evolution of singular values, seen in \cref{fig:membrane_error_decay}, shows that the decay is slow compared to the primal variable. 

Also, the discussions in previous sections highlight the fact that dual basis $\DualRB$ cannot be truncated without a significant loss of accuracy. This is also evident from \cref{tab:membrane_basis_ranks}, even a truncation tolerance as high as $0.2$ leads to truncation of just four vectors in dual RB\@.

\begin{figure}[!htb]
	\centering
    \includegraphics[width=0.8\linewidth]{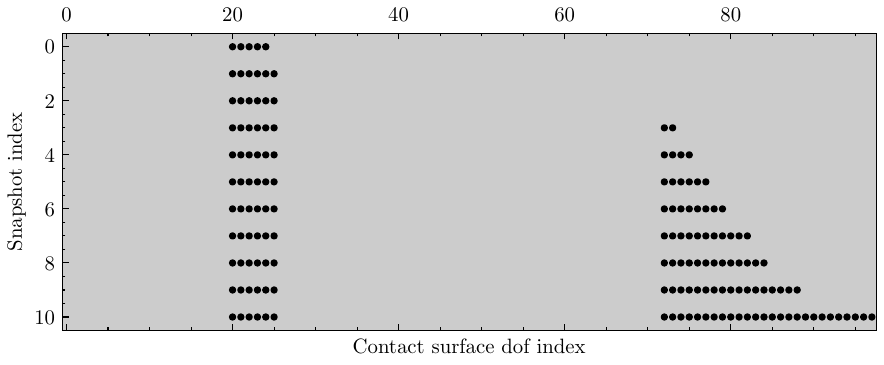}
    \caption[Sparsity pattern of dual snapshots for the rope-obstacle problem.]{Sparsity pattern of dual snapshots for the rope-obstacle problem. Each row corresponds to a snapshot}\label{fig:membrane_dual_spy}
\end{figure}

The above observations are explainable by the hypothesis that the dual snapshots are highly ``inseparable'' and do not lie in a low-dimensional subcone. Each contact pressure snapshot in \cref{fig:membrane_snaps} is linearly independent of the other due to difference in the contact area. Hence, if any snapshot is picked randomly from this figure, it cannot be reasonably approximated by any linear combination of the rest of the snapshots. In other words, each snapshot vanishes outside the contact zone, and this zone is almost ``unique'' for each snapshot. For a particular snapshot to be reasonably approximated, other snapshots with similar contact zones are needed. Therefore, the elimination of snapshots from the dual basis leads to a drastic increase in the reconstruction error for the corresponding parametric values in the training set. The dual variable, therefore, not only needs special treatment due to its positivity constraints but is also highly inseparable because of sensitivity to contact position and area. The inseparability is also evident in the non-zero pattern (``\texttt{spy}'' plot) of dual snapshots visualized in \cref{fig:membrane_dual_spy}, where most rows have an almost unique sparsity pattern. It is unlikely that a subspace whose members show varying sparsity patterns will admit a low-dimensional behaviour.

More numerical examples are explored to study the lack of low-rank behaviour in the \cref{sec:demos}.

\section{Metrics for low-rank methods}\label{sec:metrics}
In this section, we introduce the validation metrics~\cite{davies2002}, that provide a quantitative evaluation of the robustness of low-rank models. Validation metrics defined in~\cite{davies2002} identify three different measures for this purpose:
\begin{itemize}
    \item \emph{Compactness}: As the name suggests, this metric is a measure of dimensionality of the reduced space. Reference~\cite{davies2002} defines compactness as the squared sum of the first $m$ singular values of the snapshot matrix. This metric is similar to the decay of singular values studied in the previous section in \cref{fig:membrane_error_decay}. However, an alternative definition given by~\eqref{eq:compactness} is used here. The new definition allows extending the concept of compactness to CPG reduced bases, where equivalents of singular values are undefined. Also, in case of an orthogonal basis, this definition is equivalent to the normalized singular values used in \cref{fig:membrane_error_decay}.
\begin{align}
    C(m) :=  \frac{\norm{\mathbf{S}_\texttt{tr} - \Pi_{\Psi[:m]}(\mathbf{S}_\texttt{tr})}{F}}{\norm{\mathbf{S}_\texttt{tr}}{F}}
	\label{eq:compactness}
\end{align}
\noindent where $\mathbf{S}_{\texttt{tr}}$ is a matrix containing (either primal or dual) snapshots in $\mathcal{P}_{\texttt{tr}}$ arranged column-wise. $\Pi_{\Psi[:m]}(\mathbf{S}_{\texttt{tr}})$ is a projection operators that projects each column in $\mathbf{S}_{\texttt{tr}}$ on the first $m$ vectors of the basis $\Psi$ (indicated in Python notation by $\Psi[:m]$). The snapshot matrix $\mathbf{S}_{\texttt{tr}}$, projection operator $\Pi$ and the basis $\Psi$ are generic symbols and are used in the following combinations, indicated by the following labels: \\[0.5em]
    \begin{centering}
        \small
        \centering
        \begin{tabular}{lccc}
            \toprule
            Label in figures & Snapshots in $\mathbf{S}_{\texttt{tr}}$  & Projection operator $\Pi$  & Basis $\Psi$   \\
            \midrule 
            primal orth & Primal & Orthogonal & POD primal basis $\PrimalRB$    \\ 
            dual orth & Dual & Orthogonal & POD dual basis $\Upsilon$    \\ 
            dual cone & Dual & Cone & CPG dual basis $\DualRB$    \\
            \bottomrule
        \end{tabular}
    \end{centering} \\[0.5em]
\noindent The Frobenius norm is chosen in~\eqref{eq:compactness}, because it is a natural choice for a POD/SVD basis that contains the best rank-wise approximations in the sense of Frobenius norm, but other norms can be also be used. This is actually equivalent to the original definition of compactness in~\cite{davies2002}, because in the case $\Psi$ consists of the left singular vectors, compactness corresponds to cumulative energy of singular values.

For the dual field, the metric is given using two kinds of projections: orthogonal and CPG\@. Though the dual CPG metric is more relevant to the inequality constrained problems, the dual orthogonal metric is also included as it is more intuitive due to correspondence with the singular values.
    \item \emph{Generalization Ability}: This metric measures the ability of the reduced basis to approximate parametric instances that are not in the training set. It is computed by using a leave-one-out approach on the training set to generate a $m$-rank basis. Then, the reconstruction error of the eliminated snapshot is used to compute generalization ability as given by~\eqref{eq:generalization},
        \begin{align}
            G(m) = \frac{1}{N} \sum_{i=1}^N \frac{\norm{\bm{s}_i- \Pi_{\Psi(\mathbf{S}_{\texttt{tr}}\setminus \bm{s}_i)[:m]}(\bm{s}_i)}{}}{\norm{\bm{s}_i}{}}
            \label{eq:generalization}
        \end{align}
        where  $\Psi(\mathbf{S}_{\texttt{tr}}\setminus \bm{s}_i)$ is the basis created using the snapshot subset $\mathbf{S}_{\texttt{tr}} $ with $ \bm{s}_i$ removed (i.e.\ $\mathbf{S}_{\texttt{tr}}\setminus \bm{s}_i$). The symbols $\mathbf{S}_{\texttt{tr}}$, $\Pi$ and $\Psi$ follow the definition in~\eqref{eq:compactness}.
    \item \emph{Specificity}: The subspace spanned by the reduced basis is expected to contain elements similar to that of the training set. This metric measures the extent of dissimilarity between randomly picked elements from the subspace and the snapshots in the training set. To compute specificity, a set of random shapes in the associated subspace are generated using random coefficients $R = \{\Psi \bm{\alpha}_i\}_{i=1}^N$, where each vector $\bm{\alpha}_i$ is drawn randomly within a predefined range, and compared to the closest snapshot in the training set. The mean value of this error is defined as Specificity.
        \begin{align}
            S(m) = \frac{1}{N} \sum_{i=1}^N \min_{\bm{s}_j \in \mathbf{S}_{\texttt{tr}}}  \frac{\norm{\bm{s}_j - \Psi[:m] \bm{\alpha}_i}{}}{\norm{\bm{s}_j}{}}
        \end{align}
        The vectors $\{\bm{\alpha}_i \ | \ \bm{\alpha}_i \in \mathbb{R}^m\}$ are drawn randomly in the range defined by $\bm{\mu}_\alpha \pm \bm{\sigma}_\alpha \subset \mathbb{R}^m $. The vectors $\bm{\mu}_\alpha$ and $\bm{\sigma}_\alpha$ contain the element-wise mean and standard deviation of vectors $\{\bm{\gamma}_k\}_{k=1}^m$, where $\bm{\gamma}_k$ is the reduced coordinate of the training set snapshot $k$. (e.g.\ if $\Psi$ is an orthogonal basis, the vectors $\bm{\gamma}_k$ are the columns of matrix $\bm{\Gamma} := \Psi^T \mathbf{S}_{\texttt{tr}}$). In the specific case where the basis $\Psi$ is the $\DualRB$, the CPG basis, only non-negative entries are allowed in $\{\bm{\alpha}_i\}$. The symbols $\mathbf{S}_{\texttt{tr}}$, $\Psi$ follows the definition in~\eqref{eq:compactness}.
\end{itemize}
\section{Numerical Examples}\label{sec:demos}
\subsection{Hertz problem}\label{sec:hertz}
\begin{figure}[!htb]
	\centering
    \begin{subfigure}[t]{0.44\linewidth}
        \includegraphics[width=1\linewidth]{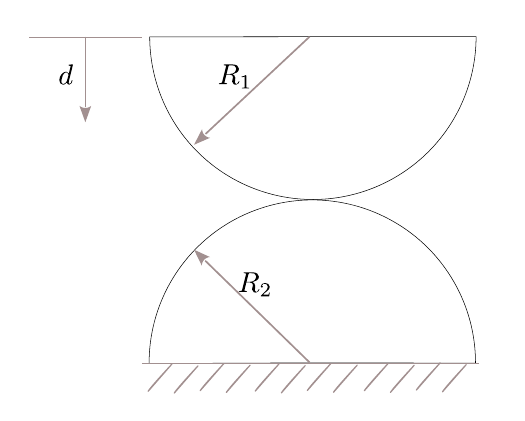}
        \caption{Hertz problem}\label{fig:hertz_half_cyl}
    \end{subfigure}
    \begin{subfigure}[t]{0.55\linewidth}
        \includegraphics[width=1\linewidth]{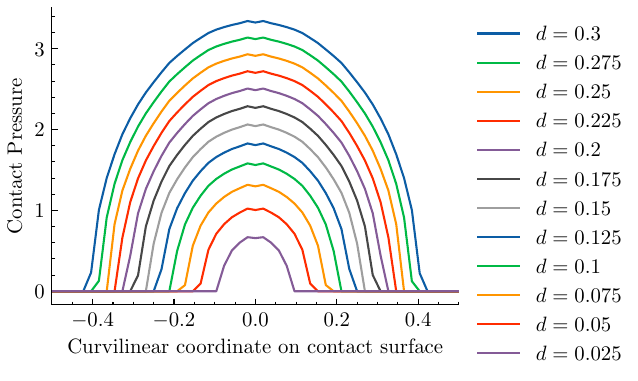}
        \caption{Contact pressure snapshots}\label{fig:hertz_snaps}
    \end{subfigure}
    \caption[Hertz problem: two half cylinders loaded against each other and the resulting contact pressure snapshots.]{Hertz problem: two half cylinders loaded against each other and the resulting contact pressure snapshots. Displacement $d$, imposed on the top cylinder, is treated as parameter in the reduced model.}
\end{figure}

This section deals with the reduced model of a Type-2 contact problem, namely the Hertz contact problem, of two half-cylinders loaded against each other, as shown in \cref{fig:hertz_half_cyl}. A \emph{loading parametrization} problem is solved, where the imposed displacement $d \in (0,0.3)$ is the parametric space, with $R_1=R_2=1.0$. The contact pressure snapshots for various imposed displacement values are shown in  \cref{fig:hertz_snaps}. \\

\noindent \textbf{High-fidelity model}:  A finite element model of the geometry in \cref{fig:hertz_half_cyl} is created using a quad mesh. Each half-cylinder has 513 nodes and 466 elements of which 78 elements lie of the potential contact surface i.e.\ the semi-circular edges. Linear shape functions are used to discretize the displacement field, while contact pressure is discretized using piece-wise constant shape functions centred at the surface nodes of the quad mesh (collocation method, as described in~\cite{Benaceur2020}). For surface integrals, a single-point gauss quadrature centred on the node is used. This formulation gives mostly smooth contact pressure profiles, except near the peak pressure.\\

\noindent \textbf{Reduced Model}: Snapshots are generated in the training set, $\mathcal{P}_\texttt{tr} \in \mathcal{P}$ consisting of 12 equidistant points in the parametric space. The full rank and truncated dual RB (using CPG) is considered, as given in \cref{tab:hertz_basis_ranks}, whereas the primal RB is not truncated. The full RB has 12 basis vectors for both primal and dual fields, whereas the truncated dual RB has 8 vectors. Note that dual RB truncates only 4 out of 12 basis vectors for a high truncation tolerance of $0.05$. 

\begin{figure}[!htb]
	\centering
	\begin{subfigure}[t]{0.48\linewidth}
		\includegraphics[width=1\linewidth]{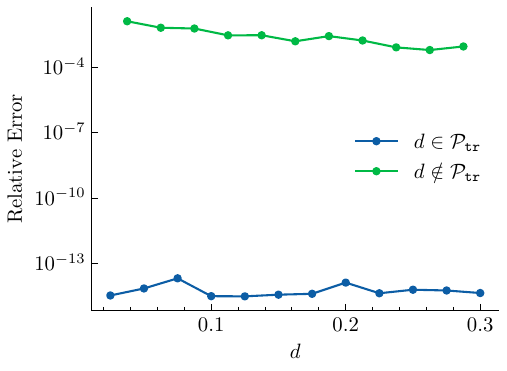}
		\caption{Full RBs}
		\label{fig:hertz_disp_relH1_error_no_trunc}
	\end{subfigure}
	\begin{subfigure}[t]{0.48\linewidth}
		\includegraphics[width=1\linewidth]{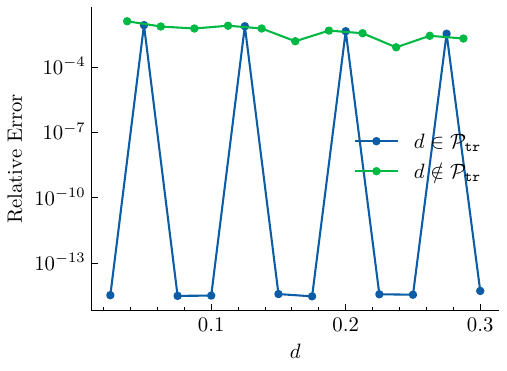}
		\caption{Full primal \& truncated dual RB}
		\label{fig:hertz_disp_relH1_error_trunc}
	\end{subfigure}
    \caption{Primal reconstruction ($\mathcal{H}^1$) errors for the Hertz problem.}\label{fig:hertz_disp_relH1_error}
\end{figure}

\begin{figure}[!htb]
	\centering
    \includegraphics[width=0.8\linewidth]{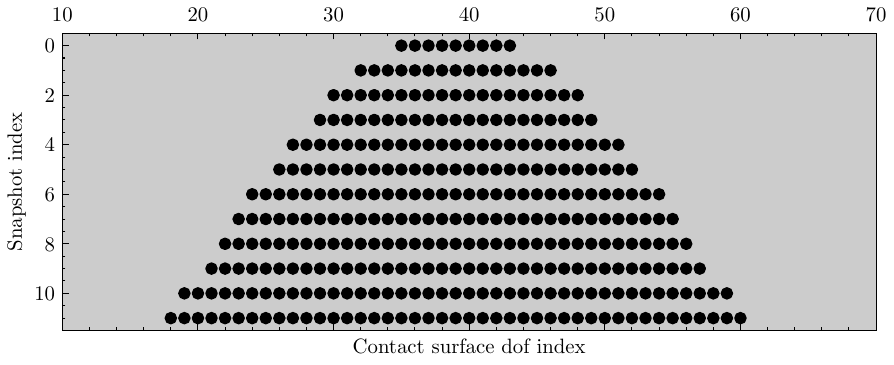}
	\caption{Sparsity pattern of dual snapshots for Hertz problem}\label{fig:hertz_dual_spy}
\end{figure}

\begin{figure}[!htb]
	\centering
	\includegraphics[width=0.5\linewidth]{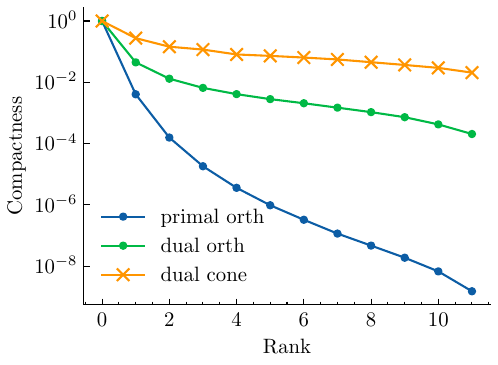}
	\caption{Compactness of reduced bases for Hertz problem}\label{fig:hertz_rb_rank_error}
\end{figure}

The reconstruction errors are shown in \cref{fig:hertz_disp_relH1_error} for full and truncated dual RBs. As expected, the points in the training set show very low reconstruction error, and the points outside the training set show a relatively high error. When the dual RB is truncated, points corresponding to truncated dual snapshots have a high error level, in the same order of error for points outside the training set. Like the illustrative case in \cref{sec:illustration}, the non-zero pattern in snapshots in \cref{fig:hertz_dual_spy} also show a unique sparsity pattern for each snapshot. Also, the comptactness shown in \cref{fig:hertz_rb_rank_error} displays a slow decay of the dual orth and dual cone, indicating a high-rank behaviour of the dual variable. The slightly slower dual cone curve than the dual orth was expected because of the non-negativity constraint that appears in approximation using the cone of the dual basis. These observations again reinforce the proposition that the dual basis does not admit a low-rank behaviour. 

\begin{table}[!htb]
    \small
	\centering
	\begin{tabular}{llcccc}
		\toprule
		&& Primal basis & Dual basis     \\
		&& (POD) & (CP-greedy)   \\
		\midrule
		\multirow{2}{*}{Full} & Tolerance & $\varepsilon$ & $\varepsilon $ \\
		\cdashline{2-4} 
		& Rank & 12 & 12  \\
		\midrule
		\multirow{2}{*}{Truncated} & Tolerance & $\varepsilon$ & $5\times10^{-2}$ \\
		\cdashline{2-4}
		& Rank & 12 & 8  \\
		\bottomrule
	\end{tabular}
    \caption[Truncation tolerances and ranks for both bases for Hertz problem.]{Truncation tolerances and ranks for both bases for Hertz problem. $\varepsilon$ indicates numerical precision}
	\label{tab:hertz_basis_ranks}
\end{table}

\begin{figure}[!htb]
	\centering
	\begin{subfigure}[t]{0.48\linewidth}
		\centering
		\includegraphics[width=1\linewidth]{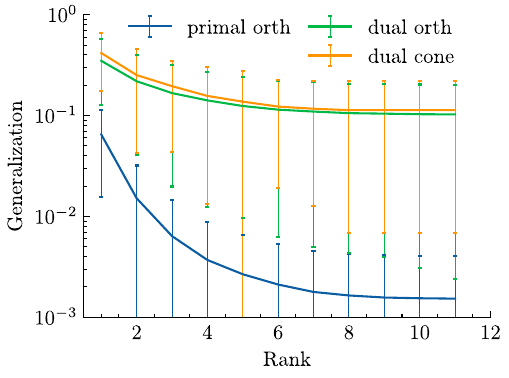}
		\caption{Generalization ability}
	\end{subfigure}
	\begin{subfigure}[t]{0.48\linewidth}
		\centering
		\includegraphics[width=1\linewidth]{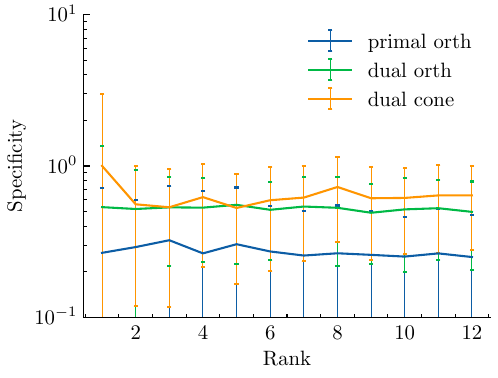}
		\caption{Specificity}
	\end{subfigure}
    \caption[Validation metrics for Hertz problem.]{Validation metrics for Hertz problem. The bars indicate $\pm 1 \sigma$ interval of corresponding metrics}\label{fig:hertz_metrics}
\end{figure}

Generalization ability and Specificity metrics introduced in \cref{sec:metrics} are plotted in \cref{fig:hertz_metrics}. The generalization ability for primal field is much better than the dual field, another reflection of separability issues with the dual field. On the other hand, primal and dual specificity are about the same order, though the primal specificity is consistently lower.

\subsection{Ironing problem}\label{sec:ironing}
\begin{figure}[!htb]
	\centering
    \includegraphics[trim=0 25 0 25, clip, width=0.8\linewidth]{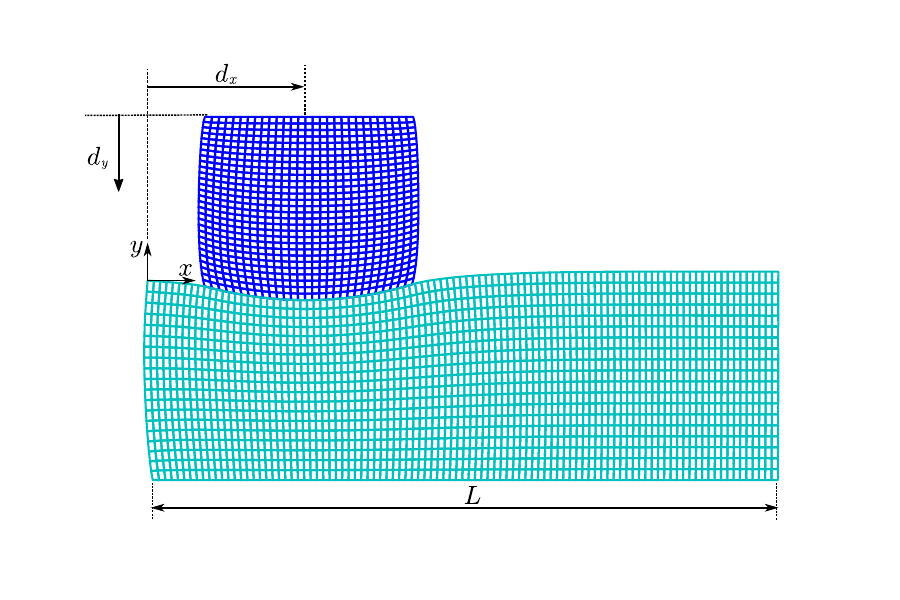}
    \caption[Ironing problem]{Ironing problem: Iron block is pressed against the flat slab by a displacement $d_y$ and moved horizontally. The horizontal displacement of the iron block $d_x \in [0,L]$ is treated as the parameter in the reduced model}\label{fig:ironing}
\end{figure}

\begin{figure}[!htb]
	\centering
	\includegraphics[width=0.6\linewidth]{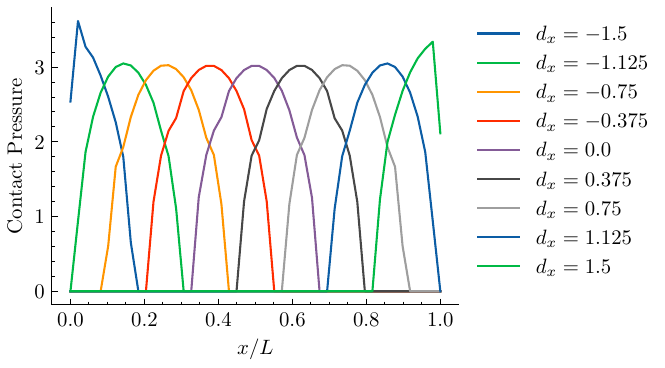}
	\caption{Sample snapshots of contact pressure for ironing problem}\label{fig:ironing_snaps}
\end{figure}

\begin{figure}[!htb]
	\centering
    \includegraphics[width=0.8\linewidth]{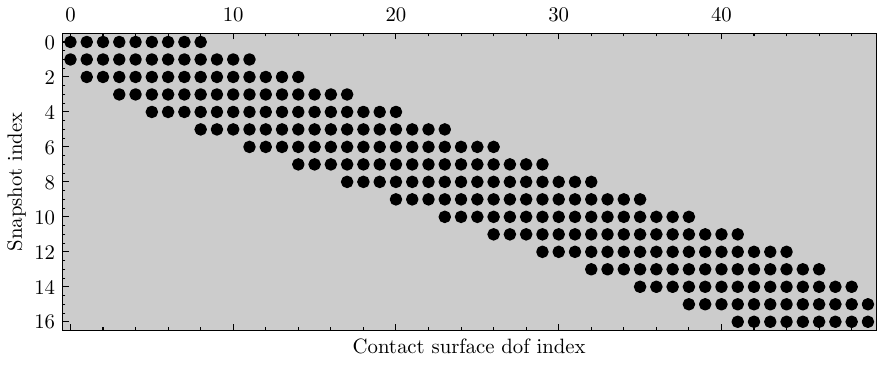}
	\caption{A typical sparsity pattern of dual snapshots for ironing problem}\label{fig:ironing_dual_spy}
\end{figure}

A Type-2 problem where the contact area changes quite significantly, namely the ironing problem~\cite{Fischer2006}, is considered. The ironing problem consists of an iron block pressed against a flat slab and moved along the length of the slab (\cref{fig:ironing}). The problem is simplified with two more assumptions: the first is that iron moves slowly enough that the problem can be considered to be \emph{quasi-static} and the second is that surfaces are frictionless. The horizontal position $d_x$ of the iron is taken as the parameter for the reduced model. In this problem, as the potential contact surface on the slab is quite larger than the actual position of contact, the contact pressure snapshots display large changes in contact position, shown in \cref{fig:ironing_snaps}. The same is reflected in \texttt{spy} pattern of the contact pressure snapshot matrix, shown in \cref{fig:ironing_dual_spy} 

\noindent \textbf{Details of FE model}: A finite element model of the ironing problem is created using a structured quad mesh for both the iron and the slab. Like the Hertz problem, the displacement field is approximated using  linear shape functions and contact pressure using piecewise-constant shape functions centred at the surface nodes. To better demonstrate the inseparability issues, two meshes: coarse and fine, are considered for iron and slab. 
\begin{center}
    \begin{tabular}{lcc}
        \toprule
        & Iron & Slab \\
        \midrule 
        Coarse Mesh & $30\times 30$ & $20\times 100$  \\
        Fine Mesh & $60\times 60$ & $40\times 200$  \\
        \bottomrule
    \end{tabular}
\end{center}

\begin{figure}[!htb]
	\centering
	\begin{subfigure}[t]{0.48\linewidth}
		\centering
		\includegraphics[width=1\linewidth]{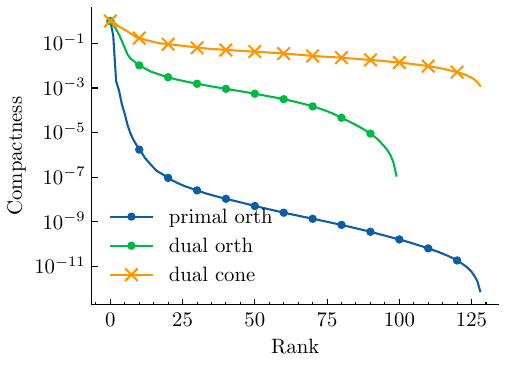}
		\caption{Coarse Mesh}
        \label{fig:irorning_rank_error_coarse_mesh}
	\end{subfigure}
	\begin{subfigure}[t]{0.48\linewidth}
		\centering
		\includegraphics[width=1\linewidth]{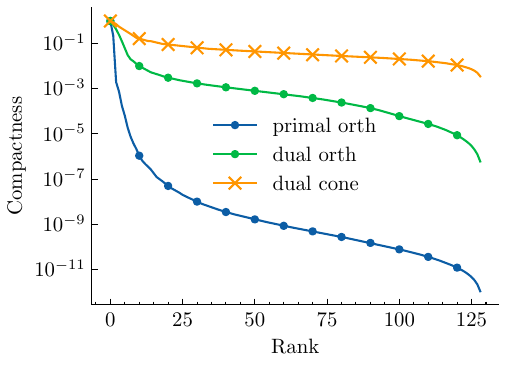}
		\caption{Fine Mesh}
        \label{fig:irorning_rank_error_fine_mesh}
	\end{subfigure}
    \caption{Compactness for reduced bases of ironing problem.}\label{fig:irorning_rb_rank_error}
\end{figure}

The compactness metric, which is equivalent to decay of truncation error as discussed before, is shown in \cref{fig:irorning_rb_rank_error} for a coarse and fine mesh, created using $128$ snapshots. The contrasting part of the two graphs is the green curve representing the dual orth compactness, as this curve vanishes at rank $100$ for the coarse mesh, but not for the fine mesh. This is natural, since coarse mesh has only $100$ nodes on the surface of the slab and thus allowing a maximum rank of $100$, whereas the fine mesh has around $200$ nodes. Since the orthogonally computed dual (dual orth) compactness is also a measure for determining the true rank of the snapshot matrix, the green curve naturally vanishes at a rank of $100$. However, the most interesting part of this figure is the blue curve representing the truncation computed using cone-projection. Even for the coarse mesh, whose true rank is bounded by the number of surface nodes, i.e.\ $100$, the dual cone compactness doesn't decay well beyond this rank. This is another interesting behaviour of the low-rank approach since the dual solution is sought in a subcone that involves positivity constraints, further degrading compactness.

\begin{figure}[!htb]
	\centering
	\begin{subfigure}[t]{0.48\linewidth}
        \centering
        \includegraphics[width=1\linewidth]{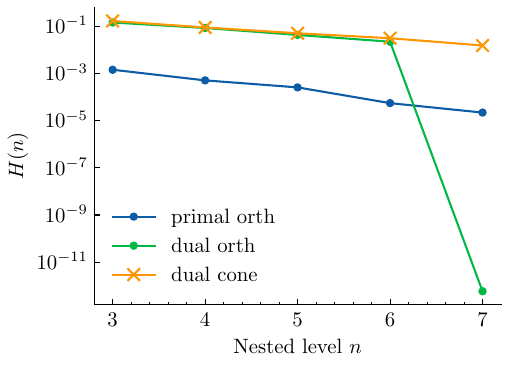}
		\caption{Coarse Mesh}
	\end{subfigure}
	\begin{subfigure}[t]{0.48\linewidth}
        \centering
        \includegraphics[width=1\linewidth]{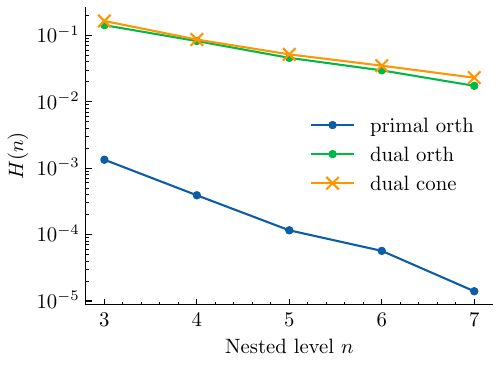}
		\caption{Fine Mesh}
	\end{subfigure}
    \caption{Projection error of nested level $n+1$ snapshots on dual RB of level $n$.}\label{fig:ironing_h_error}
\end{figure}

Apart from compactness, another way to demonstrate the dimensionality of the dual subcone is by using the projection error of nested training sets, which will be referred to as nested error subsequently. To compute the nested error, snapshots are computed in a nested set of points in the parametric space, so that $n$-th level training set is a subset of the $(n+1)$-th level training set. In each nested level, $2^n+1$ points are uniformly distributed in the parametric space $[0,L]$. This means, every alternate point in $(n+1)$-th level is the mid-point of two consecutive points in $n$-th level. The nested error $H(n)$ can then be defined as:
\begin{align}
	H(n) =  \frac{\norm{\mathbf{S}^{n+1}_\texttt{tr} - \Pi_{\Psi^n}(\mathbf{S}^{n+1}_\texttt{tr})}{F}} {\norm{\mathbf{S}^{n+1}_\texttt{tr}}{F}}
	\label{eq:h_error}
\end{align}
\noindent where the symbols $\Psi$ and $\mathbf{S}_\texttt{tr}$ carry same meaning defined in \cref{sec:metrics}. $\norm{\cdot}{F}$ indicates the Frobenius norm. The RB computed using the $n$-th level training set is indicated by  $\Psi^n$, and snapshot matrix for $(n+1)$-th level is indicated by $\mathbf{S}^{n+1}_\texttt{tr}$. Nested training sets where each level contains nearly twice the number of points compared to previous level are created and the dual RB for each level is computed. The slope of nested errors computed using orthogonal and cone projections indicate that the offline stage will require a very large number of snapshots to explore the dual subcone, which is another indication of its high-dimensionality, as shown in \cref{fig:ironing_h_error}. The effects of mesh size are similar to those seen in compactness, as nested error using dual cone projection does not decay even after training set is full rank in case of coarse mesh.

\begin{figure}[!htb]
	\centering
	\begin{subfigure}[t]{0.48\linewidth}
		\centering
		\includegraphics[width=1\linewidth]{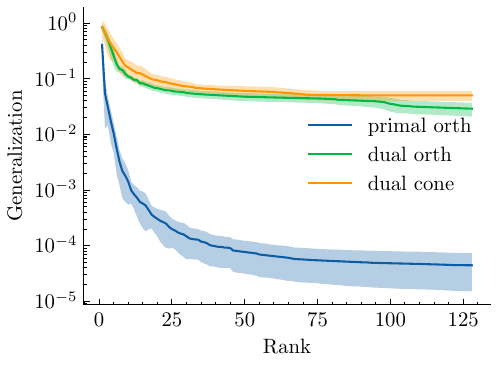}
		\caption{Generalization ability}
	\end{subfigure}
	\begin{subfigure}[t]{0.48\linewidth}
		\centering
		\includegraphics[width=1\linewidth]{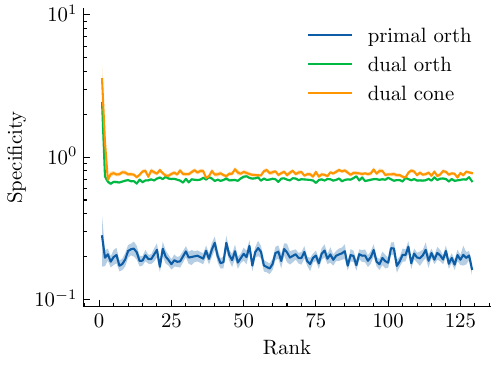}
		\caption{Specificity}
	\end{subfigure}
    \caption[Validation metrics for the ironing problem.]{Validation metrics for the ironing problem (for fine mesh). The fill regions indicate $\pm 1 \sigma$ interval of corresponding metrics.}\label{fig:ironing_metrics}
\end{figure}

Generalization ability and Specificity metrics are plotted in \cref{fig:ironing_metrics}. The generalization ability shows a similar trend as that of nested error, as they are similar quantities measuring reconstruction errors outside the training set. The primal basis outperforms its dual counterparts in this metric, as expected. The specificity metric is quite higher for dual quantities than the primal ones, which is because random shapes generated using the RB are quite different.

\section{Conclusions and Perspectives}
In this chapter, the limitations of low-rank ROMs to contact mechanics problems were discussed. Demonstrations of this approach show moderate reconstruction errors for parametric points outside the training set. To explore the possible sources of this error, we focus on the key assumption of the low-rank hypothesis i.e.\ the solution subspace is low-dimensional. It turns out that the Lagrange Multiplier associated with the inequality constraints, i.e.\ the contact pressure, does not admit a low-rank subspace. Qualitative and quantitative assessments are provided to support this argument. The linear inseparability of dual subspace/subcone is demonstrated using various metrics that show high-dimensionality of contact pressure compared to the relatively low-dimensionality of the displacement field.\\

To circumvent the limitations due to the lack of low-rank structure, several potential approaches need to be explored. One possibility is the usage of dictionary-based approximation with a relatively large number of snapshots spanning the dual subspace efficiently. Another possibility is resorting to either a non-linear transformation; potentially transforming the dual snapshots in a way that they lie in a low-dimensional subspace. 

\chapter{Dictionary approach in contact mechanics problems}\label{ch:sparse}

\section{Motivation}\label{sec:philosophy_dictionary_approx}
The contact mechanics problem has a feature, the contact zone, that exhibits a local nature. The contact zone is a function of the state of the system, and hence, may strongly depend on the loading, geometry, and other physical parameters. In general, for variational inequalities, only a subset of the constraints are usually active. An inequality constraint is referred to as active only if the current state of the problem satisfies the equality. Due to the last KKT condition in~\eqref{eq:kkt_reduced_comp_slack}, the Lagrange multipliers $\lambda$ assume a non-zero value only for active constraints i.e.\ the part of surface where contact is established. The strong dependence of the active zone on the state of the system leads to the linear inseparability of the Lagrange multiplier field. It has been shown in \cref{ch:lowrank}, that the local characteristics of contact pressure can render subspace learning ineffective and cause difficulties in creating efficient reduced models of contact problems. The separability issues  may cause the following problems: 

\begin{enumerate}
    \item \textbf{Complexity of subspace learning}: It is desirable to keep the complexity of the offline stage at affordable levels. Usually, this is achievable because a relatively low number of snapshots compared to the dofs of the high-fidelity problem are sufficient to create a reliable reduced basis, thanks to the low-rank nature of many resolution problems. However, given the inseparability of the contact pressure field, a large training set may be necessary to discover the entire solution subspace to a reasonable accuracy.\label{item:high_offline_cost}
    \item \textbf{Lack of low-rank subspace}: Discovery of a low-rank subspace is key to creating a low-cost reduced model, as the number of reduced dofs scales with the rank of the discovered solution subspace. While low-rank methods for contact problems~\cite{Balajewicz2016,Benaceur2020} define a contact pressure subcone, it is shown in \cref{ch:lowrank} using various metrics that compression of contact pressure data into a low-rank subcone or subspace leads to loss of useful information. Therefore, despite the creation of a large training set, as mentioned in the previous point, extracting a low-rank subspace might not be possible.\label{item:not_compressible}
\end{enumerate}

The degree of these effects depends on the nature of the contact problem. In this chapter, the approximation using a dictionary instead of a low-rank subcone is explored, addressing the second point as the dictionary is made of snapshots without applying any compression, like NNMF or CPG\@.\ In effect, the high complexity of the offline stage is accepted as a compromise for creating a reliable reduced model. The ideas behind over-complete dictionaries are expanded in the coming sections.

\section{Sparse methods with over-complete dictionaries}
The primary motivation for using a dictionary-based approximation is to mitigate the limitations due to linear inseparability, discussed in \cref{sec:philosophy_dictionary_approx}. The word ``mitigate'' is explicitly used here to indicate that the strategy developed in this chapter does not resolve the inseparability issue at a fundamental level, rather it attempts to limit the effects of inseparability. The intention is to accommodate the cost of computing snapshots over a large training set and focus mainly on addressing the point~\ref{item:not_compressible} of \cref{sec:philosophy_dictionary_approx} using an algorithm that can handle large dictionaries efficiently. A matrix containing the large set of snapshots is referred to as an over-complete dictionary $\mathbf{D}$. The term ``over-complete'' is used loosely in this context, as it is impossible to sufficiently sample a truly inseparable subspace for a given approximation error target and is simply used to imply that the dictionary has numerous snapshots, possibly larger than the number of contact dofs in the high-fidelity problem.

\begin{figure}[!htb]
    \centering
    \includegraphics{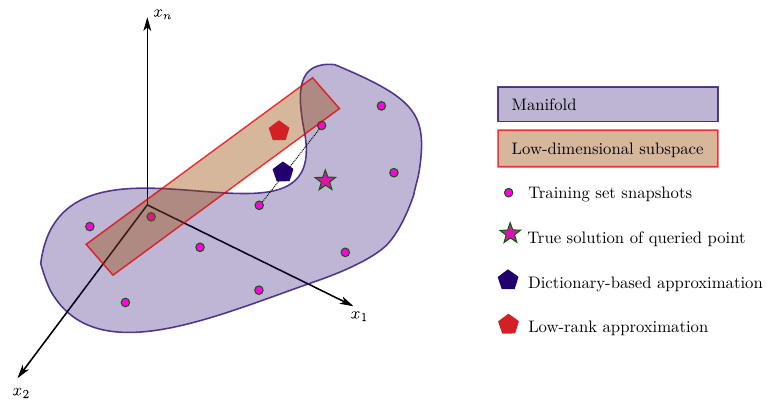}
    \caption[Illustration of low-rank vs.\ dictionary-based approximation]{Illustration of low-rank vs.\ dictionary-based approximation. Snapshots (circles) in an $n$-dimensional system lie on a manifold of dimension lower than $n$. The low-dimensional subspace computed using low-rank methods has a dimension $\ll n$. A dictionary-based approximation of the queried point (star-shaped) can be computed as linear combination of the nearest snapshots (shown by a line connecting the nearest candidates). The low-rank approximation is simply the projection of the queried point on the low-dimensional subspace.}
    \label{fig:dictionary_vs_low_rank}
\end{figure}

The dictionary-based approximation can be visualized in an abstract $n$-dimensional space in \cref{fig:dictionary_vs_low_rank}. Typically, in cases where low-rank hypothesis is not satisfied, the underlying physics of the problem assume a manifold structure. A manifold can be thought of as a low-dimensional entity that does not satisfy the properties of a linear space. The idea expressed in this figure is that a manifold, like a hypersurface with non-zero curvature, can be approximated better using a large cloud of points lying on it. If the point cloud is rich enough, it is more likely that this approximation will be better than than low-dimensional subspace which is a hyperplane\footnote{The terms hyperplane and hypersurface refer to a plane and surface in dimensions higher than 3 respectively.}. It is this large cloud of points that is referred to as over-complete dictionary in this chapter.

In dictionary-based approximation, unlike the low-rank approach~\cite{Balajewicz2016,Benaceur2020}, no compression of the contact pressure snapshots is applied. The dual solution can simply be restricted to the column-cone\footnote{column-cone is the subset of the column-space restricted to non-negative coefficients.} of the dictionary $\mathbf{D}$, which might be reasonably rich to approximate the true contact pressure subcone. In the online phase, only a few snapshots will be selected to estimate the contact pressure field. This is justified because not all information in the dictionary is necessary for the reconstruction of a particular instance and hence, the dual reduced dofs will admit a sparse solution. This feature is also seen in case of the low-rank model generated using snapshots selected by CPG\@. Therefore, enough motivation exists to use sparsity enforcing methods to choose a few columns of a dictionary $\mathbf{D}$ to approximate the dual solution. Keeping these ideas in mind, sparse regression techniques are briefly introduced at this stage and a few of these methods from current literature are discussed. Following this, their application to resolution problems will be briefly discussed.

\subsection{Sparse regression in data approximation problems}\label{sec:sparse_methods}
Sparse regression techniques have been under development  in recent decades across various fields, especially in signal processing and statistics. The main idea behind these techniques is to approximate a signal with the few most suitable elements of a precomputed dictionary of signals~\cite{Elad2010}, or to develop the simplest predictive model~\cite{Brunton2016}. A typical sparse regression can be written as the following optimization statement:
\begin{align}
    \centering
    \begin{aligned}
        \min & \norm{\bm{\alpha}}{p} \\
        \text{s.t.} & \norm{\bm{x} - \mathbf{D} \bm{\alpha}}{}  < \varepsilon
    \end{aligned}
    \label{eq:parsimony_idea}
\end{align}
where $x$ is the signal being approximated, $\mathbf{D}$ could either be a dictionary of signals or a basis that possibly admits sparse representation, $\norm{\cdot}{p}$ is a sparsity inducing norm, and $\varepsilon$ is the tolerance on the approximation error. It is also possible to write other forms of sparse regression problems where the approximation error is minimized and $p$-norm term is constrained to stay under a specific tolerance of sparsity, or the minimization of a penalty form with a weighted sum of the two terms.

It is common knowledge in the communities using sparse regression methods that the $\ell_0$-norm is a measure of cardinality, hence~\eqref{eq:parsimony_idea} with $p=0$  is the natural choice in a theoretical sense. However, this problem is NP-hard in general, and therefore, it is more common to use $\ell_p$-norm with $p\in ]0,1]$. A detailed explanation of sparsity inducing properties of $\ell_p$-norms with $p \in ]0,1]$ is given in~\cite{Mairal2014}. A few sparse regression methods are briefly detailed here:

\begin{enumerate}
    \item \textbf{LASSO}~\cite{Tibshirani1996}: The specific case of~\eqref{eq:parsimony_idea} with $p=1$ is known as the LASSO is one of the most widely known sparse regression methods. The $\ell_1$-norm, even though not a true measure of cardinality, has a constant ``driving force'' or more technically speaking, a gradient towards the origin that is constant in magnitude everywhere. This drives many entries of $\bm{\alpha}$ to zero, unlike the $\ell_2$ norm that can reduce the magnitude of $\bm{\alpha}$ entries, but not set them to zero. Also, the $\ell_1$-norm allows solving~\eqref{eq:parsimony_idea} in polynomial time.
    \item \textbf{Dantzig Selector}~\cite{Candes2007}: This method allows selection of a sparse $\bm{\alpha}$ using a modified version of~\eqref{eq:parsimony_idea}, that can be converted into a linear programming problem after with algebraic manipulations. This is achieved with $p=1$ and the $\ell_\infty$ norm for the residual.~\eqref{eq:parsimony_idea} is recast as:
\begin{align*}
    \begin{aligned}
        \min & \norm{\bm{\alpha}}{1} \\
        \text{s.t.} & \norm{\mathbf{D}^T(\bm{x} - \mathbf{D} \bm{\alpha})}{\infty}  < \varepsilon
    \end{aligned}
\end{align*}
which can be expressed in a linear programming form as:
\begin{align*}
    \begin{aligned}
        \min & \sum_i \gamma_i \\
        \text{s.t.} & - \bm{\gamma} < \bm{\alpha}  < \bm{\gamma} \\
        \phantom{\text{s.t.}} & -\varepsilon \bm{1} < \mathbf{D}^T(\bm{x} - \mathbf{D} \bm{\alpha}) < \varepsilon \bm{1}
    \end{aligned}
\end{align*}

\noindent where $\bm{\gamma}$ is a dummy variable obtained by opening the absolute operator $\| \bm{\alpha} \|$.

    \item \textbf{Orthogonal Matching Pursuit}~\cite{Pati1993}: The OMP method aims to approximate the solution of~\eqref{eq:parsimony_idea} with $p=0$ using a greedy technique. At each greedy step $k$, the OMP algorithm searches for the column of dictionary $\mathbf{D}$ that is the least orthogonal to the current residual $(\bm{x}-\mathbf{D} \bm{\alpha}^{k-1})$ and adds to the set of indices previously selected. This process is carried out until the residual is below the tolerance $\varepsilon$. This is done by projecting the residual vector on each column and selecting the largest projection.
\begin{align*}
    \mathcal{I} \gets \mathcal{I} \cup \argmax_j \left |\mathbf{D}[:,j]^T (\bm{x}-\mathbf{D} \bm{\alpha}^{k-1})\right |
\end{align*}
where $\mathcal{I}$ is the current set of selected indices and columns of dictionary $\mathbf{D}$ are normalized. Then the new set of coefficients $\bm{\alpha}^k$ are computed using the least square solution:
\begin{align*}
        \bm{\alpha}^k_\mathcal{I} = (\mathbf{D}_{\mathcal{I}}^T \mathbf{D}_{\mathcal{I}})^{-1} \mathbf{D}_{\mathcal{I}}^T\bm{x}
        \label{eq:omp_updaterep}
\end{align*}
with $\bm{\alpha}_\mathcal{I} =\bm{\alpha}[\mathcal{I}]$ and $\mathbf{D}_\mathcal{I} = \mathbf{D}[:,\mathcal{I}]$.

Update of $\mathcal{I}$ is locally optimal, but the above least-square update of $\bm{\alpha}^k$, the OMP solution is optimal w.r.t to the currently selected subset $\mathbf{D}_\mathcal{I}$ of the dictionary. This is in contrast to the predecessor of OMP, the Matching Pursuit~\cite{Mallat1992}, which used locally optimal updates for $\bm{\alpha}^k$ at the current iteration.

\item \textbf{FOCUSS}~\cite{Gorodnitsky1997}: The FOCUSS method finds a sparse approximation using an iterative process that starts from a fully dense solution and progresses towards ``localized energy solutions''. The FOCUSS method is broadly based on the following ideas:
    \begin{itemize}
        \item For a full rank dictionary of size, say $m \times n\ (m < n)$, the system $\mathbf{D} \bm{\alpha} = \bm{x}$ is under-determined. For such a system, the closed form solution for~\eqref{eq:parsimony_idea} with $p=2$ is given by $\bm{\alpha} = \mathbf{D}^\dagger \bm{x}$, where $\mathbf{D}^\dagger$ is the Moore-Penrose pseudo-inverse. Note that this solution is not sparse in general, as $\ell_2$-norm does not have sparsity inducing properties. 
        \item A modified minimization problem can be defined by replacing the $\norm{\bm{\alpha}}{2}$ with a weighted norm $\norm{\mathbf{W}^{-1} \bm{\alpha}}{2}$, or more generically $\norm{\mathbf{W}^\dagger \bm{\alpha}}{2}$ if $\mathbf{W}$ is singular. The authors of~\cite{Gorodnitsky1997} argue that by changing $\mathbf{W}$, every possible solution of the under-determined system $\mathbf{D} \bm{\alpha} = \bm{x}$ can be obtained.
        \item The basis of the algorithm lies in iterating towards a sparse solution using a weight $\mathbf{W}^k$ which induces sparsity. For the iteration $k$, the weight is chosen as $\mathbf{W}_{k} = \mathtt{diag}(\bm{\alpha}^{k-1})$. The trick lies in the fact that algorithm ends up minimizing the following weighted norm:
            \begin{equation*}
                \norm{\mathbf{W}_k^\dagger \bm{\alpha}}{2}^2 = \sum_i \left (\frac{\alpha^k_i}{\alpha^{k-1}_i} \right )^2 
            \end{equation*}
        due to which smaller entries of $\alpha^{k-1}$ tend to diminish further.
    \end{itemize}
    On combining these ideas, the FOCUSS algorithm ends up with quite a simple implementation, given in \cref{alg:focuss}. Usually, the algorithm is initialized with the $\bm{\alpha}^0$ containing all non-zeros. In~\cite{Gorodnitsky1997}, the $\norm{\bm{\alpha}}{2}$ minimizing solution $\bm{\alpha}^0 = \mathbf{D}^\dagger \bm{x}$ was used for initialization. 
    \begin{algorithm}[h]
        \caption{FOCUSS} \label{alg:focuss}
        \begin{algorithmic}[1]
            \State{Inputs: $\mathbf{D}, \bm{x}$, $\bm{\alpha}^0$}
            \State{Initialize $k = 1$}
            \While{$\bm{\alpha}$ not converge}
            \State{$\mathbf{W}^k = \mathtt{diag}(\bm{\alpha}^{k-1})$}
            \State{$\bm{\alpha}^k = \mathbf{W}^k (\mathbf{D} \mathbf{W}^k)^\dagger \bm{x}$} \label{step:focuss_alpha}
            \State{$k \gets k+1$}
            \EndWhile{}
        \end{algorithmic}
    \end{algorithm}
\end{enumerate}

\subsection{Application of sparse methods to unconstrained resolution problems} \label{sec:dictionary_resolution}
Sparse regression methods, though originally developed for application on pure data approximation or data compression problems, have seen recent interest in applications towards physics problems~\cite{Wang2011,Schaeffer2013,Brunton2016}. The methods that have been discussed previously are also meant for data approximation. At this point, application of such methods on unconstrained resolution problems is discussed. A resolution problem can be framed as residual minimization problem in order to use dictionary-based approximation and sparse regression. Consider a discretized resolution problem with the state vector $\bm{v}$ and residual $\bm{r}(\bm{v})$. The basic idea is to approximate the solution $\bm{v} \approx \mathbf{D} \widehat{\bm{v}}$ using the elements of dictionary such that the residual of the numerical model $\bm{r}(\mathbf{D} \widehat{\bm{v}})$ is below a certain tolerance, while minimizing the $\ell_p$-norm of coefficient vector, as shown in~\eqref{eq:sparse_res}. Sparsity is imposed on the coefficients of the dictionary, meaning that only a few dictionary columns are selected.

\begin{gather}
    \begin{aligned}
        \min  & \norm{\widehat{\bm{v}}}{p} \\
        \text{s.t. } &\norm{\bm{r}(\mathbf{D} \widehat{\bm{v}})}{} < \varepsilon
    \end{aligned}
    \label{eq:sparse_res}
\end{gather}

Dictionary-based sparse approximation for resolution problems has been discussed in~\cite{Balabanov2021}, where the solution space is set to the span of dictionary columns. They also discuss the application of random sketching to residual evaluation, which provides an efficient way of solving the minimization problem. This is similar to the approach of randomized SVD proposed in~\cite{Halko2011}. The idea is to project the residual on a low-dimensional subspace in $\text{colsp}(\mathbf{D})$.\footnote{colsp indicates column-space} The matrix $\mathbf{B}$, containing a basis corresponding to a random combination of dictionary columns, is used for this projection in~\eqref{eq:sparse_res_randomized}. This essentially reduces the complexity of enforcing the residual constraint. Given a dictionary $\mathbf{D}$ of size $m\times n$, a smaller matrix that has information from all columns of $\mathbf{D}$ can be generated by multiplying a matrix $\mathbf{R}$ of size $n\times L$ with $L \ll \min{(m, n)}$. After orthogonalization, this matrix $\mathbf{B}$ can be used for projecting the residual $\bm{r}$.
\begin{subequations}
    \begin{align}
        \min & \norm{\widehat{\bm{v}}}{p} \\
        \text{s.t. }&\norm{\mathbf{B}^T \bm{r}(\mathbf{D} \widehat{\bm{v}})}{} < \varepsilon \\
        \text{where } &\mathbf{B} = \text{orth}(\mathbf{D} \mathbf{R}) \ , \ \ \mathbf{R} \in \mathcal{U}_{[0,1]}^{n \times L} \label{eq:randomized_orth_matrix}
    \end{align}
    \label{eq:sparse_res_randomized}
\end{subequations}
The number $L$ can be tuned to not only control the complexity of the minimization problem, but also to crudely constrain sparsity of the solution. For $L \ll n$, $\text{colsp}(\mathbf{B}) \subset \text{colsp}(\mathbf{D}$) is low-dimensional and hence, the projection $\mathbf{B}^T \bm{r}$ is a ``weaker'' evaluation of $\bm{r}$, which allows for a sparser $\widehat{\bm{v}}$ at the price of permitting residuals orthogonal to $\mathbf{B}$.

Another related approach in which snapshots are directly used for reducing the model without any explicit compression is the CUR decomposition~\cite{LeGuennec2018}. In this method, not only columns, but also rows from the snapshot matrix are used for developing a parametric regression. It works by algorithmically selecting a column and row subset from the snapshot matrix, and using them for developing a regression model. A contrasting feature of the CUR approach is the use of regression without invoking the partial differential equation or its weak form, resulting in a non-intrusive method.

\subsubsection{Diffusion-reaction problem}
\begin{figure}[htpb]
    \centering
    \includegraphics[width=0.5\linewidth]{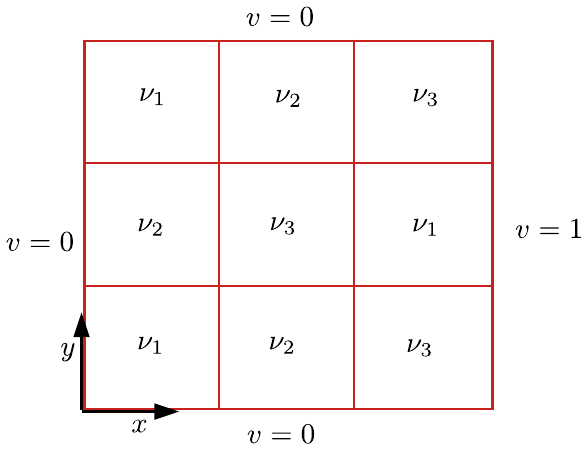}
    \caption{Two-dimensional domain for the transport problem~\eqref{eq:transport}}
    \label{fig:transport_domain}
\end{figure}

A transport problem of a quantity $v$ with two-dimensional domain of size $3 \times 3$ units, divided into a grid with cells of size $1\times 1$ is considered. The cells have a piecewise uniform diffusivity value. A three-dimensional parametrization is considered with diffusivity values $\nu_{1}, \nu_{2}$ and $\nu_{3}$, as shown in \cref{fig:transport_domain}. The reaction coefficient $\sigma$ is considered to be uniform across all sectors.
\begin{gather}
    \begin{aligned}
        - \nu(x,y) \nabla^2 v(x,y) + \sigma v(x,y) = f \\
        v(x,0) = 0 \\
        v(x,3) = 0 \\
        v(0,y) = 0 \\
        v(3,y) = 1 \\
    \end{aligned}
    \label{eq:transport}
\end{gather}

\noindent where $\nu(x,y)$ describes the piecewise uniform diffusivity values.

\begin{figure}[!htpb]
    \centering
    \begin{subfigure}[t]{0.8\linewidth}
        \includegraphics[width=1\linewidth]{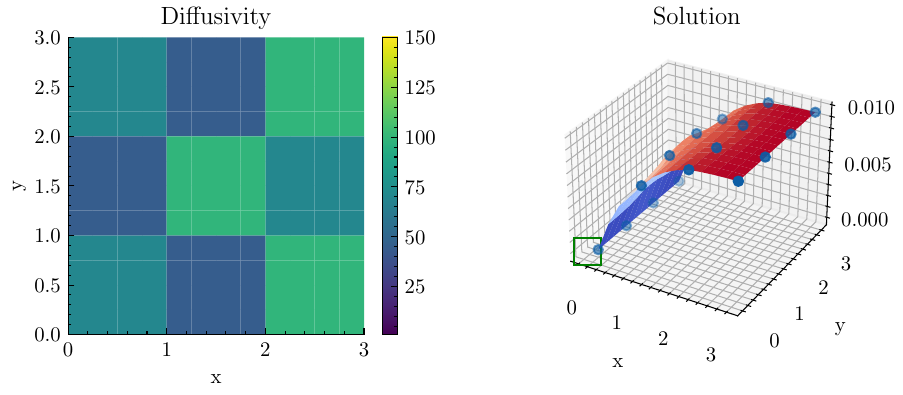}
        \caption{}
    \end{subfigure}
	\begin{subfigure}[t]{0.8\linewidth}
        \includegraphics[width=1\linewidth]{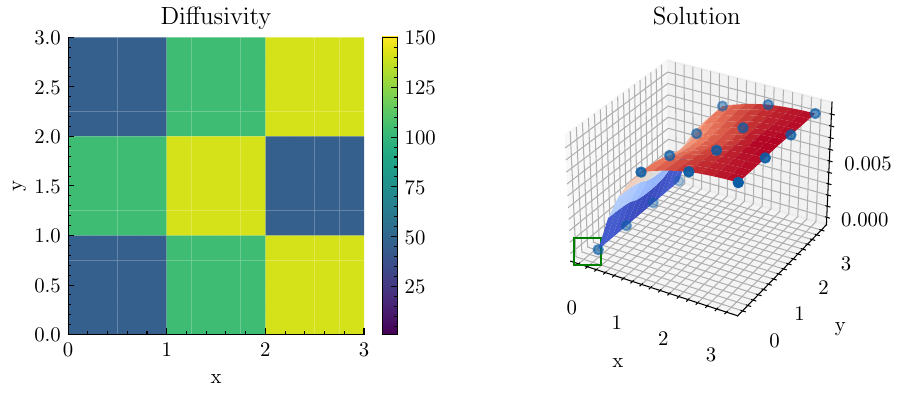}
        \caption{}
    \end{subfigure}
	\begin{subfigure}[t]{0.8\linewidth}
        \includegraphics[width=1\linewidth]{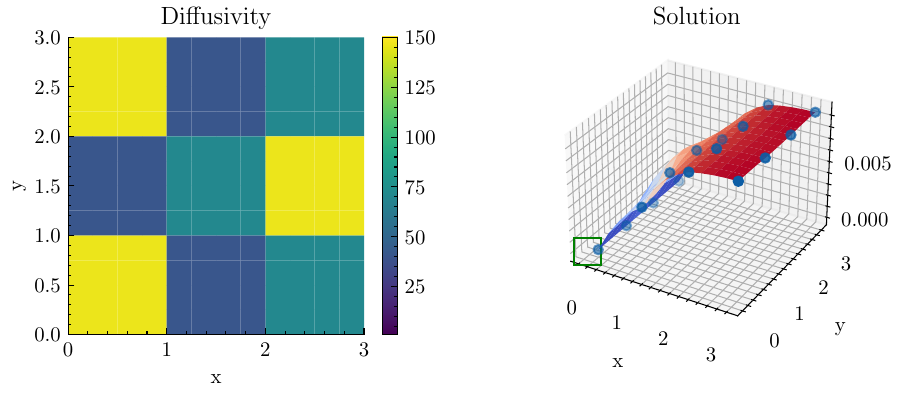}
        \caption{}
    \end{subfigure}
    \caption{Snapshots of the transport problem~\eqref{eq:transport} for randomly sampled parametric values}
    \label{fig:transport_snaps}
\end{figure}

The 2D domain is discretized into a structured $15\times15$ mesh. The diffusivity of each block $\nu_{i} \in (1,150)$ is considered for the parametric model. The reaction coefficient is chosen as $\sigma = 10$. The reaction is included mainly to add some local effects to an otherwise globally diffusive problem. A dictionary of 320 snapshots is created using randomly sampled diffusivity values. \cref{fig:transport_snaps} shows some sample snapshots from the training set.

\begin{figure}[!htb]
    \centering
    \begin{subfigure}[t]{0.95\linewidth}
        \includegraphics[width=1\linewidth]{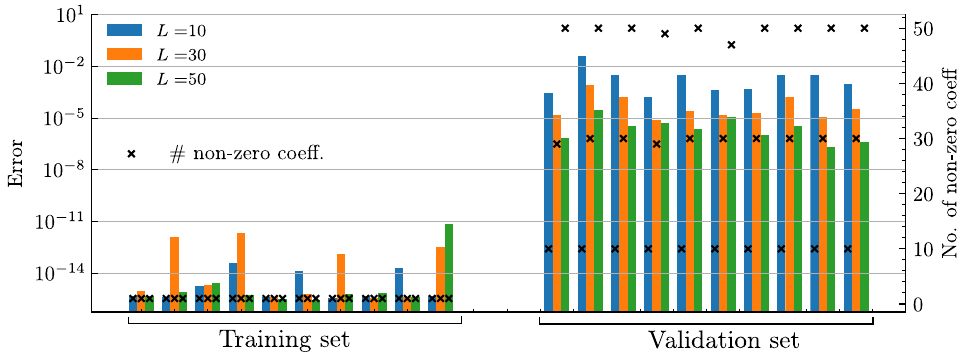}
        \caption{FOCUSS}
    \end{subfigure}
	\begin{subfigure}[t]{0.95\linewidth}
        \includegraphics[width=1\linewidth]{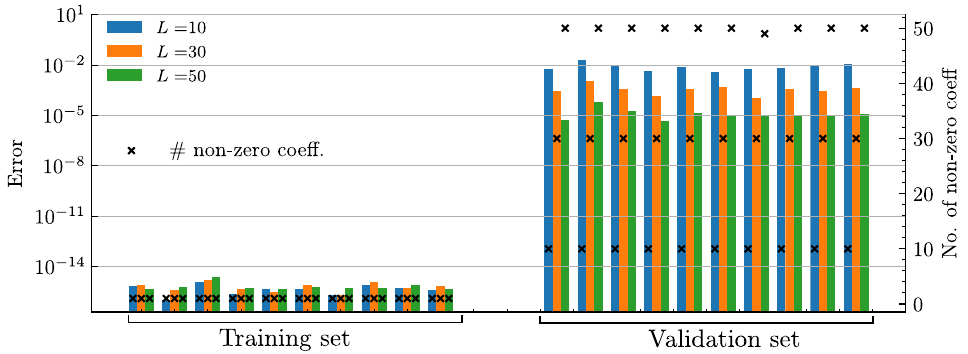}
        \caption{OMP}
    \end{subfigure}
    \caption[Reconstruction errors of the transport problem using dictionary-based approximations.]{Reconstruction errors (bars) of the transport problem using dictionary-based approximations and number of non-zero entries ({\tiny ${\boldsymbol \times}$}) on training and validation sets. $L$ indicates the target sparsity}
    \label{fig:transport_recon}
\end{figure}

\cref{fig:transport_recon} shows reconstruction errors using the dictionary-based approximation methods using~\eqref{eq:sparse_res_randomized}. FOCUSS and OMP are used for sparse regression. As expected, the training set snapshots are almost accurately reconstructed with a single non-zero coefficient. In the validation set, reconstruction errors are observed in the order of $10^{-2}$ to $10^{-1}$ for sparsity target $L=10$ and errors improve about an two orders for $L=30$. The outcome of sparsity for each $L$, shown using a cross ($\boldsymbol{\times}$), is almost the same as the parameter $L$ for the validation set. This is due to the algorithm trying to select as many vectors permitted by the user. 

Now that the dictionary methods have been demonstrated on unconstrained methods, we would like to extend these methods to inequality-constrained problems. However, it is not straightforward to apply random sketching to problems involving inequality constraints. The projection of inequality constraints on a randomized matrix like $\mathbf{B}$ creates linear combinations of constraints, which makes it harder to identify which constraints are being violated and which are not. In upcoming sections, two dictionary-based approximation schemes will be studied, but the projection of residual on $\mathbf{B}$ is applied in only one of them.

\section{Dictionary methods for contact problems}\label{sec:sparse_dict}
Now that the feasibility of dictionary-based approximations have been established in an unconstrained problem, we move towards inequality constrained problems. To this end, two dictionary-based approximation methods for contact mechanics problems have been formulated, namely the \emph{greedy active-set} approach and \emph{convex-hull approximation} approach. To provide a broad overview of the features of these methods, \cref{tab:dictionary_methods} gives a list of features of these methods. The important differences between the two methods are bases types and the manner of enforcing constraint. These details will be discussed along with their formulations and application to contact problems are discussed in this section. 

\begin{table}[!htb]
	\setlength\extrarowheight{2pt}
    \small
	\centering
    \begin{tabular}{m{0.3\textwidth}m{0.24\textwidth}m{0.05\textwidth}m{0.26\textwidth}}
        \toprule
                                & Greedy active-set     & & Convex hull approximation \\
        \midrule                                          
        Primal Basis type       & Truncated low-rank    & & Dictionary \\
        Dual Basis type         & Dictionary            & & Dictionary \\
        Monolithic dictionary   & No                    & & Yes\\              
        Iteration method        & Fixed point           & & Fixed point \\ 
        Constraint enforcement  & Active-set            & & Convex combinations of non-penetrating snapshots \\
        Sparsity induction      & Greedy activation and disactivation of dual dofs
                                                        & & Non-negative FOCUSS\\ 
        Additional assumptions  & None                  & & Convexity of the feasible region\\
        \bottomrule
    \end{tabular}
    \caption{Characteristics of the two dictionary-based approximation methods for contact mechanics problems}
	\label{tab:dictionary_methods}
\end{table}

\subsection{A greedy active-set method for dual dictionary element selection}
To solve the contact problem using a dictionary of contact pressure snapshots, a greedy method inspired from the OMP and active set method is devised. From this point onwards, $\DualDict$ indicates a dictionary of dual snapshots of the contact problem. Consider the following sparse problem:

\begin{subequations}
	\begin{align}
		&\min \norm{\widehat{\bm{\lambda}} }{0} \\ 
		& \text{s.t. }  \begin{bmatrix}
			\mathbf{\PrimalRB}^T \mathbf{K} \mathbf{\PrimalRB}  & \mathbf{\PrimalRB}^T \mathbf{C}^T \DualDictActive \\
			\DualDictActive^T \mathbf{C} \mathbf{\PrimalRB} & \mathbf{0}
		\end{bmatrix} 
		\begin{bmatrix}
			\widehat{\bm{u}} \\
			\widehat{\bm{\lambda}}_{\mathcal{I}}
		\end{bmatrix} =  
        \begin{bmatrix}
			\mathbf{\PrimalRB}^T \bm{f} \\
			\DualDictActive^T \bm{g}
		\end{bmatrix} \label{eq:sparse_greedy_system} \\
		&  \text{and } \DualDict^T  \left ( \mathbf{C} \mathbf{\PrimalRB} \widehat{\bm{u}} - \bm{g} \right ) \leq \tau  \label{eq:sparse_greedy_kkt_nonpen}\\
        &   \phantom{\text{and }} \widehat{\bm{\lambda}} \geq 0 \label{eq:sparse_greedy_kkt_nnpr}\\
		& \text{where $\mathcal{I}$ indicates the active set and }  \DualDictActive = \DualDict[:,\mathcal{I}] \nonumber
	\end{align}
\label{eq:sparse_greedy}
\end{subequations}

\noindent\eqref{eq:sparse_greedy_system},~\eqref{eq:sparse_greedy_kkt_nonpen}~and~\eqref{eq:sparse_greedy_kkt_nnpr} are based on the same principle as the reduced KKT conditions~\eqref{eq:kkt_reduced}, but applied over a dictionary $\DualDict$ rather than a low-rank dual basis. Comparing~\eqref{eq:sparse_greedy} with the generic form of sparse regression statement~\eqref{eq:sparse_res}, the balance equation is imposed as an equality constraint, instead of an inequality constraint where the residual norm is limited to certain tolerance. Imposing limits on the residual norm is usually more relevant when using random sketching, like in~\eqref{eq:sparse_res_randomized} when the system is under-determined. Also, it is not straightforward to use standard sparse regression algorithms in presence of the KKT inequality constraints~\eqref{eq:sparse_greedy_kkt_nonpen}~and~\eqref{eq:sparse_greedy_kkt_nnpr}. Equality constraint can be used as long as the system is determined, which is currently the case in~\eqref{eq:sparse_greedy_system}. In~\eqref{eq:sparse_greedy_kkt_nonpen}, $\tau$ indicates a small penetration that has been permitted in the algorithm, which is necessary to ignore spurious penetrations created by the low-rank displacement basis $\PrimalRB$. Unless specified, the $\tau$ is set to a value in the same order of truncation tolerance $\delta$ (see \cref{ch:tau_delta}).

The primary notion behind greedy algorithms is to enrich the approximation by a rank-1 term in each step. For standard sparse regression problems like~\eqref{eq:sparse_res}, the OMP algorithm starts with an empty coefficient vector and adds a single non-zero entry in one step. The position of the new non-zero entry is chosen by searching the vector in the dictionary with the highest correlation to the residual. The greedy active-set algorithm devised here performs the greedy enrichment in an almost similar way but with a modification, i.e.\ the orthogonality with the violation of inequality constraint~\eqref{eq:sparse_greedy_kkt_nonpen} is used instead of the residual of system~\eqref{eq:sparse_greedy_system}. More precisely, the algorithm starts with a zero $\widehat{\bm{\lambda}}$ and in each enrichment step, the dictionary column with highest correlation to the violations of the non-penetration condition is added in the greedy enrichment. Thus, enrichment index of the dictionary $p_{\texttt{add}}$ can be computed as:

\begin{align}
    p_{\texttt{add}} &:= \argmax_p \DualDict[:,p]^T  \left ( \mathbf{C} \PrimalRB \widehat{\bm{u}}^{k} - \bm{g} \right )^{+, \tau}
    \label{eq:sparse_enrich_concept}
\end{align}
 where  $
	(\bm{z})^{+,\tau}_i :=\left \lbrace\begin{matrix}
		z_i & \text{if } z_i>\tau \\
		0 & \text{otherwise}
    \end{matrix}  \right . $, which is a hard-thresholding operator. \\

\noindent\eqref{eq:sparse_enrich_concept} can be computed more efficiently using the reduced operators $\widehat{\mathbf{C}}(\mu,  \widehat{\bm{u}}^{k-1})$ and $\widehat{\bm{g}}(\mu,  \widehat{\bm{u}}^{k-1})$:
\begin{align}
    p_{\texttt{add}} &:= \argmax_p  \left ( \widehat{\mathbf{C}} \widehat{\bm{u}}^{k} - \widehat{\bm{g}} \right )^{+,\tau}[p] 
\end{align}
where 
\begin{align*}
    \widehat{\mathbf{C}}(\mu,  \widehat{\bm{u}}^{k-1}) &= \DualDict^T \mathbf{C}(\mu,  \widehat{\bm{u}}^{k-1}) \PrimalRB  \\
    \widehat{\bm{g}}(\mu,  \widehat{\bm{u}}^{k-1}) &= \DualDict^T\bm{g}(\mu,  \widehat{\bm{u}}^{k-1})
\end{align*}

However, to ensure that inequality constraints~\eqref{eq:sparse_greedy_kkt_nnpr} that prohibits non-negative pressures, it is necessary to not only enrich but also to eliminate terms that violate the constraints in each greedy step, which is simply done by eliminating the largest negative $\widehat{\bm{\lambda}}$ i.e.\ setting the largest negative entry to zero. The elimination index $p_{\texttt{rem}}$ is computed as:
\begin{align}
    p_{\texttt{rem}} &= \argmin_p \left ( \widehat{\bm{\lambda}}^{k} \right )^{-} [p]
\end{align}

\noindent where $
	(\bm{z})^-_i :=\left \lbrace\begin{matrix}
		z_i & \text{if } z_i<0 \\
		0 & \text{otherwise}
	\end{matrix}  \right . $

Though the motivation for this algorithm was based on the OMP to perform the dual dictionary element selection, the greedy active-set algorithm ends up being quite similar to the fixed-point active-set method. In this approach, only a maximum of one vector is activated in a given iteration thereby maintaining a small active set in intermediate iterations. The algorithm is given in \cref{alg:sparse_greedy}, where \crefrange{step:enrich_start}{step:enrich_end} contain the enrichment and the elimination process. In each iterative step, only one of the two operations is performed. Due to the architecture of the algorithm, the only steps that scale with the size of the dictionary are enrichment/elimination \crefrange{step:enrich_start}{step:enrich_end}. The size of mixed system in \cref{step:invert} is not significantly influenced by the size of active set during the intermediate iterations, thanks to the small sizes of active set $\mathcal{I}$ which is closely related to the sparsity of $\widehat{\bm{\lambda}}$.

\begin{algorithm}[!htb]
    \caption{Greedy active-set algorithm} \label{alg:sparse_greedy}
    \begin{algorithmic}[1]
        \State{Input: Queried value of parameter $\mu$}
        \State{Given: Primal basis $\PrimalRB$ and dual Dictionary $\DualDict$}
        \Statex{\hspace{3em} Reduced operators $\mathbf{\PrimalRB}^T \mathbf{K} \mathbf{\PrimalRB} \text{ and } \mathbf{\PrimalRB}^T \bm{f}$} \Comment{can be built offline}
        \State{Initialize: $k = 0, \ \ \mathcal{I} = \emptyset$}
        \While{$\widehat{\bm{u}}$ and $\mathcal{I}$ not converge} \label{step:loop_start}
        \State{Build reduced constraint operators $\widehat{\mathbf{C}}(\mu,  \widehat{\bm{u}}^{k-1})$ and $\widehat{\bm{g}}(\mu,  \widehat{\bm{u}}^{k-1})$} \label{step:build_operators}
        \State{Solve the linear system~\eqref{eq:sparse_greedy_system}} \label{step:invert}
        \State{Set $\widehat{\bm{\lambda}}[\mathcal{I}] \gets \widehat{\bm{\lambda}}_{\mathcal{I}}$ and $\widehat{\bm{\lambda}}[\mathcal{I}^{\mathsf{c}}] \gets \bm{0}$} \Comment{${\mathcal{I}}^{\mathsf{c}}$ is the complementary set} 
        \Statex{}
        \Statex{\hspace{1em} Compute the largest constraint violations and enrich/eliminate accordingly:}
        \If{$ \widehat{\bm{\lambda}}^{k}  \geq \bm{0}$ :} \label{step:enrich_start}
        \Statex{\centerline{$ \begin{array}{rcl}
                p_{\texttt{add}} &=& \argmax_p  \left ( \widehat{\mathbf{C}} \PrimalRB \widehat{\bm{u}}^{k} - \widehat{\bm{g}} \right )^{+,\tau}  \\
                \mathcal{I} &\gets&  \mathcal{I} \cup \{ p_{\texttt{add}}  \}   
        \end{array} $ }} 
        \Else:
        \Statex{\centerline{$ \begin{array}{rcl}
                p_{\texttt{rem}} &=& \argmin_p \left ( \widehat{\bm{\lambda}}^{k} \right )^{-} [p] \\
                \mathcal{I} &\gets&  \mathcal{I} \setminus  \{ p_{\texttt{rem}} \}
        \end{array} $ }}
        \EndIf{} \label{step:enrich_end}
        \State{$k \gets k+1$}
        \EndWhile{} \label{step:loop_end}
        \State{Reconstruct $\bm{u} = \PrimalRB \widehat{\bm{u}}$ and $\bm{\lambda} = \DualDict \widehat{\bm{\lambda}}$}
        \State{Output: $\bm{u}, \bm{\lambda}$}
    \end{algorithmic}
\end{algorithm}

\subsubsection{Application to the Hertz problem}\label{sec:sparse_hertz}
This section deals with the implementation of the greedy active-set algorithm for the Hertz contact problem discussed in \cref{sec:hertz} of previous chapter (\cref{fig:hertz_half_cyl}). The parametrization also remains same, i.e.\ the imposed displacement on the top cylinder $d \in (0,0.3)$. The finite element model also remains the same. \\

\begin{figure}[!htb]
	\centering
    \includegraphics[width=0.5\linewidth]{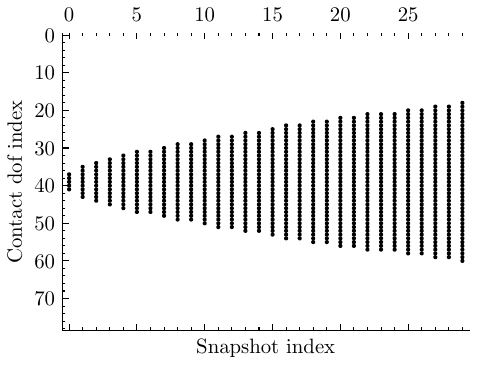}
    \caption[A spy pattern of dual dictionary with 30 elements for the Hertz problem.]{A spy pattern of dual dictionary with 30 elements for the Hertz problem. The number of rows, 79, is same as the number of dofs on the potential contact surface. The dictionary columns are arranged in the increasing order of loading parameter $d$.}
	\label{fig:hertz_dictionary}
\end{figure}

\noindent \textbf{Reduced model}: Primal and dual snapshots on three training sets of sizes 12, 30, 60 and 120 are computed in the offline stage. The parametric points in each training set $\mathcal{P}_\texttt{tr}$ are distributed uniformly in parametric space $\mathcal{P} = ]0,0.3]$. The primal basis is computed using SVD and truncated using cumulatively normalized singular values using a tolerance $\delta$. For the dual variable, a  dictionary consisting contact pressure snapshots is used instead of a reduced basis. A spy plot of the dual dictionary for the training set of size 30 is shown in \cref{fig:hertz_dictionary} for visualization purposes. A validation set $\mathcal{P}_\texttt{val}$ with 119 points, which are the mid-points of the 4th training set is used to compute reconstruction errors.

\begin{figure}[!htb]
	\centering
	\begin{subfigure}[t]{0.48\linewidth}
        \centering
        \includegraphics[width=1\linewidth]{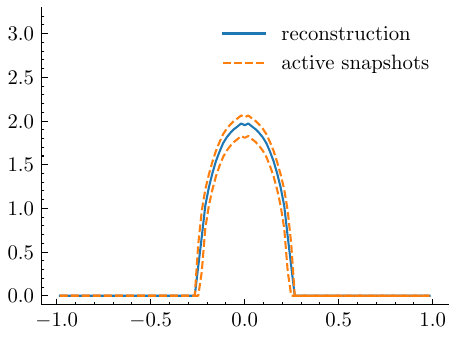}
        \caption{{\tabular[t]{@{}l@{}}$d = 0.14$\\ Recon.\ error $1.6\times 10^{-2}$\endtabular}}
	\end{subfigure}
	\begin{subfigure}[t]{0.48\linewidth}
        \centering
        \includegraphics[width=1\linewidth]{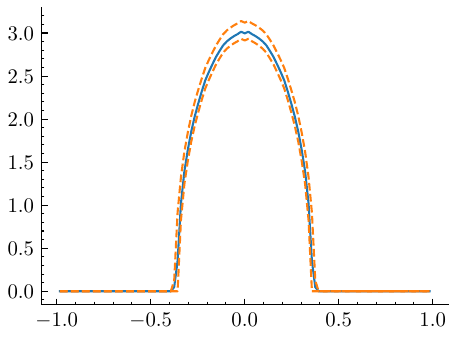}
        \caption{{\tabular[t]{@{}l@{}}$d = 0.26$\\ Recon.\ error $2.5\times 10^{-3}$\endtabular}}
	\end{subfigure}
    \caption[Greedy active-set reconstruction of contact pressure for certain parametric instances of the Hertz problem]{Greedy active-set reconstruction of contact pressure for certain parametric instances outside the training set of the Hertz problem. The elements of dictionary chosen by the algorithm are shown and the reconstruction errors are given. Dictionary of size 12 is used in these examples.}\label{fig:hertz_reconstruction_examples}
\end{figure}

Two instances of contact pressure reconstructions by the greedy active-set algorithm for parametric points outside the training set are shown in \cref{fig:hertz_reconstruction_examples}. The dictionary snapshots selected by the algorithm for the reconstruction are also shown. As expected, only a few dictionary elements, two in this case, are selected. The two selected snapshots are the ones whose contact position is similar to the required contact position for the specific instance. It indicates the greedy enrichment and elimination algorithm is selecting the correct snapshots for reconstruction. 

\begin{figure}[!htpb]
	\centering
    \begin{subfigure}[b]{1.0\linewidth}
        \centering
        \includegraphics[width=0.98\linewidth]{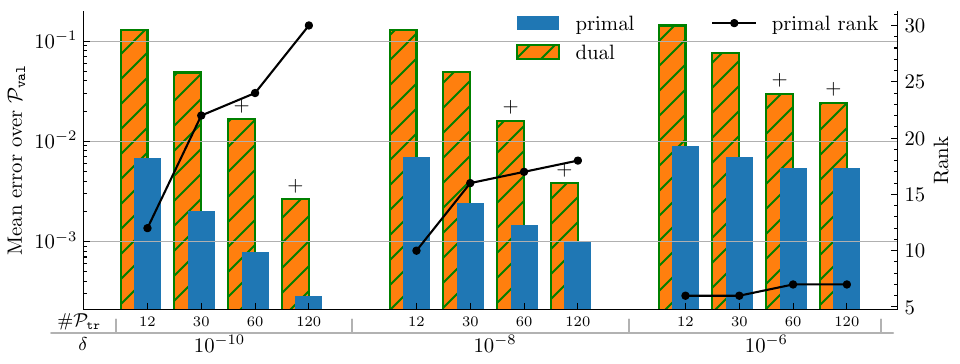}
        \caption{Mean errors (bar plot, left $y$-axis) and primal rank (line plot, right $y$-axis)\label{fig:hertz_mean_error_valset}}
    \end{subfigure}
    \begin{subfigure}[b]{0.48\linewidth}
        \centering
        \includegraphics[width=1.0\linewidth]{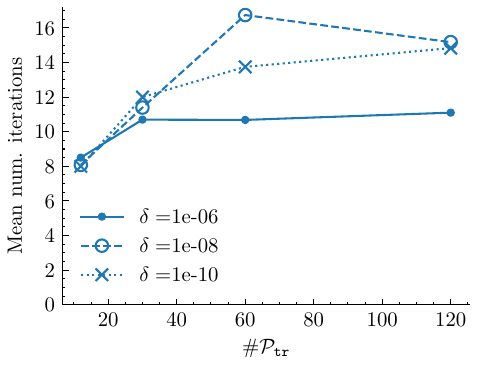}
        \caption{Mean number of iterations\label{fig:hertz_num_iters}} 
    \end{subfigure}
    \begin{subfigure}[b]{0.48\linewidth}
        \centering
        \includegraphics[width=1.0\linewidth]{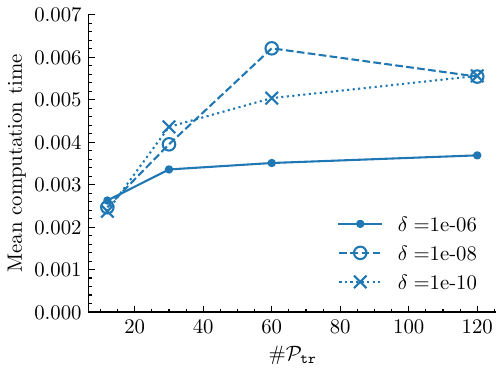}
        \caption{Mean computation time\label{fig:hertz_computation_time}} 
    \end{subfigure}
    \begin{subfigure}[b]{0.48\linewidth}
        \centering
        \includegraphics[width=1.0\linewidth]{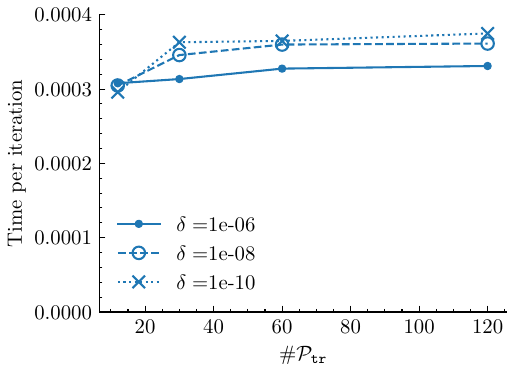}
        \caption{Mean computation time per iteration\label{fig:hertz_time_per_iter}} 
    \end{subfigure}
    \caption[Reconstruction errors, primal rank, number of iterations and computation time vs.\ training set size for the Hertz problem]{Evolution of \textbf{(a)} Mean reconstruction (relative) errors and primal rank, \textbf{(b)} Mean number of iterations, \textbf{(c)} Mean computation time and \textbf{(d)} Mean computation time per iteration with training set size for the Hertz problem. The mean reconstruction error is computed by taking mean of reconstruction errors over the validation set $\mathcal{P}_\texttt{val}$. Primal and dual errors are computed using $\mathcal{H}^1$ and $\mathcal{L}^2$ norm, respectively. Curves correspond to different truncation tolerance  $\delta$ for the primal basis, for values $10^{-6}, \ 10^{-8}$ and $10^{-10} $. Note that computation time excludes time for construction of non-linear operators.}
	\label{fig:hertz_mean_error_and_rank}
\end{figure}

Mean reconstruction errors over the validation set vs.\ size of the training set are shown in \cref{fig:hertz_mean_error_valset}, along with the rank of primal basis. The reconstruction errors improve by more than an order of magnitude as dictionary size is increased from $12$ to $120$ for the case $\delta=10^{-10}$. However, such gains with over-complete dictionaries diminish for the case $\delta=10^{-6}$, especially the primal error improves only slightly. This is probably because at $\delta=10^{-6}$, the primal truncation starts to become dominant source of error, rather than just the linear inseparability of contact pressure which seems to the case for $\delta=10^{-10}$.

The computational time is computed as the total time taken for computation of \cref{alg:sparse_greedy}, excluding the time taken in \cref{step:build_operators}. The operator construction time is excluded because their efficient construction is not a part of this work. Mean computational time on the validation set (\cref{fig:hertz_computation_time}) shows an initial increase with the dictionary size, but later settles down, but it is evident that the computational time is mostly influenced primarily by the number of iterations (\cref{fig:hertz_num_iters}). This can also be seen in \cref{fig:hertz_time_per_iter}, where the computational time per iteration is more or less flat w.r.t.\ the size of dictionary. The small increases in computational time per iteration with dictionary size can be explained due to the increase in the dimension of primal basis $\PrimalRB$, especially the case of $\delta=10^{-10}$ where the primal rank is strongly dependent on the training set size.

\begin{figure}[!htpb]
	\centering
    \begin{subfigure}[t]{1.0\linewidth}
        \centering
        \includegraphics[width=0.48\linewidth]{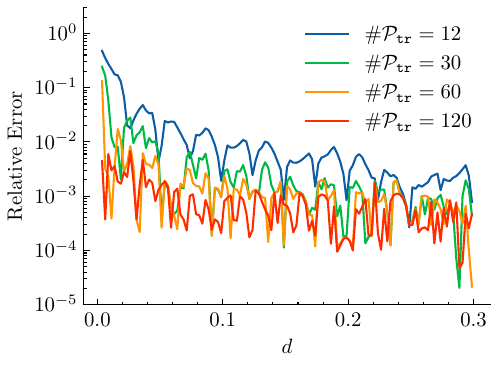}
        \quad
        \includegraphics[width=0.48\linewidth]{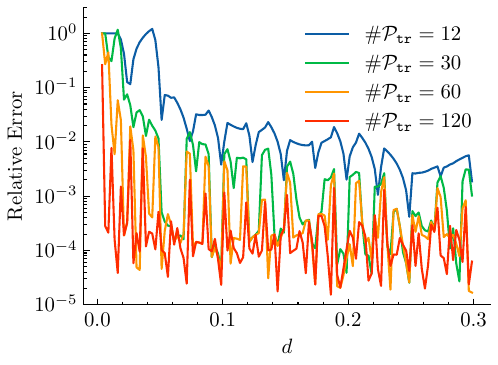}
        \caption{Primal and dual errors for $\delta = 10^{-10}$}
    \end{subfigure}
    \begin{subfigure}[t]{1.0\linewidth}
        \centering
        \includegraphics[width=0.48\linewidth]{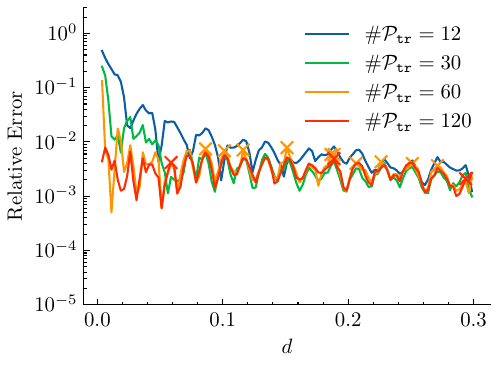}
        \quad
        \includegraphics[width=0.48\linewidth]{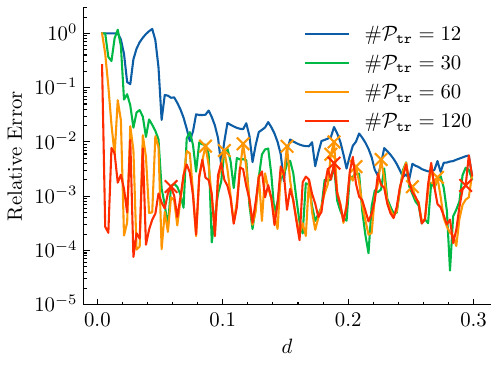}
        \caption{Primal and dual errors for $\delta = 10^{-8}$}
    \end{subfigure}
    \begin{subfigure}[t]{1.0\linewidth}
        \centering
        \includegraphics[width=0.48\linewidth]{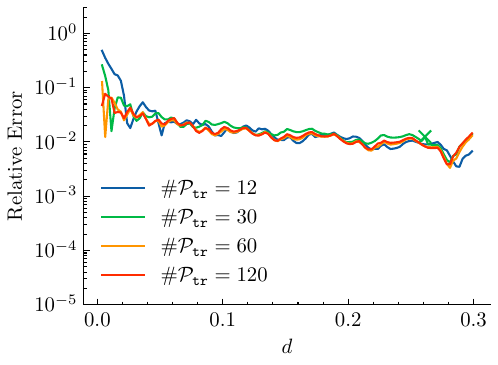}
        \quad
        \includegraphics[width=0.48\linewidth]{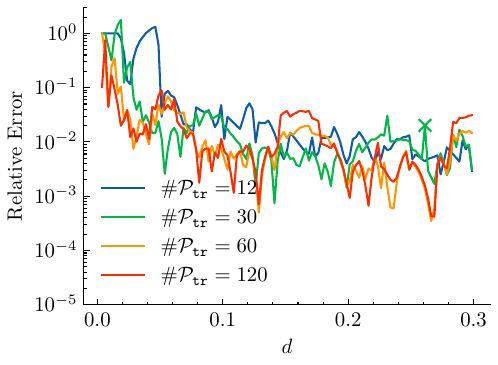}
        \caption{Primal and dual errors for $\delta = 10^{-6}$}
    \end{subfigure}
    \caption[Detailed primal and dual reconstruction errors for the Hertz problem]{Detailed primal and Dual reconstruction errors for the Hertz problem using the training sets 1,2 and 3. Validation set is common to all curves, consisting of 60 equispaced points in $\mathcal{P}$. Each subfigure corresponds to a different truncation tolerance $\delta$ for the primal basis. Primal and dual errors are computed using $\mathcal{H}^1$ and $\mathcal{L}^2$ norm, respectively. Cross-marks $\bm{\times}$ indicate points where the greedy active-set algorithm did not converge to the defined tolerance $10^{-5}$}
	\label{fig:hertz_error_valset}
\end{figure}

Detailed plots of reconstruction errors over the validation set $\mathcal{P}_\texttt{val}$ are shown in \cref{fig:hertz_error_valset} for primal truncation tolerance $\delta \in \{10^{-10}, 10^{-8}, 10^{-6}\}$. Each curve corresponds to one of the training sets. Naturally, a larger training set is expected to have a lower (or at least the same) reconstruction error level. This anticipation is marginally satisfied in the case of $\delta=10^{-10}$, but not in the case of $\delta=10^{-6}$, where there are certain regions where the dual error for larger training sets is higher than a smaller training set. This happens due to the spurious penetrations discussed previously (and in \cref{ch:tau_delta}) and implies that even after relaxing the constraint by value of $\tau$ (see~\eqref{eq:sparse_greedy_kkt_nonpen}), some spurious selections are still possible for higher values of $\delta$. Also, for all cases, the errors tend to be larger in general near $d=0.0$, as the contact is sensitive to the parameter $d$ around this region.

\begin{figure}[!htb]
	\centering
    \includegraphics[width=0.5\linewidth]{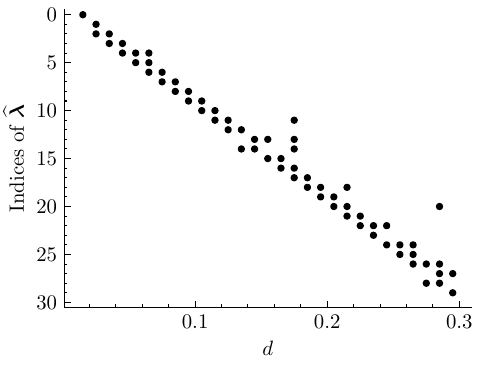}
    \caption[Sparsity of dual reduced dofs for the validation set of Hertz problem]{Sparsity of dual reduced dofs $\widehat{\bm{\lambda}}$ selected by greedy active-set method for reconstructions in validation set of the Hertz problem using dual dictionary of size 30. Dots indicate the non-zero positions.}
	\label{fig:hertz_dualDof_sparsity}
\end{figure}

Sparsity pattern of the dual reduced dofs $\widehat{\bm{\lambda}}$  (i.e.\ the coefficients of the dual dictionary $\DualDict$) resulting from the greedy active-set approach for the training set of size 30 is shown in \cref{fig:hertz_dualDof_sparsity}. As the snapshots in $\DualDict$ are arranged in increasing order of $d$, most of the columns of the sparsity pattern have two close non-zeros, indicating that nearest two snapshots to the targeted reconstruction were chosen by the algorithm. Spurious selections are seen in some rows, but the values of corresponding coefficients are quite small and do not have significant impact on solution.

\subsubsection{Application to the Ironing problem}\label{sec:sparse_ironing}
Unlike the Hertz problem, other contact problems might exhibit larger changes in the contact zone, like the ironing problem from \cref{sec:ironing} (\cref{fig:ironing}). The effects of the strong variations of contact zone on the inseparability of the contact pressure were also discussed. Therefore, over-complete dictionaries might be more significant in generating an effective reduced model.

The parametrization of the horizontal displacement $d_x$ of the iron over the slab is considered for the reduced model. The vertical displacement $d_y = 0.3$ is a fixed value. Given length of the slab $L$, the parametric space is defined by $\mathcal{P}:= \{d_x\ |\ 0\leq d_x\leq L\}$. To study the influence of size of training sets on the reconstruction quality, training sets $\mathcal{P}_{\texttt{tr}}$ are built in a nested manner which allows studying the evolution of reconstruction errors as a function of training set size i.e.\ size of the over-complete dictionary. Each nested level has uniformly distributed points over the parametric space $\mathcal{P}$. The number of points is nearly doubled in the next nested level which also contain the points from the current level. Therefore, at the $n$-th nested level, $2^n+1$ snapshots are computed, and the contact pressure snapshots constitute the dual dictionary $\DualDict$. The primal basis, however, is truncated by filtering the singular values by a fixed tolerance $\delta$, and therefore its rank must be less than $2^n$. Training sets of nested levels $n=3$ to $n=7$ are used to build the primal RB and the dual dictionary. The validation set $\mathcal{P}_{\texttt{val}}$ contains all points from nested level $n=8$ that are not present in level $n=7$. 

\begin{figure}[!htb]
	\centering
	\begin{subfigure}[t]{0.48\linewidth}
        \centering
        \includegraphics[width=1\linewidth]{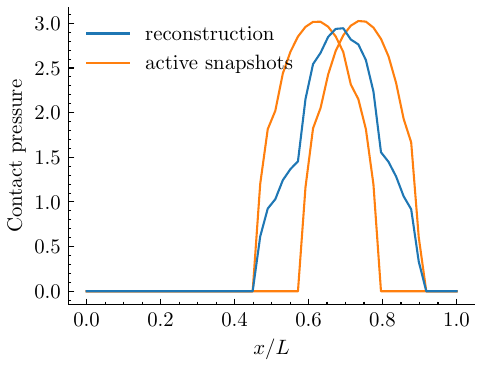}
        \caption{$d_x = 0.68 L$, nested level 3}
	\end{subfigure}
	\begin{subfigure}[t]{0.48\linewidth}
        \centering
        \includegraphics[width=1\linewidth]{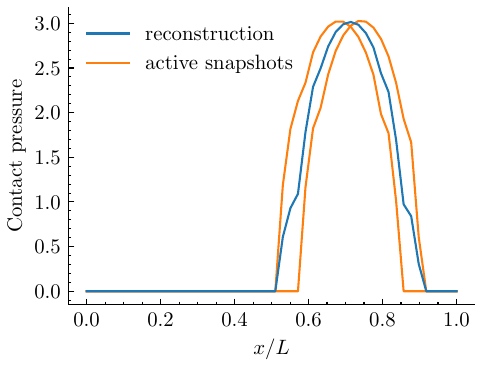}
        \caption{$d_x = 0.65 L$, nested level 4}
	\end{subfigure}
    \caption[Greedy active-set reconstruction of contact pressure for certain parametric instances of the ironing problem.]{Greedy active-set reconstruction of contact pressure for certain parametric instances of the ironing problem. The dictionary nested level, dictionary elements chosen by the algorithm and the reconstruction errors are also shown.}
	\label{fig:reconstruction_example_ironing}
\end{figure}

\begin{figure}[!htb]
	\centering
    \includegraphics[width=0.5\linewidth]{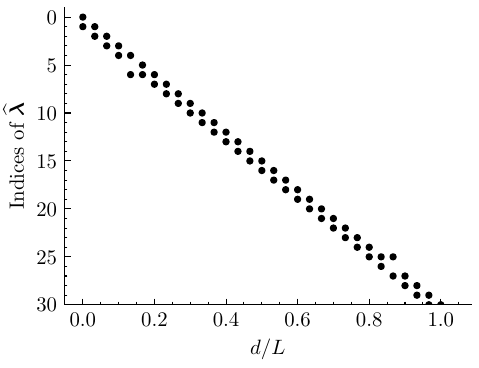}
    \caption[Sparsity of dual reduced dofs selected by greedy active-set method for the ironing problem]{Sparsity of dual reduced dofs $\widehat{\bm{\lambda}}$ selected by greedy active-set method for various reconstructions of ironing problem using dual dictionary of nested level 5}
	\label{fig:ironing_dualDof_sparsity}
\end{figure}

\cref{fig:reconstruction_example_ironing} shows two instances of reconstruction of contact pressure along with the dictionary elements that were chosen by the algorithm. The algorithm selects snapshots that are closest to the required contact position. This can also be seen in \cref{fig:ironing_dualDof_sparsity}, where the sparsity of dual dofs shows a clear selection pattern. Unlike the sparsity pattern of Hertz problem (\cref{fig:hertz_dualDof_sparsity}), the ironing problem does not display spurious snapshot selections; possibly due to the stronger dependence of contact position of the ironing problem; unlike Hertz problem where contact pressure curves were centered around the same position. Therefore, the activation process is less ambiguous and less sensitive to primal truncations for the case of ironing problem.

\begin{figure}[!htpb]
	\centering
    \begin{subfigure}[b]{1.0\linewidth}
        \centering
        \includegraphics[width=0.98\linewidth]{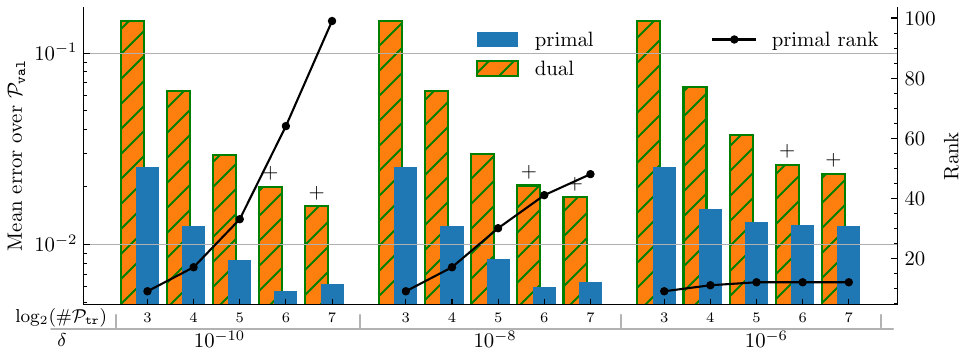}
        \caption{Mean reconstruction errors and primal rank}
	\end{subfigure}
	\begin{subfigure}[b]{0.48\linewidth}
        \centering
        \includegraphics[width=1\linewidth]{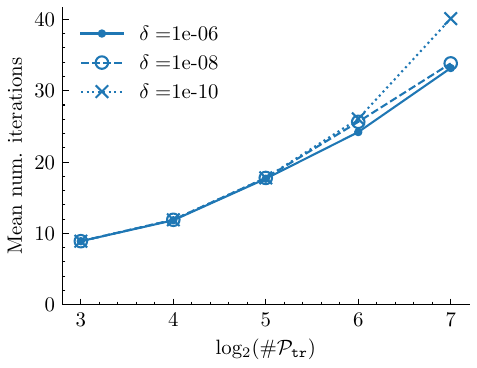}
        \caption{Mean number of iterations}
        \label{fig:ironing_hlevel_evolution_iters}
	\end{subfigure}
	\begin{subfigure}[b]{0.48\linewidth}
        \centering
        \includegraphics[width=1\linewidth]{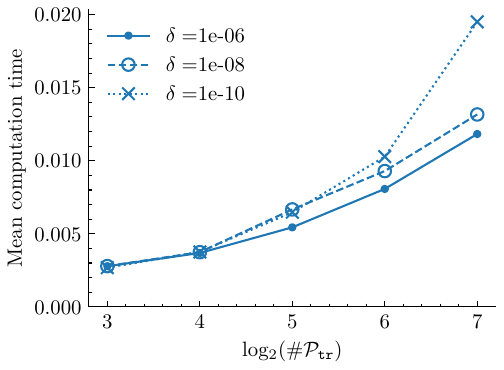}
        \caption{Mean computation time}
        \label{fig:ironing_hlevel_evolution_time}
	\end{subfigure}
	\begin{subfigure}[b]{0.48\linewidth}
        \centering
        \includegraphics[width=1\linewidth]{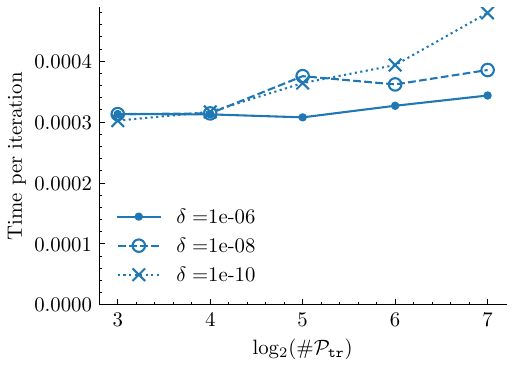}
        \caption{Mean computation time per iteration}
        \label{fig:ironing_hlevel_time_per_iter}
	\end{subfigure}
    \caption[Reconstruction errors, primal rank, number of iterations, and computation time vs.\ training set size for the Ironing problem]{Evolution of \textbf{(a)} mean reconstruction (relative) errors and rank of primal basis \textbf{(b)} mean number of iterations, \textbf{(c)} mean computation time and \textbf{(d)} mean computation time per iteration with training set size for the Ironing problem. The mean reconstruction error is computed by taking mean of reconstruction errors over the validation set $\mathcal{P}_\texttt{val}$. Primal and dual errors are computed using $\mathcal{H}^1$ and $\mathcal{L}^2$ norm, respectively. Curves correspond to different truncation tolerance  $\delta$ for primal basis, for values $10^{6}, \ 10^{-8}$ and $10^{-10}$. Note that computation time excludes time for construction of non-linear operators. \label{fig:ironing_hlevel_evolution}}	
\end{figure}

\cref{fig:ironing_hlevel_evolution} shows the evolution of reconstruction (relative) errors, primal basis rank, number of iterations and computation time over the validation set $\mathcal{P}_\texttt{val}$ {vs.} training set nested levels ($\text{log}_2 (\# \mathcal{P}_{\texttt{tr}})$). Means of the reconstruction related quantities over the validation set is shown for primal truncation tolerances of $10^{-6}$, $10^{-8}$ and $10^{-10}$. The dual error improves, although slowly, with increasing size of dictionary, but the primal error shows signs of stagnation near level $6$ and $7$ and even earlier for the case of $\delta=10^{-6}$. Errors for $\delta=10^{-8}$ and $10^{-10}$ are almost the same, possibly because the inseparability of contact pressure and the finite size of dual dictionary are important sources of error. Conversely, in case $\delta=10^{-6}$, the primal truncation seems to have more influence as primal error stagnates early, similar to the observations in case of the Hertz problem.

Moreover, the number of iterations and consequently the computation time seems to increase significantly with the dictionary size, in contrast to the Hertz problem where this relation was comparatively flat. A possible explanation for this behaviour could be that the larger dictionaries of the ironing problem, coupled with stronger parametric dependence of contact pairs, forces the greedy active-set to perform more iterations to find the appropriate dual dictionary candidates and  contact pairs. Mean computation time per iteration shows trends that are a bit similar to the primal rank i.e. both primal rank and computational time per iteration are flat for $\delta=10^{-6}$, whereas both are monotonically increasing for $\delta=10^{-10}$. This was an expected observation since a large part of the system~\eqref{eq:sparse_greedy_system} is related to the dimensionality of the primal basis $\PrimalRB$. The size of the dual dictionary, even though not truncated, does not seem to have any significant influence on the computation time per iteration.

\subsection{A non-penetrating convex hull approach for monolithic dictionaries}
\label{sec:convex_hull}
As discussed in \cref{rem:div_free}, an interesting feature of mixed ROMs for incompressible flow problems is the computation of divergence-free RB\@. The advantage of this feature is that incompressibility constraint does not need to be explicitly enforced in the reduced problem. Moreover, orthogonality of a divergence-free velocity subspace and the irrotational pressure gradient term causes the pressure terms in the reduced problem to completely vanish (see~\cite{Veroy2005,Liberge2010} for details). Computation of a divergence-free RB does not require any special treatment, as the divergence-free property of the training set snapshots is easily preserved in the resulting reduced subspace. 

In this section, we attempt to extend a similar notion of constraint satisfying solution space to contact mechanics. But the approach is quite different compared to the incompressible flow problem because the nature of constraint is different. The incompressibility constraint, an equality constraint, applies globally over the domain; whereas the non-penetration constraint, an inequality constraint, applies locally on the contact surfaces and can be either active or inactive locally. Also, the displacement and contact pressure do not possess any special mathematical properties that cause the pressure term to vanish, as in the case of incompressible flow. In fact, the approach presented in this section, a reduced solution ``space'' is defined that naturally satisfies the non-penetration condition, but does not result in elimination of the contact pressure. 

Unlike the greedy active-set algorithm previously explored which used POD basis for displacement and a dictionary for contact pressure, the following approach utilises a monolithic dictionary whose columns are snapshots consisting of both displacement and contact pressure. The reason for using a monolithic dictionary will be become evident in the formulation of this approach. Moreover, this method allows the residual to be projected on a low-dimensional space like it was done in the unconstrained problem in \cref{sec:dictionary_resolution}. It also permits the utilisation of sparse regression methods discussed in \cref{sec:sparse_methods}, provided the non-negativity constraint can be enforced. In this work, the non-negative version of the FOCUSS algorithm is used.

Consider a contact problem where the contact pairs are constant and independent of the parameter of the problem. The discrete inequality constraint can be expressed using matrices that do not change with the displacement field. The feasible region $\mathcal{K}$ in the displacement solution space can be defined as:
\begin{align}
    \mathcal{K} := \{ \bm{u} \in \PrimalSpace \ | \ \mathbf{C} \bm{u} - \bm{g} \leq 0 \}
    \label{eq:linear_ineq_constraint}
\end{align}

\begin{figure}[!htb]
    \centering
    \includegraphics{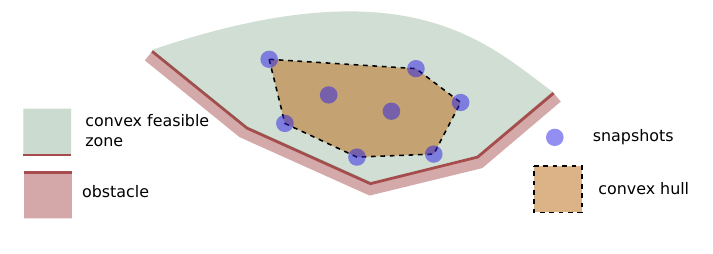}
    \caption[Illustration of non-penetrating property of convex hull]{Illustration of non-penetrating property of convex hull. The non-penetration constraint in~\eqref{eq:linear_ineq_constraint} is linear, convex and is made of segments corresponding to each element of the constraints $\mathbf{C} \bm{u} = \bm{g}$. In a high-dimensional problem, each segment of the constraint can be thought of as a hyperplane. The convex hull, which is simply a convex enclosure of snapshots, can also be thought of as a high-dimensional set.}
    \label{fig:convex_hull_demo}
\end{figure}

As the operators $\mathbf{C}$ and $\mathbf{g}$ are linear, it is easy to demonstrate the convexity of $\mathcal{K}$. This would mean that a convex combination\footnote{A convex combination is a linear combination with non-negative coefficients that sum up to unity} of non-penetrating snapshots would naturally satisfy the inequality constraint. The implication here is that if the reduced model explores the convex hull\footnote{Convex hull of a set of points (or snapshots) is the set containing all convex combinations of the points. In other words, it is also the smallest convex set that envelopes the set of points.} of the training set snapshots, the inequality constraint need not be explicitly enforced. This notion is visually shown in \cref{fig:convex_hull_demo}, where the convex hull of training set snapshots does not violate the convex non-penetration constraint. This property of convexity can be exploited to build efficient reduced models, if the following hypothesis holds: 

\begin{hypothesis}[Convex Subset hypothesis]\label{hyp:lowdim_convex_subset}
    Given a parametrized high-dimensional inequality-constrained problem with a convex feasible set $\mathcal{K}$, its solutions lie in a low-dimensional convex subset inside $\mathcal{K}$.
\end{hypothesis}

The Convex Subset hypothesis is an extension of the idea behind the low-rank hypothesis to inequality-constrained problems with convex feasible regions. This hypothesis allows us to explore a solution set limited to the convex hull defined by the snapshots without worrying about the non-penetration conditions. Note that the convexity of feasible region assures that any solution in the convex hull of displacement snapshots still lie inside the feasible region, independent of the \cref{hyp:lowdim_convex_subset}. Therefore, in this approach based on the convex hull of the dictionary, the displacement field candidates built using convex combinations of the primal dictionary $\PrimalDict$ will be explored:
\begin{align*}
    \bm{u} &\approx \PrimalDict \widehat{\bm{u}} \\
    \text{s.t. } \bm{1}^T \widehat{\bm{u}} &= 1 \\
    \widehat{\bm{u}} & \geq \bm{0}
\end{align*}

\noindent As the non-penetration constraint will be satisfied naturally by all candidates in the convex hull of dictionary snapshots, a dictionary-based approximation can be computed using coefficients $\widehat{\bm{u}}$ that satisfy equilibrium equations. Therefore, the objective is to minimize the following residual:
\begin{align*}
    \bm{r}(\widehat{\bm{u}}, \widehat{\bm{\lambda}}) = \mathbf{K}(\mu) \PrimalDict \widehat{\bm{u}} + \mathbf{C}^T(\mu, \bm{u}) \DualDict \widehat{\bm{\lambda}} - \bm{f}(\mu)
\end{align*}

\noindent As mentioned before, a monolithic dictionary will be used in this problem, in which each column consists of contact pressure snapshot stacked under the displacement snapshot. This means the coefficients $\widehat{\bm{u}}$ and $\widehat{\bm{\lambda}}$ will be same and will be denoted by $\bm{\alpha}$ here onwards. Use of monolithic dictionary not only reduces the number of unknowns, but also makes the problem of minimizing the residual $\norm{\bm{r}(\bm{\alpha}, \bm{\alpha})}{}$ well-defined. Moreover, the residual can be projected on a low-dimensional subspace of the $\text{colsp}(\PrimalDict)$, described by the matrix $\mathbf{B}$ containing its basis, while also adding sparsity constraints to the unknown $\bm{\alpha}$, like in~\eqref{eq:sparse_res_randomized}. Thus, the problem can be stated as:

\begin{subequations}
    \begin{align}
        \min \ & \norm{\bm{\alpha}}{p} \\
        \text{s.t. } & \norm{\mathbf{B}^T \bm{r}(\bm{\alpha}, \bm{\alpha})}{} < \varepsilon \\
            &  \bm{1}^T\bm{\alpha} = 1 \label{eq:convex_sum} \\
            &  \bm{\alpha} \geq \bm{0}
    \end{align} \label{eq:convex_comb_equilibrium}
\end{subequations}

As the constraint~\eqref{eq:kkt_discrete_nonpen} is not imposed directly in this approach, the non-penetration condition is satisfied only because of the convexity of the feasible region. Thus, the equality condition $\mathbf{C}_{\ActiveSet}(\mu,\bm{u}) \bm{u} = \bm{g}_{\ActiveSet}(\mu,\bm{u})$ for active constraints is also not imposed explicitly. Consequently, the complementary slackness KKT condition~\eqref{eq:kkt_discrete_comp_slack} may not be satisfied exactly. This will be seen in the upcoming numerical examples. \\

\noindent The \textbf{nnFOCUSS} Algorithm: To solve the convex hull approximation problem~\eqref{eq:convex_comb_equilibrium}, the non-negative version of the FOCUSS, the nnFOCUSS algorithm~\cite{Lauzeral2019,Manucci2022} is used. The nnFOCUSS algorithm works by computing an appropriate relaxation parameter for each iterative update that maintains non-negativity. After computation of new coefficients at \cref{step:focuss_alpha} of the \cref{alg:focuss}, the relaxation step of nnFOCUSS algorithm is performed as follows: \\

\begin{algorithmic}[0]
    \If{$\min(\bm{\alpha}^k) < 0$}
    \State{$\Delta \bm{\alpha} = -(\bm{\alpha}^k - \bm{\alpha}^{k-1})^{-}$}
    \State{$\bm{\alpha}^k \gets \bm{\alpha}^{k-1} + \min \left (\frac{\bm{\alpha}^{k-1}}{\Delta \bm{\alpha}} \right ) \Delta \bm{\alpha}$} \Comment{Element wise division}
    \EndIf{}
\end{algorithmic}

\noindent The algorithm must be initialized using a non-negative coefficients $\bm{\alpha}^0$, which is computed using the non-negative least squares~\cite{Lawson1995} solution (which is not sparse, in general). The rest of the algorithm is same as standard FOCUSS\@. The problem~\eqref{eq:convex_comb_equilibrium} can be solved by plugging the following inputs in the nnFOCUSS algorithm:
\begin{align*}
    \mathbf{D} =  \begin{bmatrix}
        \mathbf{B}^T  \mathbf{K}_{\texttt{mono}}\\
        \bm{1}^T
    \end{bmatrix} , \ \ \ \  
    \bm{x} = \begin{bmatrix}
        \mathbf{B}^T \bm{f}(\mu) \\
        1
    \end{bmatrix} , \ \ \ \  
    \bm{\alpha}^0 = \mathtt{nnls}(\mathbf{D}, \bm{x})
\end{align*}
\noindent where $\mathbf{K}_{\texttt{mono}}$ is the operator for the monolithic residual $\bm{r}(\bm{\alpha}, \bm{\alpha})$:
\begin{align*}
    \mathbf{K}_{\texttt{mono}} = \mathbf{K}(\mu) \PrimalDict  + \mathbf{C}^T(\mu, \bm{u}) \DualDict 
\end{align*}
\texttt{nnls} indicates the non-negative least squares. An open source implementation of \texttt{nnls} provided in SciPy package~\cite{scipy} is used.

\begin{figure}[htpb]
    \centering
    \includegraphics[width=0.7\linewidth]{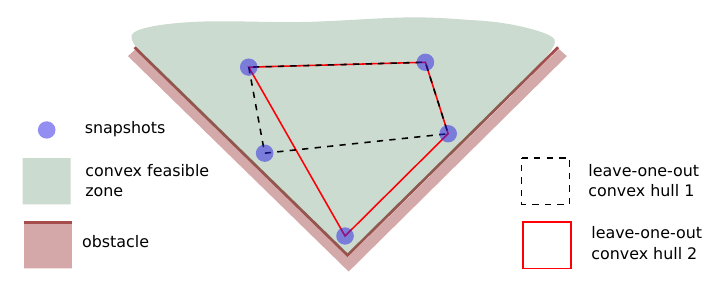}
    \caption[Illustration of the convex subset hypothesis test with a point mass moving in a convex feasible region.]{Illustration of the Test~\ref{test:lowdim_convex_subset} with a point mass moving in a V-shaped obstacle forming a convex feasible region. Two leave-one-out instances and the convex hull of the corresponding $\widetilde{\PrimalDict}$ are shown. $\epsilon_{\texttt{CHLS}}$ is a measure of distance between the left out snapshot and the convex hull of $\widetilde{\PrimalDict}$. The leave-one-out instance 1 has a higher $\epsilon_{\texttt{CHLS}}$ than the instance 2.}
    \label{fig:convex_hypothesis_test}
\end{figure}

In this thesis, analytical validation for the convexity of feasible set and \cref{hyp:lowdim_convex_subset} will not be given. Instead, the following tests, that apply the leave-one-out approach on the training set snapshots, are proposed to analyze if these conditions hold within a reasonable error:
\begin{enumerate}[Test A.]
    \item\label{test:convex_feasible_set} Even if the operators $\mathbf{C}(\mu,\bm{u})$ and $\bm{g}(\mu, \bm{u})$ are non-linear and the feasible region $\mathcal{K}(\mu)$ shows a dependence on the parameter $\mu$, it may be reasonable to assume a nearly convex feasible set if the contact pairs do not change drastically. For such cases, the convexity of feasible set can be checked numerically using the training set snapshots and the corresponding non-linear operators. Given the constraint, $\mathbf{C}(\mu,\bm{u}) \bm{u} -g(\mu,\bm{u}) \leq 0$, the cross penetration (CrPen) error of the $k$-th dictionary element $\bm{d}$ across all constraints encountered in training set is computed as follows:
        \begin{gather}
            \begin{aligned}
                \epsilon_{\texttt{CrPen}}[k] = \max \ [\mathbf{C}(\mu,\widetilde{\bm{d}}) \bm{d} -g(\mu,\widetilde{\bm{d}})]^{+,0} \ \ \forall \widetilde{\bm{d}} \in \widetilde{\PrimalDict}
            \end{aligned}
            \label{eq:crosspen_error}
        \end{gather}
        with $\bm{d} = \PrimalDict[:,k]$ and $\widetilde{\PrimalDict} = \PrimalDict \setminus \bm{d}$. If all values of $\epsilon_{\texttt{CrPen}}$ are sufficiently small, it would be reasonable to assume the feasible set $\mathcal{K}(\mu)$ to be convex. The notation $[\cdot]^{+,0}$  has the same meaning as in~\eqref{eq:sparse_enrich_concept}.
    \item\label{test:lowdim_convex_subset} The \cref{hyp:lowdim_convex_subset} can be verified by approximating the left out snapshot with the convex hull of the rest of the snapshots using least-square criteria. The approximation error, referred to here as the convex hull least square (CHLS) error, of the leave-one-out approach can be used to assess the existence of hypothetical low-dimensional convex set. Each dictionary element $\bm{d}$ is approximated using the convex hull of $\widetilde{\PrimalDict}$ and the error is recorded:
        \begin{gather}
            \begin{aligned}
                \epsilon_{\mathtt{CHLS}}[k] = \min_{\bm{\alpha}} \frac{\norm{\widetilde{\PrimalDict} \bm{\alpha} - \bm{d}}{2}}{\norm{\bm{d}}{2}} \\
                \text{s.t.}  \bm{1}^T \bm{\alpha} = 1
            \end{aligned}
        \end{gather}
        A high value of $\epsilon_{\mathtt{CHLS}}$ might occur due to lack of low-dimensional convexity or due to scattered training set data. The $\epsilon_{\mathtt{CHLS}}$ can be visualized using a low-dimensional problem shown in \cref{fig:convex_hypothesis_test}, where a point mass moves is permitted to move between a V-shaped obstacle. When leave-one-out approach is applied to a set of five snapshots, $\epsilon_{\mathtt{CHLS}}$ is high when the snapshot in the corner is left out. On the other hand, a low value of $\epsilon_{\mathtt{CHLS}}$ assures that it is possible to reconstruct reasonably using convex hull of the dictionary.
\end{enumerate}

\subsubsection{Illustrative example: elastic rope-obstacle problem}
The convex hull approximation is first demonstrated in a problem where constraint operators are constant, so that convexity of feasible space is assured. The elastic rope-obstacle problem form \cref{sec:illustration} is considered, but with a different parametrization. In this case, the elasticity $\nu(x)$ of the rope is a function of parameter $\gamma_1$.
\begin{gather}
	\begin{aligned}
        &\nu(x) \nabla^2 u(x) = f \ \ \text{ on } x \in [0,1] \\
		&u(0) = u(1) = 0 \\
		&u(x) \geq -0.2 (\sin(\pi x) - \sin(3 \pi x)) -0.5 \\
        \text{where} &  \\
        &\nu(x) = \left \{  
            \begin{matrix} \gamma_1  & \text{if } x < 0.5 \\
                                30      & \text{otherwise}
            \end{matrix} \right . \ \ \text{ on } \gamma_1 \in [10,50]
	\end{aligned}
	\label{eq:membrane_elastic_parameter}
\end{gather}

\noindent $
\begin{aligned}
    &\text{Training set:} && \gamma_1 \in \mathcal{P}_{\texttt{tr}}\hspace{0.5em} = \{10, 15, 20 \dots 45, 50\}\\
    &\text{Validation set:} && \gamma_1 \in \mathcal{P}_{\texttt{val}} = \{12.5, 17.5, 22.5 \dots  42.5, 47.5\}
\end{aligned}
$

\begin{figure}[htpb]
    \centering
    \includegraphics[width=0.5\linewidth]{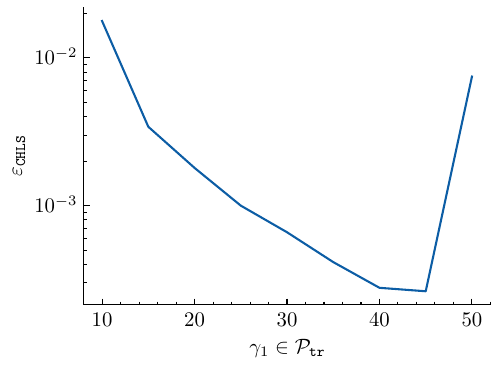}
    \caption[Testing the convex subset hypothesis for the elastic rope-obstacle problem]{Test~\ref{test:lowdim_convex_subset}: testing the \cref{hyp:lowdim_convex_subset} for the given training set of elastic rope-obstacle problem~\eqref{eq:membrane_elastic_parameter}. Test~\ref{test:convex_feasible_set} for convexity of feasible set is not necessary as the constraints are linear.}
    \label{fig:membrane_lowdim_convex}
\end{figure}

The monolithic dictionary with 9 snapshots corresponding to training set defined above is computed in offline stage. The constraint is linear and independent of the parameter, and therefore, the feasible region is convex, which eliminates the necessity of Test~\ref{test:convex_feasible_set}. Hence, only the \cref{hyp:lowdim_convex_subset} is tested using leave-one-out approach (Test~\ref{test:lowdim_convex_subset}). The test results are shown in \cref{fig:membrane_lowdim_convex}, indicating that hypothesis holds to a reasonable level.

\begin{figure}[!htb]
    \centering
    \begin{subfigure}[t]{0.45\linewidth}
        \includegraphics[width=1\linewidth]{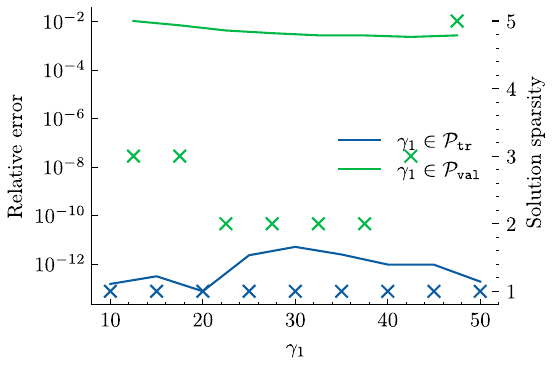}
        \caption{Displacement}
    \end{subfigure}
	\begin{subfigure}[t]{0.45\linewidth}
        \includegraphics[width=1\linewidth]{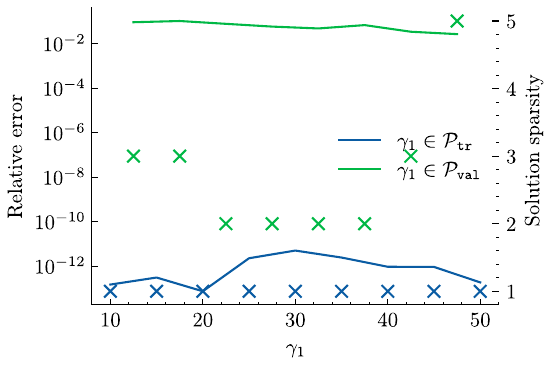}
        \caption{Contact Pressure}
    \end{subfigure}
    \caption[Convex hull reconstruction errors for elastic rope-obstacle problem.]{Convex hull reconstruction errors for elastic rope-obstacle problem with parameter $\gamma_1$ using $\mathbf{B} \gets \text{svd}(\PrimalDict, \delta=10^{-7})$. Crosses $\times$ indicate achieved sparsity (w.r.t.\ right y-axis). Primal and dual sparsity are naturally equal as the dictionary is monolithic.}
    \label{fig:convex_membrane_gamma1}
\end{figure}

The low-rank matrix $\mathbf{B}$ on which the residual will be projected is computed using left singular vectors of the truncated SVD $\mathbf{B} \gets \text{svd}(\PrimalDict, \delta)$, with tolerance $\delta = 10^{-7}$. The reconstruction errors and resultant sparsity of the solution in the training set and a validation set is shown in \cref{fig:convex_membrane_gamma1}. The training set is reconstructed within near numerical precision. In the validation set, the algorithm chooses the training set points that are nearest to $\gamma_1$ to create the best possible reconstruction with the given dictionary. 

As discussed earlier, the KKT conditions, except the non-negativity of contact pressure, is not explicitly imposed. However, in this particular case with perfectly convex feasible region, one would expect that at least the non-penetration would be satisfied accurately. However, the linear combinations are not perfectly convex as the condition~\eqref{eq:convex_sum} is satisfied only to a certain precision. The observed orders of these quantities during reconstruction are given in \cref{tab:convex_membrane_kkt_values}. The validation set reconstructions in fact show penetration and non-zero complementary slack in the order similar as that of the deviation from convex combinations i.e.\ $|\bm{1}^T \bm{\alpha}-1|$.

\begin{table}[!htb]
	\setlength\extrarowheight{2pt}
    \small
	\centering
    \begin{tabular}{lcc}
        \toprule
        &  $\gamma_1 \in \mathcal{P}_{\texttt{tr}}$   & $\gamma_1 \in \mathcal{P}_{\texttt{val}}$  \\
        \midrule
        $|\bm{1}^T \bm{\alpha}-1|$ & $\mathcal{O}(10^{-13})$ & $\mathcal{O}(10^{-3})$ \\
        Penetration & $\mathcal{O}(10^{-16})$ & $\mathcal{O}(10^{-4})$ \\
        Complementary Slack & $\mathcal{O}(10^{-15}) $ & $\mathcal{O}(10^{-3})$\\
        \bottomrule
    \end{tabular}
    \caption{Observed values of various compliances for the elastic rope-obstacle problem with parameter $\gamma_1$.}
	\label{tab:convex_membrane_kkt_values}
\end{table}

\subsubsection{Demonstration on over-complete dictionaries of the Hertz problem}
Hertz problem with same parametrization and training sets from \cref{sec:sparse_hertz} is solved using the convex hull approximation. For projecting the residual on a low-rank subspace, the matrix $\mathbf{B}$ is computed in two different ways:
\begin{align*}
    \mathbf{B} \gets \text{svd}(\PrimalDict, \delta) \ \ \ 
\text{or} \  \ \ \
    \mathbf{B} \gets \text{orth}(\PrimalDict \mathbf{R})
\end{align*}
\noindent where $\mathbf{R}$ is a matrix with random entries from a uniform distribution~\eqref{eq:randomized_orth_matrix}. As discussed in \cref{sec:dictionary_resolution}, this matrix is used to create a smaller matrix by randomly combining dictionary columns.

\begin{figure}[!htb]
    \centering
    \includegraphics[width=0.5\linewidth]{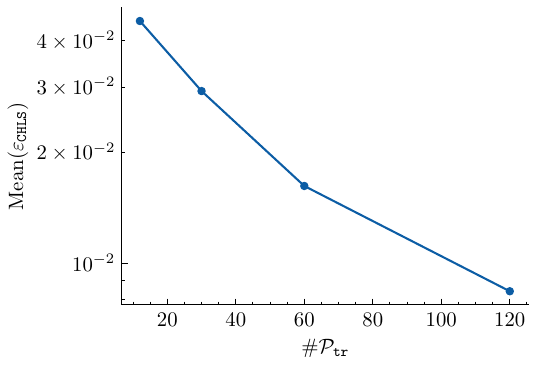}
    \caption[Testing the convex subset hypothesis for the Hertz problem]{Test~\ref{test:lowdim_convex_subset}: testing the \cref{hyp:lowdim_convex_subset} for the various training sets of the Hertz problem}
    \label{fig:hertz_lowdim_convexity}
\end{figure}

As the constraint operators $\mathbf{C}$ and $\bm{g}$ show non-linear dependence w.r.t. $\bm{u}$, it is necessary to test the convexity of the feasible set $\mathcal{K}(\mu)$ (Test~\ref{test:convex_feasible_set}). All CrPen errors ($\epsilon_{\texttt{CrPen}}$) were found to be quite small, in the order of $\mathcal{O}(10^{-7})$ for all combinations of $\bm{d}$ and $\widetilde{\bm{d}}$ in~\eqref{eq:crosspen_error}. Test~\ref{test:lowdim_convex_subset} is performed to check the \cref{hyp:lowdim_convex_subset}, with results shown in \cref{fig:hertz_lowdim_convexity}. The reconstruction errors for the reduced model can be expected to be more than the $\epsilon_{\texttt{CHLS}}$ shown in the figure. 

\begin{figure}[!htb]
    \centering
    \begin{subfigure}[t]{1.0\linewidth}
        \includegraphics[width=0.98\linewidth]{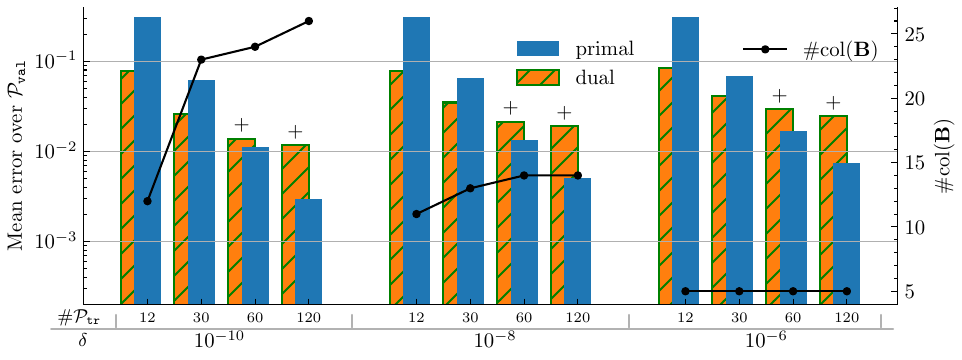}
        \caption{$\mathbf{B} \gets \text{svd}(\PrimalDict, \delta)$}
    \end{subfigure}
    \begin{subfigure}[t]{1.0\linewidth}
        \includegraphics[width=0.98\linewidth]{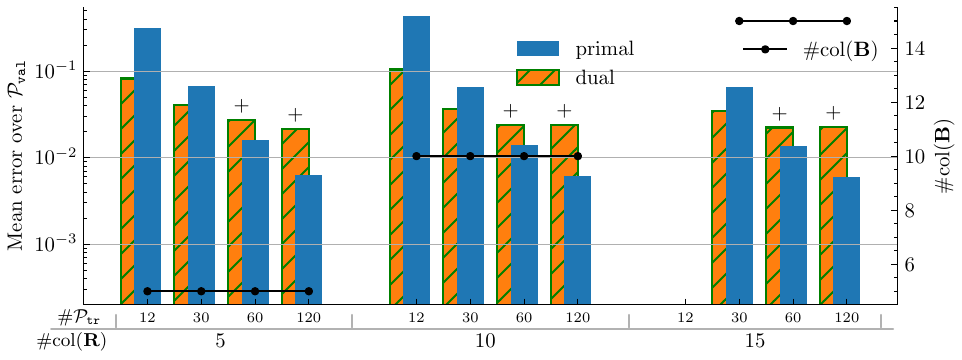}
        \caption{$\mathbf{B} \gets \text{orth}(\PrimalDict \mathbf{R})$}
    \end{subfigure}
    \caption[Convex hull reconstruction errors for the Hertz problem.]{Convex hull reconstruction errors for the Hertz problem using different dimensions of $\mathbf{B}$. The penetration in all reconstruction cases is observed to be in order of $10^{-8}$. The case of $\text{rank}(B)= 15$ and $\# \mathcal{P}_{\texttt{tr}}$ is not possible, hence not shown in \textbf{(b)}. Primal and dual errors are computed using $\mathcal{H}^1$ and $\mathcal{L}^2$ norm, respectively.}
    \label{fig:convex_hertz_recon_err}
\end{figure}

For large dictionaries, it makes sense to compare the performance of the method using various $\delta$ and sizes of $\mathbf{R}$. The mean reconstruction errors over the validation set vs size of dictionaries can be seen in \cref{fig:convex_hertz_recon_err}. The use of randomized matrices $\mathbf{R}$ naturally generates more error than truncated SVD, as the truncated SVD selects the best representative subspace for the given training set. Similar to the greedy active-set method of \cref{alg:sparse_greedy}, the improvement in error flattens out for large dictionaries. On the other hand, the number of iterations and computation time seems to be nearly independent of the dictionary size.

\begin{figure}[!htb]
    \centering

    \begin{subfigure}[t]{0.45\linewidth}
        \includegraphics[width=1\linewidth]{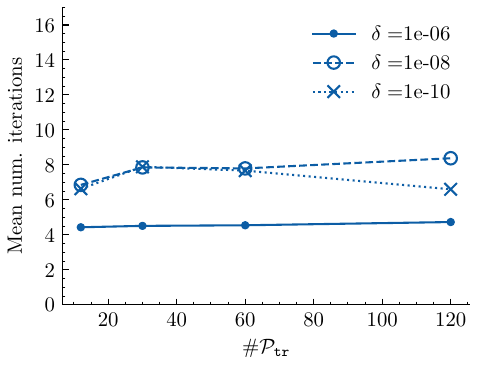}
        \caption{Num.\ iterations, $\mathbf{B} \gets \text{svd}(\PrimalDict, \delta)$}
    \end{subfigure}
	\begin{subfigure}[t]{0.45\linewidth}
        \includegraphics[width=1\linewidth]{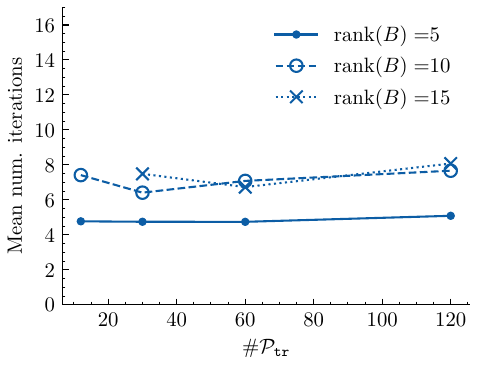}
        \caption{Num.\ iterations, $\mathbf{B} \gets \text{orth}(\PrimalDict \mathbf{R})$}
    \end{subfigure}\\[1em]
    \begin{subfigure}[t]{0.45\linewidth}
        \includegraphics[width=1\linewidth]{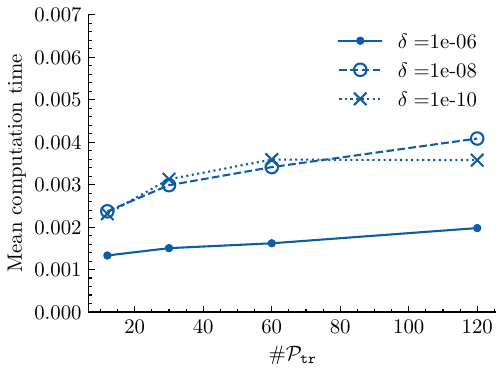}
        \caption{Comp. Time, $\mathbf{B} \gets \text{svd}(\PrimalDict, \delta)$}
    \end{subfigure}
	\begin{subfigure}[t]{0.45\linewidth}
        \includegraphics[width=1\linewidth]{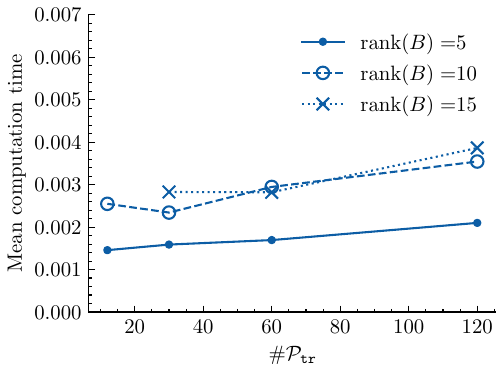}
        \caption{Comp. Time, $\mathbf{B} \gets \text{orth}(\PrimalDict \mathbf{R})$}
    \end{subfigure}
    \caption[Mean number of iterations and computation times for Hertz problem using convex hull exploration.]{Mean number of iterations and computation times for Hertz problem using convex hull exploration with different dimensions of $\mathbf{B}$. Note that computation time excludes time for construction of non-linear operators.}
    \label{fig:convex_hertz_computation_time}
\end{figure}

\subsubsection{The non-convex case of the Ironing problem}\label{sec:nonconvex_ironing}
\begin{figure}[!htb]
    \centering
	\begin{subfigure}[t]{0.4\linewidth}
        \includegraphics[width=1.0\linewidth]{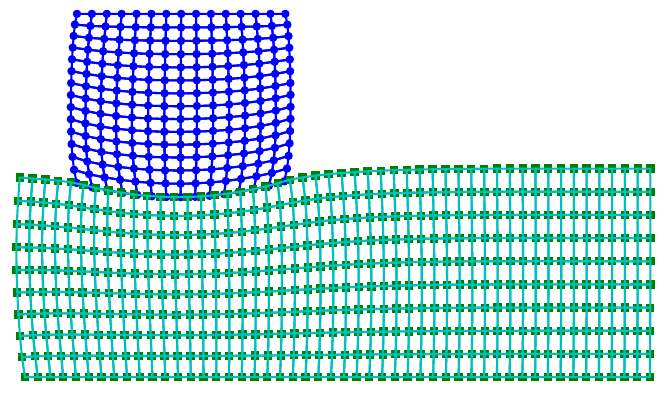}
        \caption{Snapshot A}
    \end{subfigure}
    \quad
	\begin{subfigure}[t]{0.4\linewidth}
        \includegraphics[width=1.0\linewidth]{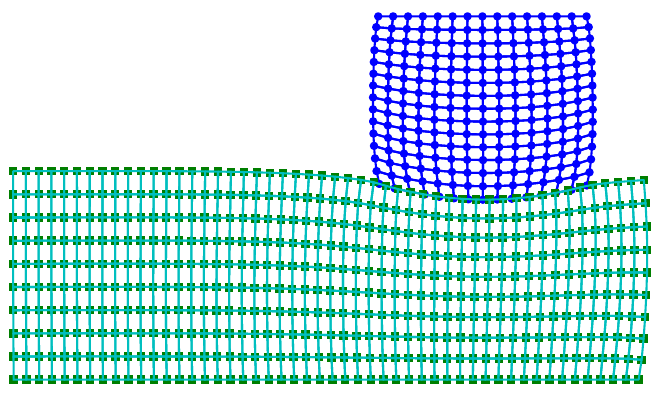}
        \caption{Snapshot B}
    \end{subfigure}
	\begin{subfigure}[t]{0.4\linewidth}
        \centering
        \includegraphics[width=1.0\linewidth]{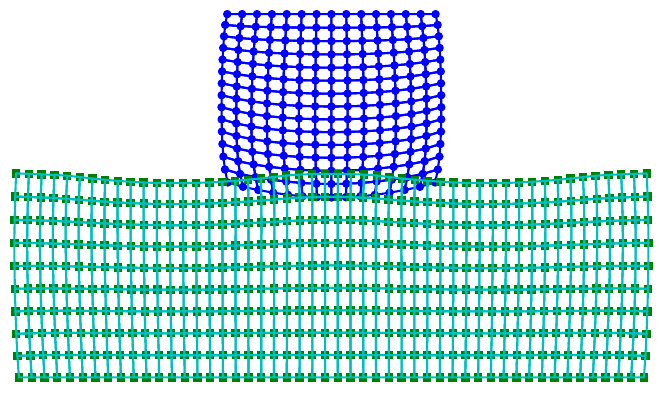}
        \caption{A convex combination of snapshots A and B}
    \end{subfigure}
\caption[Demonstration of the non-convex nature of the ironing problem.]{Demonstration of the non-convex nature of the ironing problem. Two snapshots from training set of the ironing problem and a convex combination of the snapshots that violates the non-penetration condition}
    \label{fig:nonconvex_ironing}
\end{figure}

The convexity test for the feasible set (Test~\ref{test:convex_feasible_set}) is performed on the Ironing problem. Since, the contact zone changes significantly, the test performed using the training set snapshots fails as expected. The CrPen error was found in $\mathcal{O}(10^{-1})$. Therefore, the ironing problem cannot be solved using the current formulation of convex hull approach. The non-convexity can also be demonstrated by a convex combination that violates the non-penetration condition significantly, as shown in \cref{fig:nonconvex_ironing}. For the two snapshots shown in this figure, the slab surface is deformed at different positions. The convex combination where the two displacement field snapshots were equally weighted, generated a deformed configuration with the slab being deformed at the two locations corresponding to each of the snapshots. But no deformation was generated at the position of the iron, causing penetration at that position. In this case, local behaviour of the displacement field leads to the non-convexity of the ironing problem.

\section{Discussions}
The application of dictionary methods was proposed to mitigate the challenges posed by linear inseparability of the contact pressure. As the lack of local information prevents the creation of a robust RB\@, over-complete dictionaries have the potential to preserving local information computed in the offline stage, thereby providing a better chance at resolving the position of contact more effectively. The same is demonstrated by building reasonably accurate reconstructions using dictionary methods.

For the reconstructions of Hertz problem using the greedy active-set algorithm, the reconstruction errors decay consistently when richer training sets are used, and the computational time per iteration seems to be fairly constant, but change slightly with the primal rank.

The case of ironing problem where the contact position changes significantly is also considered. However, the decay of reconstruction errors vs.\ size of dictionary curves were found to be flatter compared to the Hertz problem. This is probably due to stronger dependence of contact position on the parameter of the problem. Though the dictionary size didn't have any significant influence on the computational time per iteration, the number of iterations increased significantly for larger dictionaries in case of ironing problem. 

Comparing the two dictionary methods: the greedy active-set algorithm and the convex hull approximations, the reconstruction errors of the Hertz problem were better using the greedy active-set algorithm (see \cref{fig:hertz_mean_error_valset,fig:convex_hertz_recon_err}). This observation can be possibly attributed to the use of the monolithic dictionary in case of the convex hull approximations; which forced common interpolation coefficients for the displacement as well as contact pressure. In particular, the convex hull approach shows large errors for small training sets, as the convex hull of a small set of snapshots might be a very small subregion in the entire solution space. 

Another possible source of error in convex hull approximations is that the snapshots in the dictionary may not be the best ones to explore the low-dimensional convex subset. On the other hand, the computational time (and number of iterations) for the convex hull approach seems to be more independent of dictionary size, compared to the greedy active-set algorithm (see \cref{fig:hertz_num_iters,fig:hertz_computation_time,fig:convex_hertz_computation_time}), thanks to the elimination of active set selection process.  

In case of the convex hull approach, the reconstruction error seems to be more sensitive to the size of dictionary than the truncation tolerance/rank of the orthogonal basis~$\mathbf{B}$. Also, for the greedy active-set method, the sensitivity to truncation seems to be low between $10^{-10}$ and $10^{-8}$. Moreover, in both approaches, the computational time per iteration is nearly flat when the rank is low. Thus, one can infer that that large dictionaries and moderate truncation tolerances / rank of $\mathbf{B}$ are an optimal choice. 

One could argue that the number of iterations increases significantly with size of dictionary in case of the ironing problem solved using the greedy active-set method (\cref{fig:ironing_hlevel_evolution_iters,fig:ironing_hlevel_evolution_time}), however, this is related to the inherent non-linearity of the contact problem. The number of iterations could possibly be reduced by using smarter initialization of the contact pairs according to the queried parametric value rather than the initializing with reference state. For the convex hull approach, the primary disadvantage is the assumption of convexity of the feasible region, limiting the applicability to convex problems. The non-convex nature of ironing problem meant that convex hull approximation method could not be used in this case.

\section{Perspectives}
\subsection{Efficient construction of constraints}
In the \cref{alg:sparse_greedy}, \cref{step:operators}  can be the most expensive step if full operators are built from scratch and then projected on the reduced bases. The constraint operators $\widehat{\mathbf{C}}$ and $\widehat{\bm{g}}$ show non-linear dependence on the displacement field $\widehat{\bm{u}}$ and therefore must be reconstructed in every iteration. An efficient construction of the constraint operators for low-rank methods is discussed in~\cite{Benaceur2020} using the Empirical Interpolation Method (EIM), where affine decompositions of the distance functions are computed offline. This enables efficient construction of constraint operators in online phase by sampling distance functions at a few points on the contact surface, followed by computing the coefficients of the affine form.

The distance function $k(\mu,\bm{v},\bm{v})$ of~\eqref{eq:distance_functions} when expressed in discrete form (as discussed in~\cite{Benaceur2020}) turns out to be a sparse vector. The \texttt{spy} pattern of this vector can be a strong function of the parameter when contact area changes significantly. This is caused by the dependence of the distance function on the contact pairs established using the current configuration. Hence, it is bound to exhibit some form of linear inseparability just like contact pressure. Therefore, it is possible that linear inseparability issues might also appear in generating an effective affine decomposition.

\subsection{Application of sketching methods to dictionary-based approximation}
As discussed briefly in \cref{sec:dictionary_resolution} and in detail in~\cite{Balabanov2021}, random sketching can be used to project the residual of a dictionary-based approximation problem on a randomly sampled subspace of the dictionary. By controlling the size of the subspace, one can control the cost of solving the dictionary-based approximation problem, relatively independent of the dictionary size. The difficulties in applying these techniques to contact pressure dictionary appear due to the inequality constraints. If random sketching methods are applied, the inequality constraints are linearly combined in random fashion leading to loss of information on violation (penetration) and satisfaction (non-penetration) of constraints. 



\subsection{Smart snapshot selection for Convex hull explorations}
RB methods equipped with error estimators sample the parametric space efficiently. Similarly, it might be possible to extend the same notion to efficiently sampling the convex feasible set. The convex hull approximations might be more robust if the snapshots lie on the edges/vertices of the low-dimensional convex subset presumed in \cref{hyp:lowdim_convex_subset}. For e.g.\ convex hull approximation for the low-dimensional problem in \cref{fig:convex_hypothesis_test}, it is important  to include the snapshot near the corner of the feasible region in the training set for a robust reconstruction. However, deducing the corners in a high-dimensional problem may not be geometrically intuitive as the low-dimensional problem in the figure.

\chapter{Exploring the contact pressure manifold}\label{ch:nlInterp}
In \cref{sec:philosophy_dictionary_approx}, the philosophy of dictionary-based approximation methods for contact problems was discussed. The principal argument in favour of over-complete dictionaries was the lack of low-rank structure due to linear inseparability of the contact pressure field, which has two consequences: (a) a large number of snapshots might be necessary to truly discover the underlying structure of the contact pressure field, and (b) linear model order reduction methods cannot capture local effects of contact pressure in a low-rank subspace. The focus of \cref{ch:sparse} was on efficient handling large dictionaries of contact pressure snapshots in the online phase, but does not address the high offline cost of creating such dictionaries. This chapter focuses on the use of non-linear dimensionality reduction methods as means to improve reconstruction errors without paying the large computational cost of creating over-complete dictionaries.

\section{Linear Subspaces vs.\ Low-dimensional Manifolds}
\label{sec:manifold}
Many of the dimensionality reduction methods work on the principal assumption that the underlying physics occupies a subset of the solution space defined in the full model. At this point, the question arises about the nature of the aforementioned subset of solution space. Low-rank methods assume that this subset is a low-dimensional linear subspace\footnote{usually the term ``linear'' is dropped in unambiguous cases} of the full solution space. The assumption that such a low-dimensional subspace exists implies that the solution of any point in the particular parametric space may be approximated using a linear combination of a limited number of precomputed basis vectors that naturally lie in this subspace. As seen in \cref{ch:literature}, methods like Proper Orthogonal Decomposition (POD) are capable of extracting a basis for the linear subspace that describes the underlying physics.

However, in some problems the parametric solutions occupy a low-dimensional manifold, but not necessarily linear subspace. Consequently, it is not possible to compute a robust low-rank basis, as shown for the contact pressure field in \cref{ch:lowrank}. Nonetheless, there exists a paradigm called non-linear dimensionality reduction that aims at computing the underlying structure by finding the associated non-linear behaviour. A subclass of non-linear dimensionality reduction known as manifold learning is quite extensively used in the analysis and visualisation of high dimensional datasets~\cite{Bregler1995,Ghodsi2006}. The fundamental idea behind manifold learning is similar to low-rank methods and is known as the manifold hypothesis~\cite{Fefferman2016}. Roughly speaking, a manifold can be described as an $n$-dimensional space that is globally non-linear in general but locally resembles a linear vector space. 


The manifold hypothesis assumes that high-dimensional data sourced from practical applications often lie in low-dimensional manifolds. The low-rank hypothesis is a special case of the manifold hypothesis involving the assumption of the low-dimensional manifold being a linear subspace. Manifold learning methods rely on finding the non-linear low-dimensional manifold within the high-dimensional data. A well-known application of manifold learning is the case of a data set containing gray-scale images of an object rotated by various angles. Linear interpolation methods like PCA fail to reduce this type of data, but manifold learning methods can extract the inherent one-dimensional space~\cite{Bregler1995,Ghodsi2006}.

To better approximate the manifold structure of a parametrized problem, partitioning of the parametric domain~\cite{Dihlmann2011} and a local reduced basis methods~\cite{Amsallem2012} have been explored. Interpolation in the tangent space to the Grassmann manifold structure have been used to construct reduced models with a wide parametric range~\cite{Amsallem2008}. 

One of the most commonly applied non-linear dimensionality reduction techniques in the parametrized ROMs is the kernel methods~\cite{Haasdonk2018,Gonzalez2018,Diez2021}. Kernel interpolation methods along with greedy algorithms to create nested training sets to create efficient surrogate models are discussed in~\cite{Haasdonk2018}. The kernel-PCA (kPCA) approach can also be used for the computation of intrinsic parameters that are not explicitly included in the mathematical model~\cite{Gonzalez2018}. An application of kPOD to an advection-diffusion problem with inseparable nature can be found in~\cite{Diez2021}, where the effectiveness of kernel-based method is shown using the higher concentration of variance in the first singular values of kPOD compared to the case of linear POD\@.

A known drawback of kernel-based methods and some other manifold learning techniques is that backward mapping, i.e.\ reconstruction of a high-dimensional entity from low-dimensional data is not as trivial as it is done in linear dimensionality reduction methods such as PCA\@. Usual the reconstruction is an approximation to the local tangent space of the manifold, computed using a linear combination of snapshots with weights computed by employing a minimization of a discrepancy functional criteria~\cite{GarciaGonzalez2020}. An improved approximation based on a quadratic tangent space is discussed in~\cite{Diez2021}.

Hyperbolic systems are often associated with moving discontinuities, which also leads to lack of a linear subspace. An approach involving the transformation of the physical variable has been used in such problems to separate the discontinuous and continuous behaviour~\cite{Bansal2021}. A related approach aimed at generating transformed snapshots thereby aligning the discontinuities has been demonstrated on Burger's equation~\cite{Welper2015}. As it will be seen, the ideas given in this chapter have similar motivations of aligning the contact pressure snapshots, but without transforming the physical variable.

\subsection{Contact pressure as a non-linear manifold}
In \cref{ch:lowrank}, qualitative and quantitative arguments about the linear inseparability of contact pressure have been made. The lack of low-dimensional linear subspace, i.e.\ a low-rank subspace for contact pressure was demonstrated. In this section, we discuss the idea of contact pressure lying in a low-dimensional non-linear manifold. Looking back at \cref{fig:hertz_snaps} and \cref{fig:ironing_snaps}, the snapshots in respective figures have similarities that are evident intuitively. Given these similarities, there is a strong possibility that these snapshots exist in a low-dimensional non-linear manifold.

\begin{figure}[!htb]
    \centering
    \includegraphics{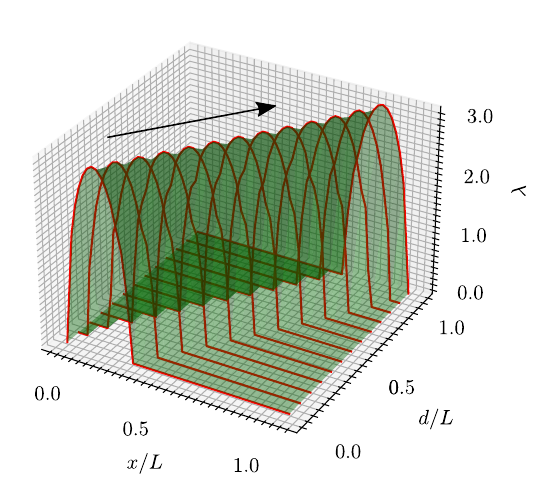}
    \caption[Visualization of contact pressure non-linearity in ironing problem]{Visualization of contact pressure non-linearity in ironing problem. The dimension $x\in \Gamma$ is the coordinate on contact surface, and $d \in \mathcal{P}$ is the coordinate on the parametric domain. The arrow shown in the figure points to a direction in which snapshots are ``aligned'' along which the snapshots are, possibly, linearly separable.}
    \label{fig:ironing_surface}
\end{figure}

This intuitive hypothesis can be supported by \cref{fig:ironing_surface}. Contact pressure snapshots of the ironing problem are shown by the lines as a function of the parameter $d$ (the horizontal displacement of the iron), and the green translucent surface is generated using the surface plot of the same data. One can intuitively see that the contact pressure is a non-linear function of just one parameter, which happens to be also the physical parameter $d$ in this case. The fallacy of low-rank methods is that they attempt to perform a linear interpolation in direction $d$, when in fact the relationship is non-linear. The idea that will be exploited in this chapter is that an appropriate non-linear transformation applied on this surface in a way that allows interpolation along the ``aligned'' direction shown in the figure. This will possibly improve the linear separability and reducibility.

\section{Dimension warping}
The core idea in manifold learning methods (such as kPCA) is that linearly inseparable data can be appropriately mapped into a higher dimension where it is relatively more separable. The mapping to high-dimensional space has to be essentially non-linear, since a linear mapping will not resolve the inseparability issues.

Consider a non-linear mapping from the spatial domain of contact surface $\Gamma$ to a warped domain $\WarpedSpace$. The warped domain $\WarpedSpace$ need not admit a physical interpretation.
\begin{align}
    \begin{aligned}
        \varphi(\mu,\cdot): \Gamma &\mapsto \WarpedSpace \\
        \WarpedCoord &\gets \varphi(\mu,x)
    \end{aligned}
\end{align}
\noindent where $\WarpedCoord$ is a coordinate in $\WarpedSpace$. The notation $\widetilde{x}(\mu,\WarpedCoord)$ is used to indicate $x$ that maps to $\WarpedCoord$.

A qualitative argument can be made regarding an effective linear behaviour in $\WarpedSpace$ space. The mapping $\varphi$ must warp the spatial dimension so that the regions with non-zero contact pressure must be mapped to a common warped region, and the gap (non-contact) regions must be mapped to other regions in warped space. In other words, let the warped space be a union of disjoint sets, $ \WarpedSpace =  \WarpedSpace_{c}  \cup \WarpedSpace_{/c} $, then for every pair $(x, \WarpedCoord)$ the mapping $\varphi$ must satisfy the following:
\begin{align}
    \varphi(\mu,\cdot): \Gamma &\mapsto \WarpedSpace \ \big | \left \lbrace
    \begin{matrix}
        \WarpedCoord \in \WarpedSpace_c \ &\text{if} \ \lambda(\mu,x)>0 \\
        \WarpedCoord \in \WarpedSpace_{/c} \ &\text{if} \ \lambda(\mu,x)=0
    \end{matrix} \right .
    \label{eq:mapping_property}
\end{align}
where $\WarpedSpace_{c}$ and $\WarpedSpace_{/c}$ are the regions where contact and no-contact zones are mapped respectively. Any function $f$ defined on $\Gamma$ has an equivalent $\widetilde{f}$ defined on $\WarpedSpace$ domain and can be related to each other with the following relation:
\begin{align*}
    \widetilde{f}(\WarpedCoord) = f(x)
\end{align*}
\begin{remark}[Preserving the split of contact and gap zones]{}{}
    Extending the above qualitative argument to the case where the contact pressure snapshots have multiple disjoint regions of contact and gap zones, each subregion should be mapped to corresponding subregions in warped space. For e.g., in the ironing problem, each snapshot can be split into three zones: one contact zone and two gap zones. Thus, there should be three disjoint zones in warped space corresponding to each zone of the original space.
\end{remark}

\subsection{Dynamic Time Warping}
\label{sec:dtw}
It is not straightforward to compute the mapping $\varphi(\mu, \cdot)$ and the warped space $\WarpedSpace$ with the desired properties discussed in the previous section. Algorithms that identify regions of similarities between two signals can prove to be useful to generate mappings that improve separability. One such method, known as Dynamic Time Warping (DTW), is an algorithm built to measure similarity between two discrete signals containing similar but displaced features. The algorithm is named so because it was built for measuring similarities in temporal sequences involving speech or motion~\cite{Rabiner1993,Muller2007}. Other methods exist that deals with sequence similarities, such as elastic matching~\cite{Talmi2017}; but they are not studied in this thesis. In this section, we illustrate the application of DTW to generate mappings that satisfy the conditions specified in~\eqref{eq:mapping_property}.\footnote{Though the DTW algorithm was originally intended for temporal data, here it will be used to warp spatial dimension.} Though only unidimensional DTW has been explored here, it is extendable to higher dimensions~\cite{Shokoohi2017}.

The DTW algorithm works on two given signals and computes similarities based on distances between pairs of points. It matches each point on the first signal to at least one point on the second, and vice-versa. A detailed explanation of DTW can be found in~\cite{Rabiner1993} and Python implementations can be found in the \texttt{dtw} module~\cite{Giorgino2009}. The mapping $\varphi(\mu_1,x)$ and $\varphi(\mu_2,x)$ generated by DTW is described by indices of the discretization of the coordinate $x$. This is illustrated using a simplified example. \\

\noindent \textbf{An illustration of DTW}: Consider a set of $n$ discrete points $\bm{x}$ in $[0,1]$ and two discrete vectors $\bm{\lambda}_k \in \mathbb{R}^n, k \in \{1,2\}$. For a simple illustration, the vectors are chosen such that each has a peak (with value $1.0$) at different positions surrounded by vanishing areas (values are given in \cref{tab:dtw_illus_values}). The DTW alignments (computed using~\cite{Giorgino2009}) are given by indices $I_k \in \mathbb{R}^{\widetilde{n}} $. The index at $I_1[m]$ is matched to the index at $I_2[m]$ of respective $\bm{\lambda}$ vectors (in Python notation). Usually $\widetilde{n}\geq n$, where the equality holds if alignment happens to be a one-to-one mapping.

\begin{figure}[!htb]
    \centering
    \includegraphics[width=0.48\linewidth]{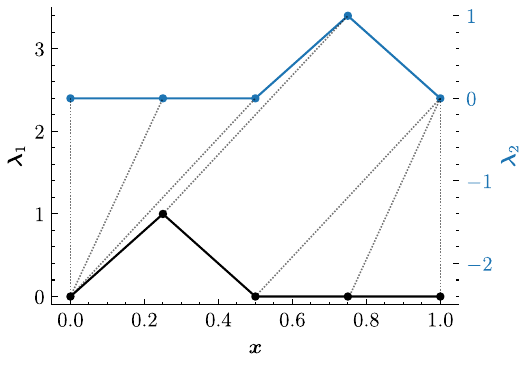}
    \caption[DTW illustration: Vectors $\bm{\lambda}_1$ and $\bm{\lambda}_2$ and their DTW alignments]{DTW illustration: Vectors $\bm{\lambda}_1$ and $\bm{\lambda}_2$ and their DTW alignments (dotted). Note that the vectors are offset in $y$-direction for illustration purposes.}
    \label{fig:DTWillus_lambda1_and_2}
\end{figure}

\begin{figure}[!htb]
    \centering
    \begin{subfigure}[t]{0.48\linewidth}
        \includegraphics[width=1\linewidth]{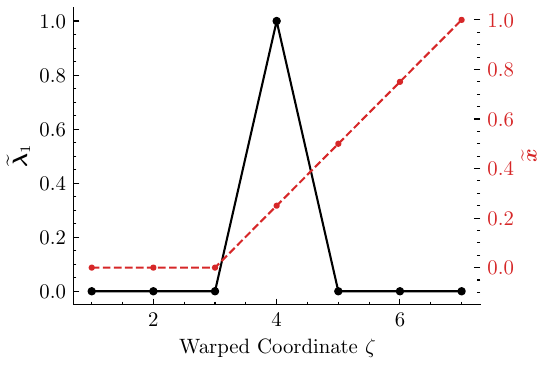}
        \caption{$\widetilde{\bm{x}}_1$ and $\widetilde{\bm{\lambda}}_1$}
    \end{subfigure}
    \begin{subfigure}[t]{0.48\linewidth}
        \includegraphics[width=1\linewidth]{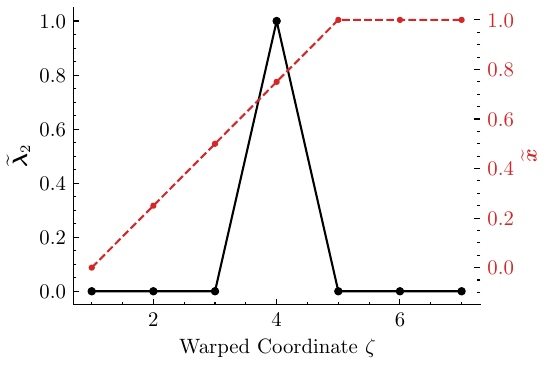}
        \caption{$\widetilde{\bm{x}}_2$ and $\widetilde{\bm{\lambda}}_2$}
    \end{subfigure}
    \caption{DTW illustration: $\widetilde{\bm{\lambda}}$ and $\widetilde{\bm{x}}$ vs. warped coordinate $\WarpedCoord$}
    \label{fig:DTWillus_v1_v2_warped}
\end{figure}

The alignments are given in \cref{tab:dtw_illus_maps} and are shown graphically in \cref{fig:DTWillus_lambda1_and_2}. One can see that the alignments connect not only the peaks of the two vectors, but also the vanishing areas on the either sides as well. The alignments $I_k$ are equivalent to $\varphi(\mu_k, \cdot)$ of~\eqref{eq:dtw_interp}, the mapping to warped coordinate $\WarpedCoord$. The warped snapshots $\widetilde{\bm{\lambda}}_k = \bm{\lambda}_k[ I_k] $ and the warped spatial coordinate $\widetilde{\bm{x}}_k = \bm{x}[ I_k]$ of dimension $\widetilde{n}$ are shown in \cref{fig:DTWillus_v1_v2_warped}. The DTW warped coordinate $\WarpedCoord$ can be simply expressed as set of indices $\lbrace 1, 2 \dots \widetilde{n} \rbrace$.

\begin{table}[!htb]
    \setlength\extrarowheight{2pt}
    \small
    \centering
    \begin{subtable}[t]{0.48\linewidth}
        \centering
        \begin{tabular}[t]{ll}
            \toprule
            Variables  & Input Values  \\
            \midrule
            $\bm{x}$ & $\left[0, 0.25, 0.5, 0.75, 1 \right]$ \\
            $\bm{\lambda}_1$ & $\left[0, 1.0, 0, 0, 0 \right]$ \\
            $\bm{\lambda}_2$ & $\left[0, 0, 0, 1.0, 0 \right]$ \\
            \bottomrule
        \end{tabular}
        \caption{Data}
        \label{tab:dtw_illus_values}
    \end{subtable}
    \begin{subtable}[t]{0.48\linewidth}
        \centering
        \begin{tabular}[t]{ll}
            \toprule
            Alignments  & Observed values  \\
            \midrule
            $I_1$ & $\left[1, 1, 1, 2, 3, 4, 5 \right]$ \\
            $I_2$ & $\left[1, 2, 3, 4, 5, 5, 5 \right]$ \\
            \bottomrule
            & \\
        \end{tabular}
        \caption{DTW alignments}
        \label{tab:dtw_illus_maps}
    \end{subtable}
    \caption{Input variables values and the observed DTW alignments for the illustrative example}
    \label{tab:dtw_illus}
\end{table}

\noindent Though a detailed explanation of the DTW algorithm is beyond the scope of this thesis, the major features of uni-dimensional DTW are listed here:
\begin{enumerate}
    \item The first (and the last) points of the two input vectors are always connected, even if they are connected to other points. In other words, $I_k[1] = 0$ and $I_k\left[\widetilde{n}\right]=n$ for any inputs.
    \item The alignments give indices in a monotonically non-decreasing manner i.e. $I_k[j+1] \geq I_k[j]$, meaning that the alignments shown in \cref{fig:DTWillus_lambda1_and_2} never cross each other.
    \item Roughly speaking, the alignments are generated such that the distance between warped vectors $\widetilde{\bm{\lambda}}_1$ and $\widetilde{\bm{\lambda}}_2$ is minimized. This is evident in the alignments and the warped vectors shown in \cref{fig:DTWillus_lambda1_and_2} and \cref{fig:DTWillus_v1_v2_warped}.
\end{enumerate}

\begin{figure}[!htb]
    \centering
    \begin{subfigure}[t]{0.49\linewidth}
        \includegraphics[width=1\linewidth]{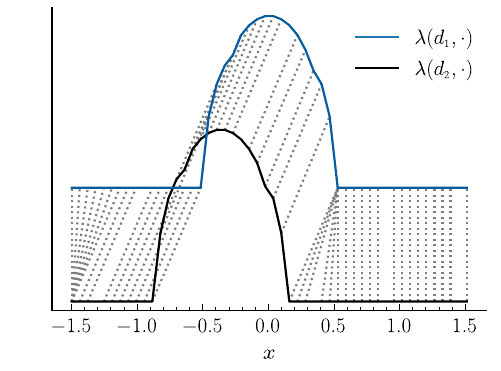}
        \caption{}\label{fig:dtw_twoway}
    \end{subfigure}
    \begin{subfigure}[t]{0.49\linewidth}
        \includegraphics[width=1\linewidth]{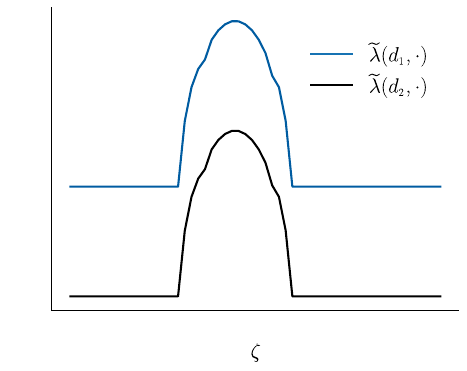}
        \caption{}\label{fig:dtw_warping}
    \end{subfigure}
    \caption[Alignments and warped snapshots computed by DTW between two contact pressure snapshots of ironing problem]{\textbf{(a)} Alignments and \textbf{(b)} warped snapshots computed by DTW between two contact pressure snapshots of ironing problem for parametric values $\mu_1=d_1$ and $\mu_2=d_2$. Notice that alignments in original space are not bijective. $y$-offset between the two curves is given only for illustrative purpose. Values on $\WarpedCoord$ axis are not given as they do not carry a physical significance, and also because DTW computes warped space on indices of the discretization rather than physical dimension.}\label{fig:dtw_twoway_and_warped}
\end{figure}

Now, the application of DTW algorithm to the contact pressure curves of the ironing problem is considered. The alignments for two snapshots corresponding to parametric values $d_1$ and $d_2$ are shown in \cref{fig:dtw_twoway} and the warped contact pressure curves are shown in \cref{fig:dtw_warping}. Notice that the alignments were able to connect respective contact and no-contact zones of both curves. Also, in the warped spatial dimension, respective contact and no-contact zones of both the warped snapshots are the same. In other words, condition~\eqref{eq:mapping_property} is satisfied. Also, looking back at the ``aligned'' direction in \cref{fig:ironing_surface}, we can say that the warping has rendered this direction orthogonal to the warped coordinate $\WarpedCoord$, thereby making linear combinations useful.

\begin{remark}[Non-linear nature of DTW mapping]{}{}
    The mapping generated by DTW is linear w.r.t.\ individual snapshots. However, the mapping is not only different for each snapshot, it depends on the pair of snapshots taken into account. These properties mean that the mapping is essentially non-linear.
\end{remark}

\subsection{Interpolating in warped space}\label{sec:warped_interp}
Although the DTW algorithm is built to compute similarities between signals, we intend to use it as an interpolation tool. Interpolations using DTW have been previously performed to reconstruct missing frames in fast motion detection~\cite{Almog2005}. The objective of performing DTW on pairs of contact pressure snapshots is to use the DTW alignments and warped space for exploring the contact pressure manifold without generating intermediate snapshots directly with FE solver.

\begin{figure}[!htb]
    \centering
    \begin{subfigure}[t]{0.49\linewidth}
        \includegraphics[width=1\linewidth]{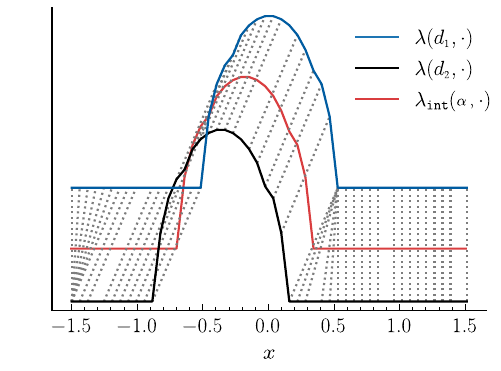}
        \caption{}
        \label{fig:dtw_interp_demo_original}
    \end{subfigure}
    \begin{subfigure}[t]{0.49\linewidth}
        \includegraphics[width=1\linewidth]{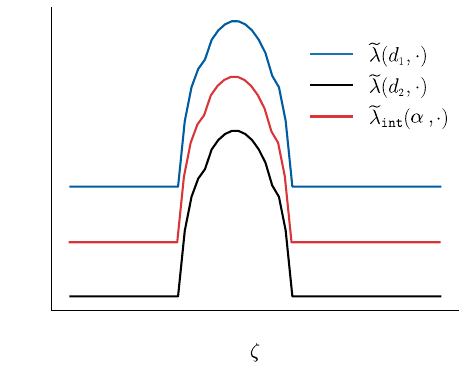}
        \caption{}
        \label{fig:dtw_interp_demo_warped}
    \end{subfigure}
    \caption[Interpolation along DTW alignments and in warped space]{Interpolation along \textbf{(a)} DTW alignments and \textbf{(b)} interpolation in warped space, for parametric values $\mu_1=d_1$ and $\mu_2=d_2$. Notice that unlike in \textbf{(a)}, the interpolation in \textbf{(b)} is a linear operation.}
    \label{fig:dtw_interp_demo}
\end{figure}

\cref{fig:dtw_interp_demo} shows a demonstration of an interpolation in the original and the warped spatial dimension. Though interpolating along the alignments in the original dimension looks intuitive (\cref{fig:dtw_interp_demo_original}), only the interpolation in warped space is a straightforward linear operation (\cref{fig:dtw_interp_demo_warped}). Coefficients $\alpha$ and $1-\alpha$ are used as weights for each snapshot in this operation, and the case for $\alpha=0.5$ is shown in the figure. In principle, the value of the parameter $\mu$ that corresponds to the interpolated snapshot is not known.

As interpolation is performed in the warped dimension, it means that a backward mapping is also needed to transform the interpolated curve back into the original spatial dimension. However, the mapping $\varphi(\mu, \cdot)$  is a function of the parameter and the mapping for the interpolated curve needs to be computed.

\noindent Formally, the interpolations on the two snapshots, and the mappings for parametric values $\mu_1 $ and $\mu_2$ can be written as:
\begin{subequations}
    \begin{align}
        \widetilde{\lambda}_\texttt{int}(\alpha,\WarpedCoord) := \alpha \widetilde{\lambda}(\mu_1,\WarpedCoord) + (1-\alpha)\widetilde{\lambda}(\mu_2,\WarpedCoord) \label{eq:lambda_interp} \\
        \widetilde{x}_\texttt{int}(\alpha,\WarpedCoord) := \alpha \widetilde{x}(\mu_1,\WarpedCoord) + (1-\alpha)\widetilde{x}(\mu_2,\WarpedCoord) \label{eq:x_interp} \\
        \varphi_\texttt{int}(\alpha,x) := \alpha \varphi(\mu_1,x) + (1-\alpha)\varphi(\mu_2,x) \label{eq:phi_interp}
    \end{align}
    \label{eq:dtw_interp}
\end{subequations}
where $\lambda_\texttt{int}$ and $\varphi_\texttt{int}$ are interpolated contact pressure and mapping respectively and $\alpha$ is the linear interpolation parameter between the two snapshots. The interpolated contact pressure in original and warped space can be related as:
\begin{align}
    \lambda_\texttt{int}(\alpha,x) = \lambda_\texttt{int}(\alpha,\varphi_\texttt{int}^{-1}(\alpha,\WarpedCoord)) = \widetilde{\lambda}_\texttt{int}(\alpha,\WarpedCoord) 
    \label{eq:dtw_back_map}
\end{align}
Once DTW alignments are available,~\eqref{eq:dtw_interp} and~\eqref{eq:dtw_back_map} can be easily implemented in a short routine \texttt{dtw\_interp} to perform the interpolation. 

\begin{notation}
    There is a slight abuse of notation in~\eqref{eq:dtw_interp},~\eqref{eq:dtw_back_map} and in interpolated quantities here onwards. The first argument of quantities {\normalfont{} $\circ_{\texttt{int}}(\cdot,\cdot)$} is the interpolation parameter $\alpha$ and not the physical parameter $\mu$. 
\end{notation}

\begin{notation}
    Interpolated quantities {\normalfont{} $\circ_{\texttt{int}}(\cdot,\cdot)$} also depend implicitly on $\mu_1$ and $\mu_2$, which are not shown in the arguments to simplify the notation.
\end{notation}


\texttt{dtw\_interp} is given in \cref{alg:dtw_interp}. It takes two discrete contact pressure snapshots $\bm{\lambda}_k \in \mathbb{R}^{n}, k \in \{1,2\}$ and the corresponding spatial coordinate $\bm{x} \in \mathbb{R}^{n}$ i.e.\ there are $n$ discrete dofs on the contact surface. The DTW package computes alignments and a warped space of dimension $\widetilde{n}$. Indices $I_k$ provide the mapping to the warped space for respective $\bm{\lambda}_k$. The warped vectors $\widetilde{\bm{\lambda}}_k \in \mathbb{R}^{\widetilde{n}}$ can be linearly combined to create the interpolated vector $\widetilde{\bm{\lambda}}_{\texttt{int}}\in \mathbb{R}^{\widetilde{n}}$ and the corresponding spatial coordinate $\widetilde{\bm{x}}_{\texttt{int}}\in \mathbb{R}^{\widetilde{n}}$. The \texttt{interp1d} function in \cref{step:xt_to_x} is borrowed from the Python module \texttt{scipy.interpolate}~\cite{scipy}. This function is used for changing the discretization from $\widetilde{\bm{x}}_{\mathtt{int}}$ to $\bm{x}$ using linear piecewise interpolations. \cref{step:lambda_int,step:x_int} corresponds to~\eqref{eq:lambda_interp} and~\eqref{eq:x_interp} and \cref{step:xt_to_x} corresponds to~\eqref{eq:dtw_back_map}.


\begin{algorithm}[!htb]
    \caption{\texttt{dtw\_interp}}\label{alg:dtw_interp}
	\begin{algorithmic}[1]
        \State{Input: $\bm{x}, \bm{\lambda}_1, \bm{\lambda}_2, \alpha$}
        \State{$I_1, I_2 = \mathtt{dtw}(\bm{\lambda}_1, \bm{\lambda}_2)$} \Comment{$I_1$ and $I_2$ map indices to warped dimension}
        \State{$\widetilde{\bm{\lambda}}_{\texttt{int}} = \alpha \bm{\lambda}_1[I_1] + (1-\alpha) \bm{\lambda}_2[I_2] $ } \label{step:lambda_int} \Comment{$\bm{\lambda}_j[I_j]$ is the warped vector $\widetilde{\bm{\lambda}}_j$}
        \State{$\widetilde{\bm{x}}_{\texttt{int}} = \alpha \bm{x}[I_1] + (1-\alpha) \bm{x}[I_2] $} \label{step:x_int} \Comment{interpolating in warped spatial coordinate}
        \State{$\bm{\lambda}_{\texttt{int}} \gets \texttt{interp1d}(\widetilde{\bm{\lambda}}_{\texttt{int}}, \widetilde{\bm{x}}_{\texttt{int}}) (\bm{x})$} \label{step:xt_to_x} \Comment{linear interpolation from $\widetilde{\bm{x}}$ to $\bm{x}$}
        \State{Output: $\bm{\lambda}_{\texttt{int}}$}
	\end{algorithmic}
\end{algorithm}

\begin{figure}[!htb]
    \centering
    \begin{subfigure}[t]{0.48\linewidth}
        \includegraphics[width=1\linewidth]{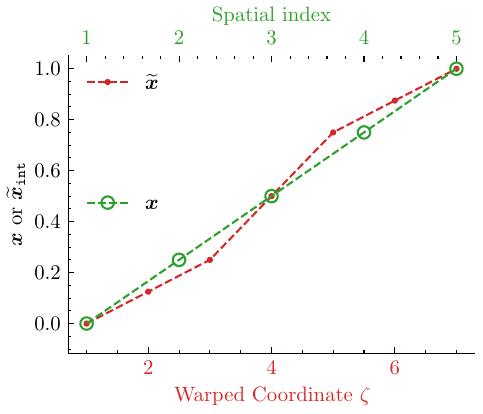}
        \caption{$\bm{x}$ and $\widetilde{\bm{x}}_\texttt{int}$ vs. $\WarpedCoord$}
        \label{fig:DTWillus_interp_space}
    \end{subfigure}
    \begin{subfigure}[t]{0.48\linewidth}
        \includegraphics[width=1\linewidth]{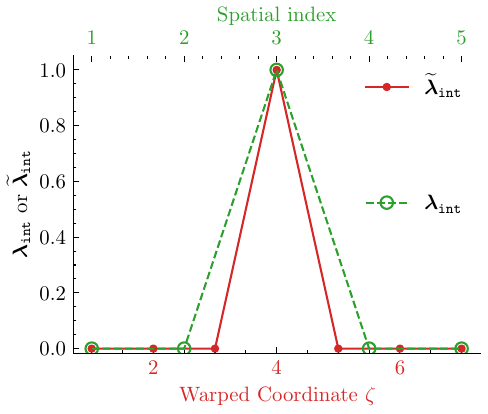}
        \caption{$\bm{\lambda}_\texttt{int}$ and $\widetilde{\bm{\lambda}}_\texttt{int}$ vs. $\WarpedCoord$}
        \label{fig:DTWillus_interp_vector}
    \end{subfigure}
    \caption{Interpolations of the illustrative example for $\alpha=0.5$ in the spatial indices (corresponding to $x$) and the warped space $\WarpedCoord$}
    \label{fig:DTWillus_interp}
\end{figure}

\noindent The application of \texttt{dtw\_interp} is demonstrated on the illustrative example of \cref{tab:dtw_illus} in \cref{sec:dtw}. The interpolated quantities for this example are shown in \cref{fig:DTWillus_interp}. 

\begin{remark}[Interpolation using the alignments $I_k$]{}{}
    In \cref{alg:dtw_interp}, the computation of $\bm{\lambda}_{\texttt{int}}$ from $\widetilde{\bm{\lambda}}_{\texttt{int}}$ can be implemented in an alternate way. The interpolation of the mapping~\eqref{eq:phi_interp} can be computed by interpolating between $I_1$ and $I_2$ to generate a set of interpolated indices. The steps would read as:
\begin{algorithmic}[0]
    \State{$I_{\texttt{int}} = \lfloor(\alpha I_1 + (1-\alpha) I_2) + \dfrac{1}{2}\rfloor$}
    \State{$\bm{\lambda}_{\texttt{int}}[I_{\texttt{int}}] \gets \widetilde{\bm{\lambda}}_{\texttt{int}}$}  \Comment{same as applying inverse mapping $\varphi^{-1}_{\texttt{int}}$}
\end{algorithmic} \vspace{0.5em}
involving a round off operation to the nearest integer, expressed here using the greatest integer function $\lfloor {\cdot} \rfloor$ and the offset of $\frac{1}{2}$. This is equivalent to using the inverse mapping $\varphi^{-1}$ in~\eqref{eq:dtw_back_map}. However, this is approach is not followed in \cref{alg:dtw_interp}, as the round off operation is sensitive to truncation at machine precision level, especially for $\alpha=0.5$, leading to improper interpolation. Instead, $\widetilde{\bm{x}}_{\texttt{int}}$ that corresponds to the warped coordinate $\WarpedCoord$ is used for interpolating back to $\bm{x}$. Also, the resulting $I_{\texttt{int}}$ contains non-unique entries, leading to non-injective inverse map $\varphi^{-1}_{\texttt{int}}$.
\end{remark}

\begin{figure}[!htb]
    \centering
    \begin{subfigure}[t]{0.45\linewidth}
        \includegraphics[width=1\linewidth]{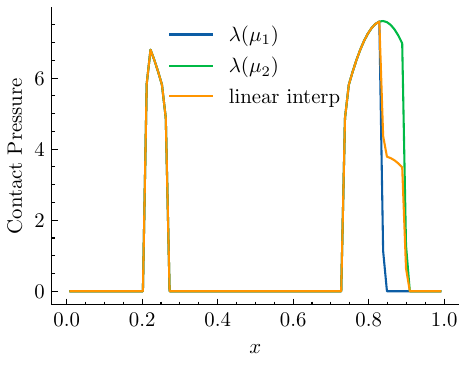}
        \caption{Linear interpolation}
    \end{subfigure}
    \begin{subfigure}[t]{0.45\linewidth}
        \includegraphics[width=1\linewidth]{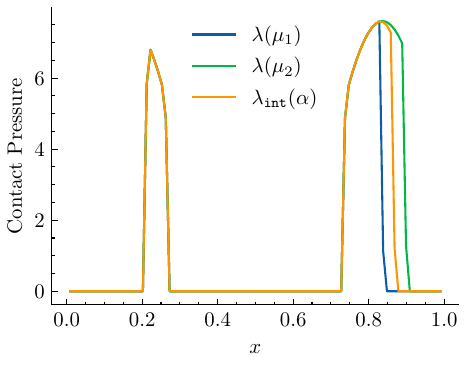}
        \caption{DTW-based interpolation}
    \end{subfigure}
    \caption[Linear and DTW-based interpolations of contact pressure for the elastic rope-obstacle problem]{Linear and DTW-based interpolations of contact pressure for the elastic rope-obstacle problem.} 
    \label{fig:interps_membrane}
\end{figure}

\begin{figure}[!htb]
    \centering
    \begin{subfigure}[t]{0.45\linewidth}
        \includegraphics[width=1\linewidth]{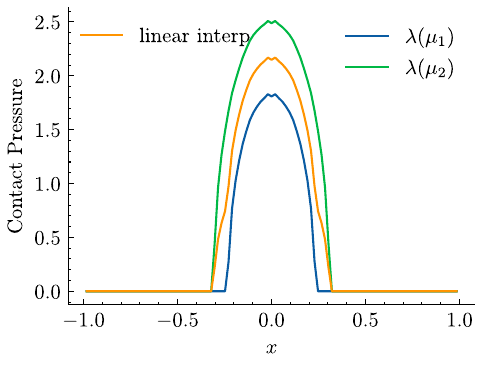}
        \caption{Linear interpolation}
    \end{subfigure}
    \begin{subfigure}[t]{0.45\linewidth}
        \includegraphics[width=1\linewidth]{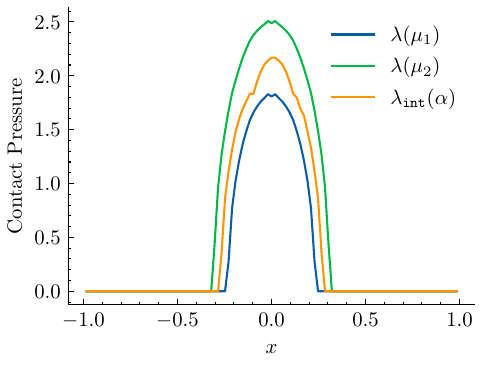}
        \caption{DTW-based interpolation}
    \end{subfigure}
    \begin{subfigure}[t]{0.40\linewidth}
        \includegraphics[width=1\linewidth]{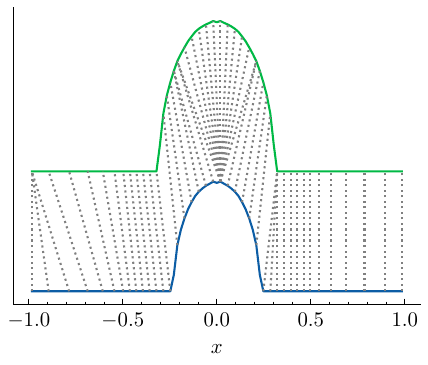}
        \caption{DTW alignments}
    \end{subfigure}
    \caption[Linear and DTW-based interpolations of contact pressure for the Hertz problem]{\textbf{(a)} Linear and \textbf{(b)} DTW-based interpolations of contact pressure for the Hertz problem.  To explain distorted interpolation, DTW alignments is also shown in~\textbf{(c)}.} 
    \label{fig:interps_hertz}
\end{figure}

The quality of interpolation using DTW on the ironing problem is quite good as the snapshots have a similar contact pressure profile, with only the contact zone being translated in the $x$-direction. DTW-based interpolations of the elastic rope-obstacle problem and the Hertz problem are also shown in \cref{fig:interps_membrane,fig:interps_hertz}. In these cases contact zone does not translate but magnifies. DTW is able to interpolate the contact zone correctly, but the interpolated pressure profile has some aberrations near the peak pressure of the Hertz case. This is explained by the concentration of the DTW alignments near the peak i.e., many points from $\lambda(\mu_2)$ are mapped to the peak of $\lambda(\mu_1)$.

\section{DTW-based enrichments for contact mechanics}
As mentioned at the beginning of the chapter, the idea is to use manifold exploration to reduce the offline cost. Since the main cause of linear inseparability is the geometrically local nature of contact pressure field, the ability to scale and translate local effects across the domain of interest seems promising. Also, the displacement field does not have linear separability issues and a robust linear subspace seems to be readily available (\cref{fig:ironing_h_error}) even if parametric space is relatively unexplored in offline stage. In an optimistic scenario, if the warped interpolations can reduce the size of training set to the same level necessary for capturing the displacement field, the bottleneck posed by linear inseparability will be effectively removed.

\begin{figure}[!htb]
    \centering
    \includegraphics{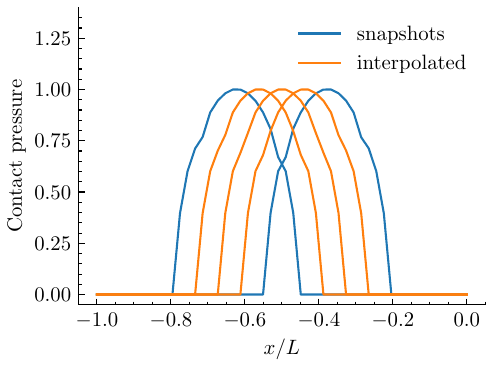}
    \caption[Interpolation of contact pressure snapshots using DTW for the ironing problem]{Interpolation of contact pressure snapshots using DTW for the ironing problem. Interpolated curves for $\alpha \in \{0.25, 0.5, 0.75 \}$ are shown}
    \label{fig:nlInterp_demo}
\end{figure}

We propose the use of DTW as a warped interpolation tool that can generate intermediate contact pressures that haven't been explored in the offline stage. This is an inexpensive way of enriching the contact pressure dictionary; since the computation cost of warped interpolations is negligible compared to generating proper snapshots. \cref{fig:nlInterp_demo} shows a demonstration of such an interpolation for three different values of $\alpha$. It is noticeable that the peak position of contact pressure of the interpolated curves relative to the peaks of the original snapshots seem to have a nearly linear relation with $\alpha$, though this is a particular feature of the ironing problem and may not apply in general.

For an effective enrichment, the interpolation should not be done on every pair of contact pressure curves. Rather, each pair of curves should be reasonably similar in some sense. Similarity can be computed using a measure of distance which can be either the $\mathcal{L}^{2}$ norm of the difference or the distance measure defined by the DTW algorithm. Snapshots that can be considered similar according to the distance measure will be referred to as \emph{neighbours} here onwards. For each pair of neighbours, the DTW alignments are computed and interpolations are performed with a set of $M$ parameter values $\lbrace \alpha_m \rbrace_{m=1}^M$. This process is shown in \cref{alg:dtw_enrich}, where neighbours of currently active vectors are computed using the distances between vectors of the dual dictionary in a matrix $\texttt{Dist}$. The entry $\texttt{Dist}[i,j]$ has the distance between $\DualDict[:,i]$ and $\DualDict[:,j]$, and hence, \texttt{Dist} is a symmetric matrix with main diagonal containing only zeros.

\begin{algorithm}[!htb]
    \caption{\texttt{adapt\_dictionary}}\label{alg:dtw_enrich}
	\begin{algorithmic}[1]
        \State{Input: $\DualDict,\ \mathcal{I},\ \lbrace \alpha_m \rbrace_{m=1}^{M}$} \Comment{$\alpha_m \in \ ]0,1[$}
        \State{Parameters: \texttt{Dist}, $N$}
        \State{Initialize: $\DualDict^{\texttt{adapt}} = \DualDict[:,\mathcal{I}],\ i=0,\ j=0$} \Comment{$\mathcal{I}$ is the current active set}
        \For{$i \in \mathcal{I}$}
            \State{$\texttt{nbrs} \gets $ indices of the least $N$ positive entries of $\texttt{Dist}[i,:]$}
            \State{Initialize $\DualDict^{\texttt{int}} = \{\cdot\}$} 
            \For{$j \in \texttt{nbrs}$} \Comment{\texttt{nbrs} contains neighbouring indices}
                \If{pair $(i,j)$ is not repeated} \Comment{$(i,j)$ is an unordered pair}
                \State{$\displaystyle \DualDict^{\texttt{int}} \gets \DualDict^{\texttt{int}} \cup \left[\ \bigcup_{m=1}^{M} \texttt{dtw\_interp}(\DualDict[:,j], \DualDict[:,i], \alpha_m) \ \right] $}  \Comment{\cref{alg:dtw_interp}}
                \EndIf{}
            \EndFor{}
            \State{$\DualDict^{\texttt{adapt}} \gets \DualDict^{\texttt{adapt}} \cup \DualDict[:,\texttt{nbrs}]  \cup \DualDict^{\texttt{int}} $}
        \EndFor{}
        \State{Output: $\DualDict^{\texttt{adapt}}$}
	\end{algorithmic}
\end{algorithm}

\begin{figure}[!htb]
    \centering
    \begin{subfigure}[t]{0.49\linewidth}
        \includegraphics[width=1\linewidth]{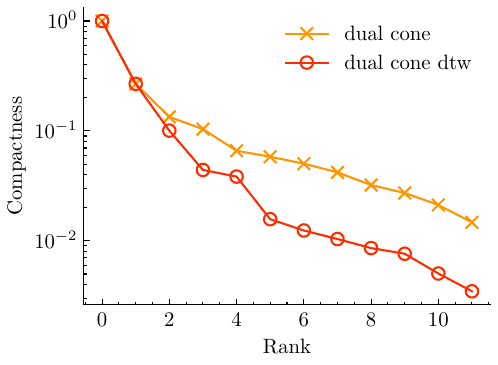}
        \caption{Hertz problem}
        \label{fig:compactness_dtw_hertz}
    \end{subfigure}
    \begin{subfigure}[t]{0.49\linewidth}
        \includegraphics[width=1\linewidth]{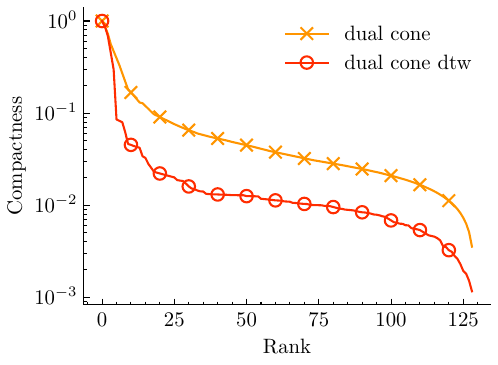}
        \caption{Ironing problem}
        \label{fig:compactness_dtw_ironing}
    \end{subfigure}
    \caption[Comparison of the compactness for original and warped snapshots]{Comparison of the compactness for original and warped snapshots. Note that the curve labelled ``dual cone'' are the same curves from \cref{fig:hertz_rb_rank_error,fig:irorning_rank_error_fine_mesh}, respectively. Curves labelled ``dual cone dtw'' show the compactness of the enriched basis.}\label{fig:compactness_dtw}
\end{figure}

The effectiveness of DTW-based enrichments can be shown using the compactness metric. Recall that the ``dual cone'' compactness metric defined in \cref{sec:metrics} can be used for as a measure of linear separability. The comparison of this metric for the original and enriched CPG dual basis is shown in \cref{fig:compactness_dtw}. The compactness of the enriched basis was computed by performing $19$ equispaced DTW interpolations (i.e. with $\alpha_m = 0.05 m$ and $M=19$ as input for \cref{alg:dtw_enrich}) between each pair of neighbours in the original dual basis. The compactness of the enriched basis is considerably lower than the original, indicating that DTW enrichment has introduced the information that was missing in the cone of the original basis.

\begin{remark}[Indefiniteness of the DTW kernel]{}{}
    A more appropriate way of testing the separability of given data in warped space is by observing the decay of singular values (equivalent to the dual orth compactness in \cref{ch:lowrank}). However, computation of singular values requires that the DTW mapping generates a symmetric positive definite (SPD) Gramian matrix $G_{ij} = \innerproduct{\widetilde{\lambda}_i}{\widetilde{\lambda}_j}$, where each entry of $G$ contains an inner product of a pair of warped snapshots. Such a comparison would have been possible if SPD kernels are used for computation of the Gramian matrix, as done in~\cite{Diez2021} for an advection-diffusion problem. However, instances of Gramian matrix computed using DTW kernel being indefinite are known~\cite{Lei2007}. The Gramian matrix computed for this case also resulted in negative eigenvalues, hence the dual orth compactness cannot be compared.
\end{remark}

\begin{figure}[!htb]
    \centering
    \includegraphics{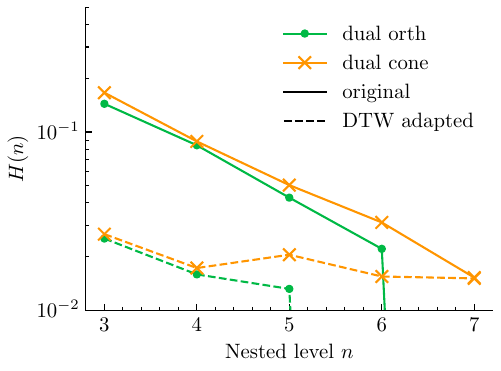}
    \caption[Projection error of nested level $n+1$ snapshots on dictionary of level $n$]{Projection error of nested level $n+1$ snapshots on dictionary of level $n$. Note that solid dual orth and dual cone lines are similar to the ones in \cref{fig:ironing_h_error}. The dual orth points not seen in the figure drop to the order of $10^{-14}$ (nearly machine precision)}
    \label{fig:nlInterp_h_error}
\end{figure}

The richness of $\DualDict^{\texttt{adapt}}$ w.r.t $\DualDict$ can also be compared using the nested errors defined for the ironing problem (see \cref{sec:ironing}, \cref{fig:ironing_h_error}). Let us recall that the nested error~\eqref{eq:h_error} was defined as the projection error of $(n+1)$-th level basis on $n$-th level basis. The comparison of nested errors with and without DTW interpolation is shown in \cref{fig:nlInterp_h_error}. For each level $n$, the dictionary neighbours are interpolated using \cref{alg:dtw_enrich} with $\mathcal{I}=\{1,2,3,\dots \# \mathcal{P}_{\texttt{tr}}\}$, $N=2$, $M=1$ and $\alpha_1=0.5$. The ironing problem has a peculiar property where the position of contact pressure peak is more or less same as the parameter $d$. Therefore, the chosen $N, M$ and $\alpha_1$ form a near best case scenario of the approximation because it suits the distribution of nested training set points i.e.\ $(n+1)$-th level contain mid-points of the points in $n$-th level. Thus, a considerable decrease in the nested errors is seen for smaller training sets. Naturally, the gains in dual orth error due to DTW decreases as the training set becomes larger, as the snapshots in large training sets are closer to the target solutions.

However, in a general scenario, the value of an optimum $\alpha$ may not be trivial. Therefore, to arrive at an appropriate value of $\alpha$, it is necessary to explore multiple DTW-based interpolations corresponding to different $\alpha$'s. This is the reason a set of $M$ values $\lbrace\alpha_m\rbrace$ is supplied to \cref{alg:dtw_enrich}. The enrichment process can be strategised in two different ways:
\begin{itemize}
    \item \textbf{Number of interpolation points} ($M$): For a given set of snapshots generated in offline phase, warped interpolation is performed between every pair of neighbouring snapshots  with a certain resolution of in $]0,1[$. For a uniformly distributed points, $\alpha_m = m/(M+1)$. The interpolations in \cref{fig:nlInterp_demo} is such a case with $M=3$.
    \item \textbf{Number of enrichment levels} ($L$): It is possible to continue the enrichment process further by creating multiple levels of interpolated vectors. Suppose, the original dictionary contains two vectors $\bm{\lambda}_1$ and $\bm{\lambda}_2$ and interpolation is performed with $M=1$ and $\alpha_1=0.5$. Thus, an interpolated vector, let us call it $\bm{\lambda}_{1.5}$, can be computed using DTW interpolation in the first level. In a second level, further interpolation can be performed on the pairs $\lbrace\bm{\lambda}_1, \bm{\lambda}_{1.5}\rbrace$ and $\lbrace\bm{\lambda}_{1.5}, \bm{\lambda}_2\rbrace$ to generate $\bm{\lambda}_{1.25}$ and $\bm{\lambda}_{1.75}$.
\end{itemize}

\noindent Further, these enrichments can be done during either the online or offline stage, depending on the desired distribution of the computational cost. If done in the offline stage, warped interpolations between every possible pair of neighbours need to be computed. On the other hand, in the online stage, interpolations can computed on a smaller set of vectors, i.e. the active dictionary vectors and their neighbours. In this thesis, the enrichment in the online stage is studied, as the focus is on the enhancement of the reconstruction accuracy using DTW and the study of computational cost  is reserved for future study. However, a preliminary comment can be made in this aspect, that the DTW-based enrichment is not a computationally intensive task. As the \cref{alg:sparse_greedy} generates a small active set, the number of interpolations performed is not large. The resulting algorithm is given in \cref{alg:dictionary_dtw}, with a parameter $L$ denoting the number of enrichment levels. In the numerical results that follow, the parameters used are $L=5$ and $M=4$ with the interpolation points $\lbrace 0.2, 0.4, 0.6, 0.8 \rbrace_{m=1}^M$.

\begin{algorithm}[!htb]
	\caption{DTW powered greedy active-set method}\label{alg:dictionary_dtw}
	\begin{algorithmic}[1]
		\State{Input : Queried value of parameter $\mu$}
        \State{Parameters : $\displaystyle \lbrace \alpha_m \rbrace_{m=1}^M$, $L$}
		\State{Given: Primal basis $\PrimalRB$ and dual Dictionary $\DualDict$}
		\Statex{\hspace{3em} Reduced operators $\mathbf{\PrimalRB}^T \mathbf{K} \mathbf{\PrimalRB} \text{ and } \mathbf{\PrimalRB}^T \bm{f}$}
        \State{Initialize: $l = 0, \ \ \mathcal{I} = \{\cdot\}, \ \ \bm{\lambda} = \bm{0}$}
        \While{$\bm{\lambda}$ not converge and $l \leq L$}
        \State{Compute $\widehat{\bm{u}}, \widehat{\bm{\lambda}}, \mathcal{I}$ using current dictionary $\DualDict$} \Comment{\crefrange{step:loop_start}{step:loop_end} of \cref{alg:sparse_greedy}}
            \State{Reconstruct $\bm{\lambda} = \DualDict \widehat{\bm{\lambda}}$}
            \State{$\DualDict \gets \mathtt{adapt\_dictionary}(\DualDict, \mathcal{I}, \lbrace \alpha_m \rbrace)$} \Comment\ \cref{alg:dtw_enrich}
            \State{$\mathcal{I} = \lbrace \cdot \rbrace$} \Comment{Reset active set after adapting $\DualDict$}
            \State{$l \gets l+1$}
        \EndWhile{}
        \State{Reconstruct $\bm{u} = \PrimalRB \widehat{\bm{u}}$}
        \State{Output: $\bm{u}, \bm{\lambda}$}
	\end{algorithmic}
\end{algorithm}

\subsection{Numerical results}
\begin{figure}[!htb]
    \centering
    \begin{subfigure}[t]{0.98\linewidth}
        \includegraphics[width=1\linewidth]{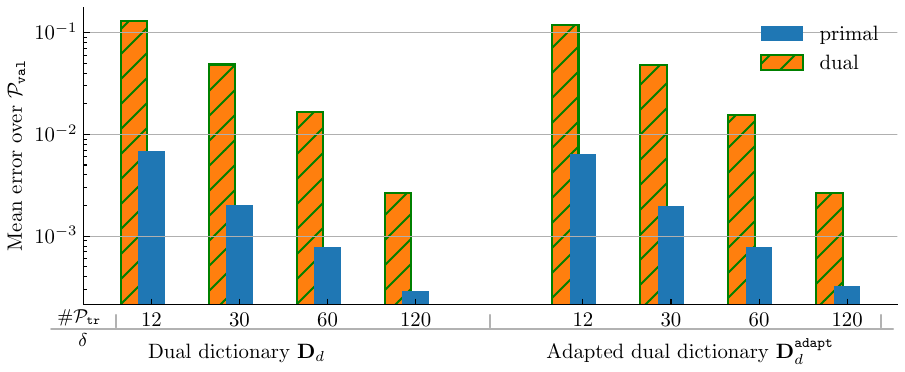}
        \caption{Hertz problem}
        \label{fig:hlevel_evolution_dtw_hertz}
    \end{subfigure}
    \begin{subfigure}[t]{0.98\linewidth}
        \includegraphics[width=1\linewidth]{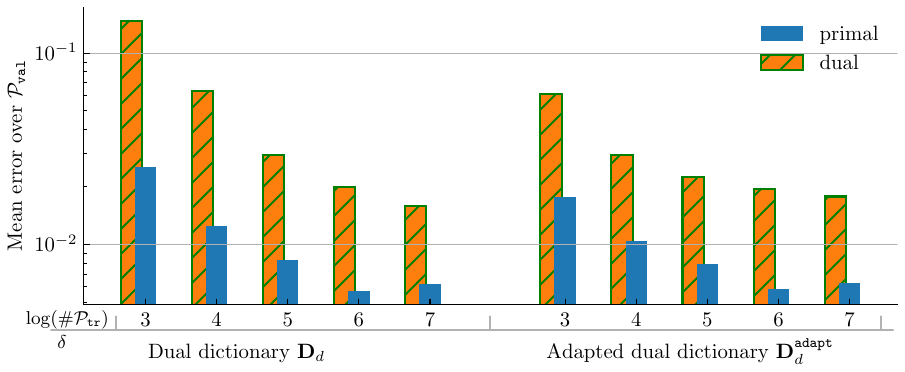}
        \caption{Ironing problem}
        \label{fig:hlevel_evolution_dtw_ironing}
    \end{subfigure}
    \caption[Comparison of reconstruction error for the Hertz and ironing problem using the original dual dictionary and DTW-adapted dual dictionary]{Comparison of reconstruction error for the Hertz and ironing problem using the original dual dictionary $\DualDict$ and DTW-adapted dual dictionary $\DualDict^{\texttt{adapt}}$. Reconstructions performed using \cref{alg:sparse_greedy} and \cref{alg:dictionary_dtw} for respective dictionaries. DTW enrichments seem to be more impactful in the ironing problem than Hertz problem. Primal and dual errors are computed using $\mathcal{H}^1$ and $\mathcal{L}^2$ norm, respectively.}
    \label{fig:hlevel_evolution_dtw}
\end{figure}

For the Hertz problem and the ironing problem training sets defined in \cref{ch:sparse}, reconstruction errors with enriched dual dictionaries using \cref{alg:dictionary_dtw} are shown in \cref{fig:hlevel_evolution_dtw}. For the ironing problem (\cref{fig:hlevel_evolution_dtw_ironing}), the DTW-based interpolation improves the quality of reconstruction for smaller training sets. The improvement asymptotically dampens out as the training set size increases. This is natural as snapshots come quite close in large training sets and DTW-based interpolations do not generate any useful information. Moreover, the improvement in reconstruction is seen more in the dual field rather than the primal. This is an expected observation as any improvement in primal reconstruction due to dual dictionary enrichments is caused indirectly through the mixed system~\eqref{eq:sparse_greedy_system}. 

On the other hand, there is no improvement in the reconstruction of Hertz problem (\cref{fig:hlevel_evolution_dtw_hertz}). One likely reason could have been that the contact pressure snapshots are relatively close, unlike the ironing problem, and therefore the DTW-based interpolations are unable to offer any significant new information. However, this contradicts the observations in \cref{fig:compactness_dtw_hertz} where improvements in dual cone compactness is observed. A more likely reason could be, that to improve in reconstruction of Hertz problem, the  primal dictionary also needs enrichment rather than just the dual dictionary. However, artificial enrichment could not be applied to the primal dictionary (as it needs a multidimensional DTW implementation), and the dual dictionary enrichments alone are ineffective at improving the errors. This reasoning is more consistent with observation that reconstruction error improves only with increase in training set size and not with DTW-based enrichments. Another possible source of error is the distorted DTW-based interpolation of the contact pressure snapshots, as shown in \cref{fig:interps_hertz}. Though, in reality, the distortion is not as pronounced as shown in the figure, since the interpolation done in \cref{alg:dictionary_dtw} is carried out between snapshots that are closer.

\section{Conclusions}
The DTW enrichments proved to be useful in improving contact pressure reconstructions for a problem with higher variations in contact area (ironing problem), especially when the dictionary is not rich enough to resolve the local effects. The improvements obtained using DTW diminish with increasing training set size, as the local information related to contact pressure start appearing in the dictionary directly. For problems with less variations in contact pressure (Hertz problem), DTW does not provide any improvements in reconstruction errors. 

Overall, it is fair to say that non-linear interpolations can prove useful to train models for problems with large variability in contact positions, especially if the complexity of high-fidelity problem makes it difficult to compute a large training set. This is a strong motivator to test other non-linear reduction techniques such as kPOD on such problems. However, whether it is trivial to enforce inequality constraints in non-linear frameworks or not is yet to be explored.

\section{Persepctives}
\subsection{The 3D case}
Extension of DTW based interpolation to 3D problems (with 2D contact surfaces) may not be straightforward. Though the DTW algorithm was originally built for a 1D signals, but a 2D version has been developed~\cite{Shokoohi2017}. However, the 2D DTW is limited to structured meshes, thereby limiting its application. Even if the geometry of the contact surfaces allow discretization using a structured mesh, it may not be the optimal mesh for generating the snapshots using finite elements. It may however, be possible to compute the FE snapshots using a generic mesh and then interpolating the solution on a structured mesh during a post-processing, thereby making it suitable for DTW-based interpolations. 

\subsection{Efficient evaluation of non-linear terms}
To solve the contact problem using the dictionary method, it is necessary to compute the operators $\widehat{\mathbf{C}}$ and $\widehat{\bm{g}}$, which are dependent on the dual dictionary $\DualDict$. However, if DTW based interpolation is used, the dictionary $\DualDict$ is dynamic. This can create additional difficulties in devising affine decompositions of these operators. If EIM approach from~\cite{Benaceur2020} is used, computation of the independent terms $\widehat{\mathbf{C}}^p$ and $\widehat{\bm{g}^q}$ in~\eqref{eq:dist_funcs_eim_operators} becomes dependent on the dual basis/dictionary. An efficient way of handling such situations is necessary. A straightforward but less efficient way might be to build an affine decomposition of high-fidelity operators $\mathbf{C}$ and $\bm{g}$, and update the reduced operators whenever the operators or the dictionary is updated.

\subsection{Non-linear interpolations of primal dictionary}
As seen in the case of Hertz problem, the dual dictionary enrichment is not effective at improving reconstruction errors. Also, in case of the ironing problem, mainly the dual reconstruction errors show improvement than the primal one. This observation opens up the possibility of enriching the primal dictionary using non-linear methods. This will also be useful to generate enriched monolithic dictionary for convex hull based approximation. 

\chapter*{Conclusions and Perspectives}
\addcontentsline{toc}{chapter}{Conclusions and Perspectives}
\markboth{Conclusions and Perspectives}{}
\markright{Conclusions and Perspectives}{}

\section*{Conclusions}
In \cref{ch:lowrank}, the robustness of low-rank hypothesis in the context of contact mechanics was analyzed. The premise of this analysis was based on the local nature of the contact pressure field which leads to an inherent lack of linear separability. As the underlying physics exhibited a high intrinsic dimensionality, reconstructions performed using ROMs had a moderate error, as expected. The robustness of contact pressure reduced bases computed using SVD and CPG methods were shown to have a slower decay of error (i.e.\ higher values of compactness) than the displacement field. The lack of low-rank structure had also negatively influenced the generalization ability of the reduced basis.

The linear inseparability of the contact pressure makes it harder to find a robust low-dimensional structure in the framework of linear dimensionality reduction. In \cref{ch:sparse}, it is shown that a large number of snapshots can provide some improvements to the approximate solution subspace that is exploited by reduced models. This is achieved by using large dictionaries of snapshots, and the reconstructed solution is restricted to a sparse linear combination of these snapshots. It is seen that dictionary-based approximations are quite useful for contact problems with large changes in contact position like the ironing problem, as local effects are quite significant in these problems. 

Two dictionary-based approximation strategies were devised. The first is based on the active-set approach implemented greedily, and the second is based on exploring the convex-hull of the snapshots to reconstruct the solution under the condition that the feasible (non-penetrative) region of the displacement space is convex. The computational time for these methods seems to be primarily dependent on the number of iterations for convergence, which may or may not change significantly with the size of the dictionary depending on the contact problem. In terms of reconstruction accuracy, the greedy active-set approach with separate dictionaries for displacement and contact pressure outperforms the convex-hull approximations with monolithic dictionaries. Moreover, the convex-hull approximations could not be applied to non-convex problems such as the ironing problem.

Non-linear dimensionality reduction aspects are particularly interesting for a problem lacking low-rank structure and are explored in \cref{ch:nlInterp}. As the computational costs of generating a large dictionary of snapshots can be very expensive, non-linear methods are used to reduce the need for a large number of contact pressure snapshots. Dynamic Time Warping (DTW) based non-linear interpolations were used to synthesize snapshots that were not present in the training set, thereby enriching the dictionary without paying the large cost of full simulations. This is found to be beneficial in cases of large variations in contact position.

\section*{Perspectives}
\begin{enumerate}
    \item The problems arising from linear inseparability of contact pressure have been addressed using over-complete dictionaries and DTW-based enrichments in this thesis. However, distance functions that are used to compute the non-penetration conditions (i.e.\ the inequality constraints) are likely to have similar local effects and consequent linear inseparability issues. Also, the current methods of computing affine decompositions of distance functions (like EIM) are based on low-rank approximations requiring linear separability. Therefore, similar strategies might be required for the efficient construction of the reduced inequality constraint operators. 
    \item Random sketching methods were applied only to the convex hull approach as the method was free of inequality constraints. It would be interesting if the random sketching could be applied to other dictionary-based strategies as well, including the greedy active-set method.
    \item The dictionary method based on the convex hull approach is currently built for monolithic dictionaries, which forces the reduced coefficients to be the same for both primal and dual variables. An approach with two uncoupled dictionaries must be explored for potential improvements in reconstruction quality. Also, an efficient sampling strategy for the offline phase is necessary, where snapshots placed at corners and edges of the presumed low-dimensional convex subset are sampled with priority.
    \item As DTW-based non-linear interpolation was found useful in particular cases with large variations of contact positions, application of other non-linear methods like kPCA/kPOD should be explored. These methods might be more relevant for three-dimensional contact problems with two-dimensional contact surfaces, as multidimensional DTW is limited to structured meshes.
    \item Currently, non-linear interpolation is only applied to contact pressure. It would be interesting to see if non-linear interpolation of displacement snapshots provides any advantages. DTW-based enrichment of displacement snapshots, along with the contact pressure, might be beneficial for creating monolithic dictionaries which need the same number of primal and dual snapshots.
    \item When the dual dictionary is updated using DTW-based enrichments, a fresh computation of the reduced inequality constraint operators is required. This will not be tenable if the DTW-based enrichments are to be deployed in a real-life application of the contact ROM\@. A strategy for efficient computation of these operators for an updated dictionary is necessary.
\end{enumerate}

\appendix
\crefalias{chapter}{appendix}

\chapter{Spurious penetrations}
\label{ch:tau_delta}
\markboth{\nameref{ch:tau_delta}}{}
\markright{\nameref{ch:tau_delta}}{}

\begin{figure}[!htb]
	\centering
    \includegraphics[width=0.5\linewidth]{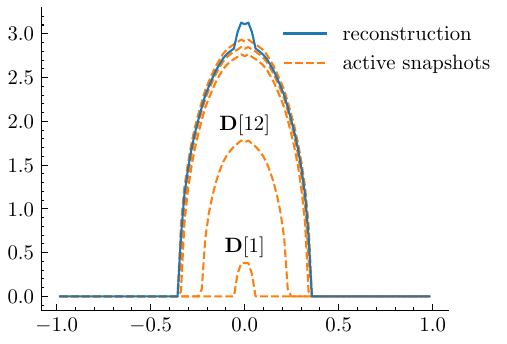}
    \caption[Reconstruction instance  with spurious snapshot selection by the greedy-active set algorithm if non-penetration condition is applied ``hardly'']
    {Reconstruction instance of Hertz problem with $d=0.25$ using $\# \mathcal{P}_{\texttt{tr}} = 30$ and primal basis truncated at $\delta=10^{-6}$, where few dictionary snapshots were spuriously selected by the greedy active-set algorithm if non-penetration condition is applied ``hardly''. The falsely selected snapshots are indicated by their index in the dual dictionary. The snapshot which contributes to most of the error is $\DualDict[1]$, corresponds to loading parameter $d=0.01$}
	\label{fig:hertz_reconstruction_spurious}
\end{figure}

In the greedy active-set method, the snapshots of the dual dictionary are selected based on the violations of the non-penetration constraints. The dictionary vector with the highest correlation with the current state of penetration is chosen. However, instead of zero penetration, a small value of penetration $\tau$ is allowed (see~\eqref{eq:sparse_greedy_kkt_nonpen} and~\eqref{eq:sparse_enrich_concept}). This is done to avoid spurious selection of dictionary elements, which happens due to a truncated primal basis. Such a reconstruction example of Hertz problem for parameter value $d=0.25$ is shown in \cref{fig:hertz_reconstruction_spurious}, where a spurious peak in the contact pressure is evident. In this case, a dual dictionary of size 30 and primal basis truncated at $10^{-6}$ are used. The algorithm selected the appropriate snapshots from the dual dictionary, but also selected some unnecessary snapshots: the first ($\DualDict[1]$) and twelfth ($\DualDict[12]$) dictionary vectors. Importantly, the $\DualDict[1]$ vector contributes to large error in reconstruction as evident in the figure. The false selection of this snapshot is linked to the truncation of primal basis, and therefore, can be avoided by setting the value of $\tau$ same as $\delta$.

\backmatter{}

\cleardoublepage{}
\phantomsection{}
\addcontentsline{toc}{chapter}{Bibliography}
\bibliography{bibfile}

\includepdf[pages={2}]{coverpages/cover.pdf}

\end{document}